\documentclass[preprint,a4paper]{article}
\usepackage{amsfonts}
\usepackage{amssymb}
\usepackage{amsmath}
\usepackage{amsthm}
\usepackage{euscript}
\usepackage{graphics}
\usepackage{multicol}

\newtheorem{theorem}{Theorem}[section]
\newtheorem*{theorem*}{Theorem}
\newtheorem{lemma}[theorem]{Lemma}
\newtheorem{proposition}[theorem]{Proposition}
\newtheorem{conjecture}[theorem]{Conjecture}
\newtheorem{corollary}[theorem]{Corollary}

\begin{document}

\title{The $q$--Gelfand--Tsetlin graph, Gibbs measures and $q$--Toeplitz matrices }

\author{Vadim Gorin\thanks{Institute for Information Transmission Problems, Bolshoy Karetny 19, Moscow 127994,
Russia. e-mail:vadicgor@gmail.com}}

\date{}

\maketitle

\begin{abstract}
The problem of the description of finite factor representations of the infinite-dimensional
unitary group investigated by Voiculescu in 1976, is equivalent to the description of all totally
positive Toeplitz matrices. Vershik-Kerov showed that this problem is also equivalent to the
description of the simplex of central (i.e.\ possessing a certain Gibbs property) measures on
paths in the Gelfand--Tsetlin graph. We study a quantum version of the latter problem. We
introduce a notion of a $q$--centrality and describe the simplex of all $q$--central measures on
paths in the Gelfand--Tsetlin graph. Conjecturally, $q$--central measures are related to
representations of the quantized universal enveloping algebra
$U_\epsilon(\mathfrak{gl}_{\infty})$. We also define a class of $q$--Toeplitz matrices and show
that extreme $q$--central measures correspond to $q$--Toeplitz matrices with non-negative minors.
Finally, our results can be viewed as a classification theorem for certain Gibbs measures on
rhombus tilings of the halfplane.

We use a class of $q$--interpolation polynomials related to Schur functions. One of the key
ingredients of our proofs is the binomial formula for these polynomials proved by Okounkov.
\end{abstract}


\tableofcontents

\section{Introduction}
\subsection{Preface}
 The infinite dimensional unitary group $U(\infty)$ is the union of the unitary groups $U(N)$
naturally embedded one into another. The study of \emph{characters} of $U(\infty)$, i.e.\
normalized positive-definite central continuous functions on the group, was initiated by Voiculescu
in 1976. He studied \emph{finite factor representations} of $U(\infty)$. These representations are
in bijection with \emph{extreme characters} of $U(\infty)$, which are extreme points of the convex
set of all characters of $U(\infty)$. In the paper \cite{Vo} Voiculescu gave a list of extreme
characters of $U(\infty)$. He also conjectured and partially proved that this list was complete.

There is a correspondence between extreme characters of $U(\infty)$ and \emph{totally positive}
Toeplitz matrices. Such matrices were studied much earlier by Shoenberg and his collaborators in
the context of classical analysis. A few years after the paper \cite{Vo}, Boyer \cite{Bo} and
Vershik-Kerov \cite{Vk} independently pointed out that the completeness of the Voiculescu's list of
characters follows from the Edrei's result \cite{Ed} on the classification of all totally positive
Toeplitz matrices.

In the same paper \cite{Vk} Vershik and Kerov suggested a new approach to the above problem. Their
method is based on the approximation of the extreme characters of $U(\infty)$ by the normalized
characters of the irreducible representations of the finite-dimensional unitary groups $U(N)$. As a
result, classification of extreme characters of $U(\infty)$ is restated in purely combinatorial
terms as the problem of the description of the \emph{boundary of the Gelfand--Tsetlin graph}.
Okounkov and Olshanski in their paper \cite{OkOlsh} gave the detailed proof and further generalized
the classification theorem for the characters of $U(\infty)$ using the Vershik-Kerov approach.

The aim of the present paper is to introduce and study a $q$--deformation of the notion of a
character of $U(\infty)$. We start from the Vershik-Kerov formulation and define a $q$--deformed
version of the Gelfand--Tsetlin graph. Our main result is the complete description of the boundary
of the $q$--Gelfand--Tsetlin graph. (See Section \ref{SubSection_Main_res} for the details.)

There are several ways to interpret our results.  We may go back to the characters and introduce
their $q$--analogues which agree with our deformation of the Gelfand--Tsetlin graph. From this
point of view, the problem that we solve in the present paper is the characterization of all
possible limits of rational Schur functions normalized in a certain (depending on $q$) way as the
number of variables grows to infinity. (See Section \ref{subsection_Intro_approximation} for the
details.)

Investigating a $q$--deformation of totally positive Toeplitz matrices, we arrive at a notion of a
$q$--Toeplitz matrix. Every point of the boundary of the $q$--Gelfand--Tsetlin graph corresponds to
a $q$--Toeplitz matrix with minors satisfying some non-negativity condition. (See Section
\ref{Subsection_intro_qToeplitz} for the details.)

Paths in the Gelfand--Tsetlin graph are in bijection with tilings of the halfplane with rhombuses
of 3 types. Through this correspondence our results turn into the classification theorem for
certain Gibbs measures on the tilings of the halfplane. (See Section 2 for the details.)

There are strong reasons to believe that the $q$--Gelfand--Tsetlin graph is related to the
representation theory of the quantized enveloping algebra $U_\epsilon(\mathfrak{gl}_\infty)$.
However we do not address this issue in the present paper.

\subsection{Statement of the main result}
\label{SubSection_Main_res}

\emph{The Gelfand-Tsetlin graph $\mathbb{GT}$} is a graded graph consisting of levels
$\mathbb{GT}_N$, $N=1,2,\dots$. The vertices of $\mathbb{GT}_N$ are $N$--tuples
$\lambda_1\ge\dots\ge\lambda_N$ of integers. Following the Weyl's book \cite{W} we call these
$N$--tuples \emph{signatures}. We join two signatures $\lambda\in\mathbb{GT}_N$ and
$\mu\in\mathbb{GT}_{N+1}$ by an edge and write $\lambda\prec\mu$ if and only if
$$
 \mu_1\ge\lambda_1\ge\mu_2\ge\dots\ge\lambda_N\ge\mu_{N+1}.
$$

Let $\mathcal T$ denote the set of all infinite paths in $\mathbb{GT}$, i.e.\
$\tau=(\tau(1),\tau(2),\dots)$ belongs to $\mathcal T$ if and only if $\tau(N)\in\mathbb{GT}_N$ and
$\tau(N)$ is connected with $\tau(N+1)$ by an edge in $\mathbb{GT}$ for every $N$.

For any probability measure $P$ on $\mathcal T$ let $P_N(\cdot)$ denote its projection on $\mathbb{
GT}_N$, in other words for any $\lambda\in\mathbb{GT}_N$
$$
 P_N(\lambda)=P(\tau\in\mathcal T\mid \tau(N)=\lambda).
$$

For any finite path $\phi=(\phi(1)\prec\phi(2)\prec\dots\prec\phi(N))$, $\phi(k)\in\mathbb{GT}_k$
let $C_\phi$ denote the corresponding cylinder set in $\mathcal T$, i.e.\
$$
 C_{\phi}=\{\tau\in\mathcal T:\, \tau(1)=\phi(1),\dots,\tau(N)=\phi(N)\}.
$$
 Let us introduce the \emph{weight} of $\phi$:
$$
 w(\phi)=q^{|\phi(1)|+|\phi(2)|+\dots+|\phi(N-1)|},
$$
where $|\lambda|=\lambda_1+\dots+\lambda_k$ for $\lambda\in\mathbb{GT}_k$. In Section
\ref{section_probabilistic_setup} we show that these weights have a simple combinatorial
interpretation: given a path $\phi$ one constructs a 3D--body and $w(\phi)$ equals $q^{vol}$, where
$vol$ is the volume of this body. For any $\lambda\in\mathbb{GT}_N$ let
$$
 {\rm Dim}_q(\lambda)=\sum_{\phi(1)\prec\dots\prec\phi(N)|\phi(N)=\lambda} w(\phi).
$$

\emph{A $q$--central} measure on $\mathcal T$ is a probability measure $P$ on $\mathcal T$
satisfying
$$
 P(C_\phi)=P_N(\phi(N))\frac{w(\phi)}{{\rm Dim}_q(\phi(N))}
$$
for any finite path $\phi=(\phi(1)\prec\phi(2)\prec\dots\prec\phi(N))$, $\phi(k)\in\mathbb{GT}_k$.
Put it otherwise, if we consider a family of finite paths $\phi$ ending at the same signature
$\phi(N)$, then the probability of $\phi$ is proportional to the weight $w(\phi)$.

Let $\Omega_q$ denote the set of all $q$--central measures on $\mathcal T$. Clearly $\Omega_q$ is a
convex set. \emph{The minimal boundary of $q$--Gelfand--Tsetlin graph} is the set ${\rm
Ex}(\Omega_q)$ of all extreme points of the set $\Omega_q$. The main result of the present paper is
the description of the minimal boundary of $q$--Gelfand--Tsetlin graph. We need to do some
preparations before stating the theorem.

We say that signature $\lambda\in\mathbb{GT}_N$ is positive if all coordinates of $\lambda$ are
non-negative, i.e.\  $\lambda_N\ge 0$. Let $\mathbb{GT}_N^+\subset\mathbb{GT}_N$ denote the subset
of all positive signatures of size $N$.

Every extreme $q$--central measure $P$ is uniquely defined by the array of numbers $P_N(\lambda)$.
In the case studied by Voiculescu (i.e.\  $q=1$) these numbers appear as the coefficients of the
expansions of certain functions in the basis of rational Schur functions $s_\lambda$. It turns out
that in the $q$--deformed case the formulas have the simplest form if instead of rational Schur
functions we use other symmetric polynomials. It is convenient to use \emph{$q$--interpolation
Schur polynomials} $s^*_\lambda(x_1,\dots,x_N;q)$ defined for every $\lambda\in\mathbb{GT}_N^+$.
These polynomials are a particular case of factorial Schur functions (see \cite{Mac2}). They are
also a particular case ($q=t$) of interpolation Macdonald polynomials (see \cite{Kn}, \cite{S},
\cite{Ok_shifted_Macdonald}). In one dimensional case we have
$$
 s^*_k(x;q)=(x-1)(x-q)\cdots(x-q^{k-1}).
$$

For any probability measure $P_N$ on $\mathbb{GT}_N^+$ we consider the following generating
function:
\begin{equation}
\label{eq_characters_q_ver2} {\mathcal S^*}(x_1,\dots,x_N;P_N)=\sum_{\lambda\in\mathbb{GT}_N^+}
P_N(\lambda)
\frac{s_{\lambda}^*(q^{N-1}x_1,\dots,q^{N-1}x_N;q^{-1})}{s_\lambda^*(0,\dots,0;,q^{-1})}
\end{equation}

Define $\mathcal N$ to be the set of all non-decreasing sequences of integers:
$$
 \mathcal N=\{\nu_1\le\nu_2\le\nu_3\le\dots\}\subset \mathbb Z^{\infty}.
$$

Finally, let  $A_\ell$ be the isomorphism of the Gelfand--Tsetlin graph taking a signature
$\lambda=(\lambda_1,\dots,\lambda_N)\in\mathbb{GT}_N$ to
$A_\ell(\lambda)=(\lambda_1+\ell,\dots,\lambda_N+\ell)$.

Now we are ready to state the main theorem of the present paper.
\begin{theorem}
\label{theorem_Main} Let $0<q<1$. We have:
\begin{enumerate}
\item ${\rm Ex}(\Omega_q)$ is parameterized by points of $\mathcal N$ with $\nu\in\mathcal N$
corresponding to the extreme $q$--central measure $\mathcal E^\nu$.
\item If $\nu_1\ge 0$ then
$$
 {\mathcal S^*}(x_1,\dots,x_N;\mathcal E^\nu_N)=H^\nu(x_1)\cdots H^\nu(x_N),
$$
where
$$
 H^\nu(t)=\dfrac{\prod_{i=0}^{\infty}(1-q^it)}{\prod_{j=1}^{\infty}(1-q^{\nu_j+j-1}t)}
$$
and $\mathcal E^\nu_N$ is the projection of $\mathcal E^\nu$ on $\mathbb{GT}_N$.
\item For general $\nu$ the measure $\mathcal E^\nu$ is the image of the measure $\mathcal E^{\nu'}$ with
$\nu'=(0,\nu_2-\nu_1,\nu_3-\nu_1,\dots)$ under isomorphism $A_{\nu_1}$ of the graph $\mathbb{GT}$.
\item $\Omega_q$ is a simplex meaning that every $P\in\Omega_q$ is a unique average of measures
$\mathcal E^\nu$. Put it otherwise, for every $P\in\Omega_q$ there exist a unique probability
measure $Q$ on $\mathcal N$ such that
$$
 P=\int_{\mathcal N} \mathcal E^{\nu} dQ.
$$
\end{enumerate}
\end{theorem}

There is a natural involution of the graph $\mathbb{GT}$ mapping signature
$\lambda_1\ge\lambda_2\ge\dots \ge\lambda_N$ to $-\lambda_N\ge-\lambda_{N-1}\ge\dots\ge
-\lambda_1$. Is is easy to see that this map provides a bijection between $q$--central measures and
$q^{-1}$--central measures. Thus, Theorem \ref{theorem_Main} also provides a description of
$q$--central measures with $q>1$.

\smallskip

The following proposition describes the support of measure $\mathcal E^\nu$ and explains the
geometric meaning of the parameter $\nu$.

\begin{proposition}
\label{Proposition_description_of_regular_point}
 Let $\tau$ be a random element of $\mathcal T$ distributed according to $\mathcal E^\nu$, then
 almost surely the last coordinates of $\tau(N)$ tend to $\{\nu_j\}$, i.e.\  for every $j$
$$
 \lim_{N\to\infty} \tau(N)_{N+1-j} = \nu_j.
$$
\end{proposition}

The proof of Theorem \ref{theorem_Main} relies on the so-called ergodic method, first proposed by
Vershik in \cite{V_ergodic}. This method was used for the graph $\mathbb{GT}$ in \cite{Vk} and for
the generalizations of $\mathbb{GT}$ in \cite{OkOlsh}. Essentially the same method of the
identification of the boundary of a graph was proposed by Diaconis and Freedman \cite{DF} in the
context of \emph{partial exchangeability}.

For the $q$--Gelfand--Tsetlin graph the ergodic method shows that every extreme $q$--central
measure on $\mathcal T$ is, in some precise sense, a limit of extreme measures on $\mathbb{GT}_N$.
Passing from the measures to their generating functions we arrive at the problem of finding all
possible limits of rational Schur functions normalized in a certain way as the number of variables
grows to infinity. (In Section \ref{subsection_Intro_approximation} we discuss this point in more
details.) It seems natural to describe extreme $q$--central measures in terms of the normalized
rational Schur functions too. However, if we try to do that, then the formulas turn out to be quite
ugly. In the contrast, if we use the polynomials $s^*$ instead of rational Schur functions, then we
arrive at simple multiplicative formulas presented in Theorem \ref{theorem_Main}.

$q$--Interpolation Schur polynomials and their properties play a crucial rule in our proofs. One of
the key ingredients is the \emph{binomial formula}, which relates $q$--interpolation Schur
polynomials with ordinary Schur polynomials and is a particular case of the binomial formula for
interpolation Macdonald polynomials proved by Okounkov in \cite{Ok_bynomial}.

\subsection{Representations of $U(\infty)$ and the origin of the Gelfand--Tsetlin graph}
 In this section we explain how the study of characters of $U(\infty)$ leads to the
Gelfand--Tsetlin graph and compare Theorem \ref{theorem_Main} with the solution of the problem
studied by Voiculescu.

There are two different approaches to the representation theory of $U(\infty)$ which lead to a
reasonable class of representations (i.e.\  finite factor representations, see e.g. \cite{Th}, or
spherical representations of Gelfand pair $(U(\infty)\times U(\infty),U(\infty))$, see
\cite{Olsh_repesentations_U_short}, \cite{Olsh_repesentations_U_long}). In both approaches
representations are in correspondence with \emph{characters}, which are central (i.e.\  constant on
conjugacy classes) positive definite continuous functions on $U(\infty)$ taking value $1$ on the
unit element of the group.

\emph{The Gelfand-Tsetlin graph $\mathbb{GT}$} is a convenient tool for studying characters of
$U(\infty)$. Recall that irreducible representations of $U(N)$ are parameterized by their dominant
weights, which are $N$--tuples $\lambda_1\ge\dots\ge\lambda_N$ of integers. (See e.g.\ \cite{Zh}.)
Thus, the vertices of $\mathbb{GT}_N$ symbolize irreducible representations of $U(N)$. Furthermore,
the edges of $\mathbb{GT}$ encode inclusion relations between irreducible representations of
$U(N+1)$ and $U(N)$.

Given any character $\chi$ of the group $U(\infty)$ we can expand its restriction to the subgroup
$U(N)$ into a convex combination of normalized characters (in the conventional sense) of
irreducible representations of $U(N)$:
\begin{equation}
\label{eq_expansion_character} \chi\rule[-2.3mm]{.4pt}{4mm}_{\,U(N)}
=\sum_{\lambda\in\mathbb{GT}_N} P_N(\lambda) \frac{\chi^{\lambda}(\cdot)}{\chi^\lambda(e)}.
\end{equation}

It is known (see e.g. \cite{Zh}) that the characters of irreducible representations of unitary
groups are rational Schur functions, i.e.\
$$
 \frac{\chi^{\lambda}(U)}{\chi^\lambda(e)}=
 \frac{s_{\lambda}(u_1,\dots,u_N)}{s_\lambda(1,\dots,1)},
$$
where $u_1,\dots,u_N$ are eigenvalues of $U$.

The coefficients $P_N(\lambda)$ determine a probability distribution on the set $\mathbb{GT}_N$.
In this way we get a bijection between characters $\chi$ and sequences $\{P_N\}_{N=1}^{\infty}$ of
probability distributions. These sequences are called $\emph{coherent systems}$ because $P_N$ and
$P_{N+1}$ satisfy a certain relation.

Any coherent system $\{P_N\}$ defines a measure $P$ on $\mathcal T$ in the following way. For a
finite path $\phi=(\phi(1)\prec\dots\prec\phi(N))$ and corresponding cylinder set $C_\phi$ set
$$
P(C_\phi)=\frac{P_N(\phi(N))}{{\rm Dim}(\phi(N))},
$$
where ${\rm Dim}(\lambda)$ is the number of paths $\tau(1)\prec\dots\prec\tau(N)$ such that
$\tau(N)=\lambda$. The coherency relations between $P_N$ and $P_{N+1}$ imply that $P$ is
well-defined. Measure $P$ is \emph{a central measure} in the sense that $P(C_\phi)$ depends only
on $\phi(N)$. In other words, projection of $P$ on the set of all finite paths
$\phi(1),\dots\phi(N)$ ending at a fixed $\phi(N)=\lambda\in\mathbb{GT}_N$ is uniform. Note that
$P_N$ is a projection of measure $P$ on $\mathbb{GT}_N$. Also note that a $q$--central measure
becomes a central measure if we set $q=1$.

Extreme points of the convex set of characters of $U(\infty)$ (in other words, extreme characters)
correspond to \emph{irreducible spherical representations} of pair $(U(\infty)\times
U(\infty),U(\infty))$ (or, again, finite factor representations). Extreme points of the convex set
of all central measures on $\mathcal T$ correspond to extreme characters via the above bijections.

\begin{theorem*}[\cite{Vo},\cite{Bo},\cite{Vk},\cite{OkOlsh}]
 Extreme characters of $U(\infty)$ are parameterized by the points $\omega$ of the
infinite-dimensional domain
$$
 \Omega\subset{\mathbb R}^{4\infty+2}={\mathbb R}^\infty\times {\mathbb R}^\infty \times {\mathbb
 R}^\infty \times {\mathbb R}^\infty \times {\mathbb R} \times {\mathbb R},
$$
where $\Omega$ is the set of sextuples
$$
\omega=(\alpha^+,\alpha^-,\beta^+,\beta^-;\delta^+,\delta^-)
$$
such that
$$
 \alpha^\pm=(\alpha_1^\pm\ge\alpha_2^\pm\ge\dots\ge 0)\in {\mathbb R}^\infty, \quad
 \beta^\pm=(\beta_1^\pm\ge\beta_2^\pm\ge\dots\ge 0)\in {\mathbb R}^\infty,
$$
$$
 \sum\limits_{i=1}^\infty(\alpha_i^\pm+\beta_i^\pm)\le\delta^\pm,\quad \beta_1^+ +\beta_1^-\le 1.
$$
The corresponding extreme character is given by the formula
\begin{equation}
\label{eq_characters_classical} \chi^{(\omega)}(U)=\prod\limits_{u\in {\rm Spectrum}(U)}
e^{\gamma^+(u-1)+\gamma^-(u^{-1}-1)}
\prod_{i=1}^{\infty}\frac{1+\beta_i^+(u-1)}{1-\alpha^+_i(u-1)}
\frac{1+\beta_i^-(u^{-1}-1)}{1-\alpha^-_i(u^{-1}-1)}.
\end{equation}
\end{theorem*}

The remarkable feature shared by the above theorem and Theorem \ref{theorem_Main} is the appearance
of the multiplicative functions. There is an independent representation-theoretic argument proving
that any extreme character of $U(\infty)$ is multiplicative (i.e.\  that it should be a product of
the values of some function over the eigenvalues of the element of $U(\infty)$). However, the
author knows no conceptual reason for the appearance of the multiplicativity in Theorem
\ref{theorem_Main} (independent of the classification theorem itself), and it would be interesting
to find one.

It is very natural to ask whether extreme central measures on $\mathcal T$ can be obtained as
$q\to 1$ limits of extreme $q$--central measures on $\mathcal T$. The answer is ``Yes'' for
extreme central measures such that $\alpha_i^+=\alpha_i^-=0$. For example, if we chose $\nu(q)$
such that
$$H^{\nu(q)}(t)=(1-q^{x_1(q)}t) (1-q^{x_2(q)}t) \cdots (1-q^{x_k(q)}t)$$
and $q^{x_i(q)}\to \beta^+_i$ as $q\to 1$, then the measures $\mathcal E^{\nu(q)}$ tend as $q\to
1$ to the extreme central measure on $\mathcal T$ parameterized by $\beta_1^+,\dots,\beta_k^+$.
Using similar simple arguments one can obtain extreme central measures on $\mathcal T$ with
arbitrary parameters $\gamma^{\pm}$, $\{\beta_i^{\pm}\}$. As for the general case, the author was
not able to point out a sequence of extreme $q$--central measures (with $q$ tending to $1$) that
converges to an extreme central measure with at least one nonzero coordinate $\alpha_i^\pm$.

\subsection{Approximation of characters and limits of symmetric polynomials}

\label{subsection_Intro_approximation}

Vershik and Kerov proved in \cite{Vk} that every extreme character of $U(\infty)$ is a limit of
normalized characters of irreducible representations of $U(N)$, i.e.\  of the functions
\begin{equation}
\label{eq_x6} \frac{s_{\lambda}(u_1,\dots,u_N)}{s_\lambda(1,\dots,1)}.
\end{equation}
  In other words, they reduced the classification problem to the following question: What are the
  possible limits of symmetric polynomials \eqref{eq_x6} as the number of variables grows to
  infinity? The answers for the similar questions for more general polynomials were obtained by
  Okounkov and Olshanski. See \cite[Theorem 1.1]{OkOlsh} and \cite[Theorem 1.4]{OkOlsh_BC}.

The results of the present paper can be also interpreted as an answer to a certain asymptotic
problem for symmetric polynomials as the number of variables grows to infinity.

Let $N(i)$ be an increasing sequence of positive integers and $\lambda(i)\in\mathbb{GT}_{N(i)}$. We
call the sequence of signatures $\lambda(i)$ \emph{regular} if for any $k$ the sequence of
functions
\begin{equation}
\frac{s_{\lambda(i)}(x_1,\dots,x_k,q^{-k},q^{-k-1},\dots,q^{1-N(i)})}
{s_{\lambda(i)}(1,q^{-1},\dots,q^{1-N(i)})} \label{eq_x14}
\end{equation}
 converges uniformly on the set $\{(x_1,\dots,x_k)\in
\mathbb C^k \mid |x_i|=q^{1-i}\}$.

The normalization $s_{\lambda}(1,q^{-1},\dots,q^{1-k})$, that we use, is well-known. See e.g.\
\cite[Examples in Section 3, Chapter 1]{Mac}.

Let $\Lambda$ be the algebra of symmetric functions, i.e.\ a projective limit of the graded
algebras $Sym(N)$ of symmetric polynomials in $N$ variables (see e.g.\ \cite{Mac}).

\begin{theorem} Let $0<q<1$. We have:
\label{theorem_limit_of_polynomials} \begin{enumerate}
\item
 Sequence of signatures $\lambda(i)$ is regular if and only if the last coordinates of $\lambda(i)$
stabilize, i.e.\  if there is a nondecreasing sequence of integers $\nu=\{\nu_j\}$ such that for
any $j>0$
\begin{equation}
\label{eq_stabilization}
 \lim_{i\to\infty}\lambda(i)_{N(i)+1-j}=\nu_j.
\end{equation}
\item Suppose that \eqref{eq_stabilization} holds and let  $Q_k^\nu(x_1,\dots,x_k)$ denote the
limit of the functions \eqref{eq_x14} as $i\to\infty$. Then $Q_k^\nu(x_1,\dots,x_k)$ can be
extended to an analytic function in $(\mathbb C^*)^k$.
\item There is a correspondence between limit functions $Q_k^\nu$ and measures $\mathcal E^\nu$ of Theorem \ref{theorem_Main}. More
precisely
$$
 Q_k^\nu=\sum_{\mu\in\mathbb{GT}_k} \mathcal E^\nu_k(\mu) \frac{s_{\mu}(x_1,\dots,x_k)}
{s_{\mu}(1,q^{-1},\dots,q^{1-k})},
$$
where $\mathcal E^\nu_k$ is the projection of $\mathcal E^\nu$ on $\mathbb{GT}_N$.
\item If $\nu$ and $\nu'$ are such that $\nu_i=\nu'_i+\ell$ for every $i$, then
$$
 Q_k^{\nu}=\frac{x_1^\ell\dots,x_k^\ell}{\left(1\cdot q^{-1}\cdots q^{1-k}\right)^\ell}Q_k^{\nu'}.
$$
\item If $\nu_1\ge 0$, then
\begin{equation}
\label{eq_xx15}
 Q_k^{\nu}=\sum_{\mu\in\mathbb{GT}_k^+}(-1)^{|\mu|}q^{n(\mu)-n(\mu')}{\rm
Spec}_{\nu}(s_\lambda){s^*_{\mu}(x_1,\dots,x_k;q)},
\end{equation}
where the series converges everywhere in $\mathbb C^k$,
$$
n(\mu)=\sum_{i=1}^k (i-1)\mu_i
$$
 and ${\rm Spec}_{\nu}$ is a specialization
of algebra of symmetric function $\Lambda$ (put it otherwise, homomorphism from $\Lambda$ to
$\mathbb C$) with $H$--generating function
$$
\sum_{j=0}^{\infty}{\rm Spec}_{\nu}(h_j) t^j= \frac{\prod\limits_{i\ge
0}(1-q^it)}{\prod\limits_{j=1}^{\infty}(1-q^{\nu_j+j-1}t)}.
$$
(Here $h_j$ stays for the complete symmetric function of degree $j$.)
\end{enumerate}
\end{theorem}

{\bf Remark 1.} Note that in \eqref{eq_xx15} we use functions $s^*_\mu(\cdot; q)$ while in the
definition of generating functions \eqref{eq_characters_q_ver2} we use $s^*_\mu(\cdot; q^{-1})$.

{\bf Remark 2.} The author knows no general explicit formulas for the functions $Q_k^\nu$. This is
the main reason why we do not use these functions as a description of measures $\mathcal E^\nu$.
Instead we use ${\mathcal S^*}(x_1,\dots,x_k; \mathcal E^{\nu}_k)$, for which a simple
multiplicative formula is provided in Theorem \ref{theorem_Main}

{\bf Remark 3.} In a special case $\nu_n=\nu_{n+1}=\dots=a$ the functions $Q_k^\nu$ admit a simple
interpretation. Let us consider only the case $a=0$. Note that if $0<q<1$, then for any
$f\in\Lambda$ sequence $f_n=f(1,q,\dots,q^k)$ converges as $n\to\infty$, and a specialization map
$$
 f(y_1,y_2,\dots)\mapsto f(1,q,q^2,q^3,\dots)
$$
is a well-defined homomorphism from $\Lambda$ to $\mathbb R$. In the same way the map
$$
 f(y_1,y_2,\dots)\mapsto f(x_1^{-1},\dots,x_k^{-1},q^{k},q^{k+1},\dots)
$$
is a homomorphism from $\Lambda$ to the ring of Laurent polynomials in $x_1,\dots, x_k$.

Now let $\lambda$ be a Young diagram with $n-1$ rows such that $\lambda_1=-\nu_1$, \dots,
$\lambda_{n-1}=-\nu_{n-1}$ and consider Schur function $s_{\lambda}$ --- element of $\Lambda$. The
image of $s_\lambda$ under the map
$$
 f(y_1,y_2,\dots)\mapsto \frac{f(x_1^{-1},\dots,x_k^{-1},q^k,q^{k+1},\dots)}{f(1,q,q^2,\dots)}
$$
coincides with $Q_k^\nu$.

{\bf Remark 4.} The functions
$$
\frac{s_{\lambda}(x_1,\dots,x_N)}{s_{\lambda}(1,q^{-1},\dots,q^{1-N})}
$$
can be viewed as quantum traces of irreducible representations of the quantized enveloping algebra
$U_\epsilon(\mathfrak{gl}_N)$. Furthermore, one can show that the definition of a $q$--central
measure is related to the branching rules for quantum characters of irreducible representations of
$U_\epsilon(\mathfrak{gl}_N)$. (This is parallel to the fact that the Gelfand-Tsetlin graph itself
is related to the branching rules of irreducible representations of the unitary group.) Thus, it is
natural to expect that the  functions $Q_k^{\nu}$ and $q$--central measures are related to certain
representations of the quantized enveloping algebra $U_\epsilon(\mathfrak{gl}_\infty)$. The author
hopes to address this issue in a later publication.

{\bf Remark 4.} It is natural to ask what happens if one replaces geometric series $q^{1-i}$ in
\eqref{eq_x14} by an arbitrary sequence $\xi_i$. The answer to this question remains unknown.

\subsection{Toeplitz and $q$--Toeplitz matrices}

\label{Subsection_intro_qToeplitz}

Let us explain the connection between extreme characters of the group $U(\infty)$ and total
positivity of Toeplitz matrices. Recall that any extreme character of $U(\infty)$ is a
multiplicative function:
$$
 \chi(U)=\prod_{u_i} \widehat  \chi(u_i),\quad u_i\in\text{Spectrum}(U).
$$
Thus, $\chi$ is uniquely defined by $\widehat \chi(u)$ which is a continuous function on $S^1$.
$\widehat \chi(u)$ can be represented as
$$
 \widehat \chi(u)=\sum_{l\in\mathbb{Z}} c_l u^l.
$$
Introduce an infinite Toeplitz matrix
$$
 c[i,j]\stackrel{def}{=} c_{i-j}.
$$
As a corollary of the fact that $\chi$ is a positive--definite function, one proves that the matrix
$c[i,j]$ is \emph{totally positive}, i.e.\  all minors of $c[i,j]$ are non-negative. Furthermore,
this correspondence is a bijection between extreme characters of $U(\infty)$ and totally positive
infinite Toeplitz matrices such that the sum of the matrix elements along the row equals $1$. (See
\cite{Vo}, \cite{Bo}, \cite{Vk}.)

\smallskip

In order to extend the correspondence between extreme measures and certain matrices to the case of
general $q$ we deform the notion of a Toeplitz matrix. We call a semi-infinite matrix $d[i,j]$,
$i>0$, $j>0$ a \emph{semi-infinite $q$--Toeplitz matrix} if
 \begin{equation}
\label{eq_q_Toeplitz2} d[i,j+1]=d[i-1,j]+(q^{1-j}-q^{1-i})d[i,j]
\end{equation}
for all $i>0$, $j>0$. Here we agree that $d[i,j]=0$ is either $i<1$ or $j<1$. Note that when
$q=1$, the relation \eqref{eq_q_Toeplitz2} turns into
$$
 d[i,j+1]=d[i-1,j].
$$
Hence, a $q$--Toeplitz matrix becomes a Toeplitz matrix.

Recall that according to Theorem \ref{theorem_Main} extreme $q$--central measures correspond to
multiplicative functions $H(x_1)\cdots H(x_N)$. Given a function $H(t)$ we construct a
lower-triangular semi-infinite $q$--Toeplitz matrix $d[i,j]$ in the following way: expand $H$ in
series
$$
 H(t)=\sum_{\ell=0}^{\infty} c_\ell \prod_{i=0}^{\ell-1}(q^{-i}-t)
$$
and let  $d[i,j]$, $i>0$, $j>0$ be a unique semi-infinite $q$--Toeplitz matrix such that
$$
  d[i,1]={c_{i-1}}, \quad i=1,2,\dots.
$$

\emph{Initial minor} of size $N$ of matrix $d[i,j]$ is a minor corresponding to either the first
$N$ columns and arbitrary $N$ rows of matrix $d[i,j]$ or the first $N$ rows and arbitrary $N$
columns of $d[i,j]$.

\begin{proposition} Let $\nu$ be a non-decreasing sequence of non-negative integers $0\le\nu_1\le\nu_2\le\dots$, let
$\mathcal E^{\nu}$ be the extreme $q$--central measure parameterized by $\nu$ and corresponding to
the function $H^\nu$. If $d^{\nu}[i,j]$ is a semi-infinite $q$--Toeplitz matrix constructed by
$H^\nu$, then all initial minors of $d^{\nu}[i,j]$ are non-negative.
\label{Proposition_q_toeplitz_non_neg}
\end{proposition}

A general theorem (see \cite{FZ}) says that if $A$ is a finite non-degenerate matrix with
non-negative initial minors, then $A$ is totally positive, in other words \emph{all} minors of $A$
are non-negative. But, alas, the $q$--Toeplitz matrices corresponding to $q$--central measures are
usually degenerate. These matrices are triangular and some elements on the main diagonal vanish.
Thus, we cannot guarantee that $q$--Toeplitz matrices corresponding to $q$--central measures are
totally positive, we can only claim that their certain top-left corners are. And, indeed,
straightforward computations show that even some matrix elements of these matrices are negative.

However, it might be still interesting to classify all $q$--Toeplitz matrices with non-negative
initial minors. We have the following conjecture here, which is a straightforward analogue of
$q=1$ case.

\begin{conjecture}
 Suppose that $d[i,j]$ is a lower-triangular semi-infinite $q$--Toeplitz matrix with non-negative initial minors
 satisfying normalization condition
 $$
  \sum_{i=1}^{\infty} q^{-(i-1)(i-2)/2}d[i,1]=1,
 $$
 then $d[i,j]$ coincides with one of the matrices from Proposition
 \ref{Proposition_q_toeplitz_non_neg}, i.e.\  $d[i,j]=d^\nu[i,j]$ for some $\nu$.
\end{conjecture}

\subsection{Other branching graphs and general formalism}
The Gelfand--Tsetlin graph and the $q$--Gelfand--Tsetlin graph are two examples of \emph{branching
graphs}. (This term was introduced by Vershik and Kerov.) There is a wide class of problems that
can be stated as a problem of identification the boundary of a branching graph.

A bunch of examples comes from the representation theory of ``big'' groups. Representations of the
infinite symmetric group $S(\infty)$ are related to the boundary of the Young graph (see
\cite{Th_S}, \cite{VK_S}, \cite{Ok_S}), projective representations of $S(\infty)$ are related to
the Schur graph (see \cite{N}, \cite{I}), and (as we already mentioned) representation of the
infinite dimensional unitary group $U(\infty)$ are related to the Gelfand--Tsetlin graph.

However, there are other examples of purely probabilistic and combinatorial nature. Perhaps, the
most known example is De Finetti's theorem (see. e.g \cite[Chapter VII, \S 4]{Feller} or \cite{A})
which states that every probability measure on $\{0,1\}^{\infty}$ invariant with respect to finite
permutations of coordinates is a mixture of Bernoulli measures. Here the underlying branching
graph is the Pascal graph.

Motivated by a problem of population genetics Kingman introduced in \cite{Kingman} the notion of a
partition structure. Kingman's classification of partition structures is equivalent to the
description of the boundary of a certain graph, that is now called the Kingman graph, see also
\cite{Kerov_Kingman}.

Other examples can be found in \cite{OkOl2}, \cite{GO1}, \cite{GO2}, \cite{GO3}, \cite{GP}.

Let us consider a subgraph of $\mathbb{GT}$ consisting of zero--one signatures, i.e.\
$\lambda_1\ge\dots\ge\lambda_N$ such that $1\ge\lambda_1\ge \lambda_N\ge 0$. $q$--central measures
on paths in this subgraph can be identified with central measures on paths in the $q$--Pascal graph
studied by Gnedin and Olshanski \cite{GO1}. These measures are related to $q$--analogues of De
Finetti's theorem. (See also earlier paper by Kerov \cite[Chapter 1.4]{Kerov_book}) The authors of
\cite{GO1} proved that for $0<q<1$ the boundary of the $q$--Pascal graph is parameterized by points
of the set $\{1,q,q^2,\dots\}\cup\{0\}$. This result agrees with our description of the boundary of
the $q$--Gelfand--Tsetlin graph. Indeed, the only extreme $q$--central measures concentrated on
zero--one signatures are those parameterized by $\nu$ with $0\le\nu_i\le 1$ for every $i$. The
extreme central measure on $q$--Pascal graph parameterized by $q^k$ corresponds to $q$--central
measure on $\mathcal T$ parameterized by $\nu=0\le 0\le\dots\le 0\le 1\le 1 \dots$ with exactly $k$
zeros.

In another article \cite{GO2} Gnedin and Olshanski studied a multidimensional generalization of
the model of \cite{GO1}. The common feature of the results of both these two papers and the
present paper is that in the contrast to $q=1$ case, parameters of the extreme $q$--central
measures (and of the measures studied in \cite{GO1}, \cite{GO2}) are discrete.

As a final remark of this section we want to mention the paper \cite{DF} where Diaconis and
Freedman introduced a notion of partial exchangeability. Both central and $q$--central measures of
$\mathcal T$ are particular cases of \emph{partially exchangeable probabilities}.

\subsection{Further connections and developments}

Representation theory of $U(\infty)$ shows numerous connections with $S(\infty)$, i.e.\  inductive
limit of symmetric groups $S_n$. (See \cite{Olsh_sym} and \cite{KOV} for reviews of the
representation theory of $S(\infty)$) For example, while extreme characters of $U(\infty)$ are
related to infinite totally positive Toeplitz matrices, extreme character of $S(\infty)$ are
similarly related to semi-infinite totally positive Toeplitz matrices. The description of the
simplex of $q$--central measures (in particular, their ``right'' definition) on the Young graph
(substitute of the graph $\mathbb{GT}$ for the group $S(\infty)$) is yet to be done.

The are two important problems in the representation theory of the infinite-dimensional unitary
group: identification of all irreducible representations and decomposition of natural
representations into irreducible ones (see \cite{Olsh}). Irreducible representations (again, we
should speak either about finite  factor representations or irreducible spherical representations
here) lead to extreme central measures on $\mathcal T$ and in the present paper we study a
$q$--analogue of these measures. The natural representations associated with $U(\infty)$, in turn,
lead to a remarkable family of  central measures on $\mathcal T$, i.e.\  so called
$(z,w)$--measures. At the moment it is unclear whether one can introduce some natural $q$--central
measures that will serve as a $q$--analogue of $(z,w)$--measures.

 Although, the study of central measures on $\mathcal T$ has a representation-theoretic origin,
later these measures led to a number of very interesting probability models related to random Young
diagrams and random stepped surfaces, see \cite{BK}, \cite{BF}, \cite{B_schur}. The author hopes
that $q$--central measures introduced in the present paper will also provide a source of new
interesting probability models.

\subsection*{Organization of the paper} In Section 2 we introduce all the basic objects of the
study and give a purely combinatorial definition of $q$--central measures. Sections  3 and 4
contain the information about the technical tools that we use for analyzing $q$--central measures.
In Section 3 we introduce various symmetric polynomials and Section 4 is devoted to the
definitions of certain special probability generating functions for $q$--central measures. In
Section 5 we prove the main theorems of the present paper. Some of the proofs are quite involved
and contain lots of technicalities, and we moved these proofs to Section 6. Finally, in Section 7
we discuss the relations with $q$--Toeplitz matrices.

\subsection*{Acknowledgements}

 The  author would like to thank G.~Olshanski for many fruitful discussions at various stages
of this work. The author is grateful to A.~Borodin, L.~Petrov and an anonymous referee for
valuable remarks. The author was partially supported by Moebius Foundation for Young Scientists,
by ``Dynasty'' foundation,  by Russian Foundation for Basic Research --- Centre National de la
Recherche Scientifique [grant 10-01-93114], by the program ``Development of the scientific
potential of the higher school'' and by the additional scholarship of Independent University of
Moscow for graduate students.

\section{Combinatorial setup}
\label{section_probabilistic_setup}
\begin{center}{\bf $\star$ Throughout the paper we assume
$0<q<1$. $\star$}\end{center}

In this section we introduce the basic definitions and give a combinatorial interpretation for the
notion of $q$--centrality. We are going to follow the notations of \cite{Mac} when it is possible.

 A \emph{Young diagram} $\mu$ is a finite collection of boxes arranged in rows with nonincreasing
row lengths $\mu_i$. The total number of boxes in $\mu$ is denoted by $|\mu|$. Every box of $\mu$
has two coordinates $(i,j)$, the first one is increasing from top to bottom and the second one is
increasing from left to right. The top-left corner of a diagram has coordinates $(1,1)$. If we
reflect Young diagram $\mu$ with respect to the diagonal $i=j$, then we obtain a \emph{transposed
diagram $\mu'$}. Row lengths $\mu'_i$ of $\mu'$ are column lengths of $\mu$.

The set of all Young diagrams has a natural partial order by inclusion relation, i.e.\  we write
$\lambda\subset\mu$ if for every $i$ we have $\lambda_i\le\mu_i$.

A \emph{signature} $\lambda$ of size $N$ is an ordered collection of integers
$$
 \lambda_1 \ge \lambda_2 \ge \dots \ge \lambda_N,
$$
we call $\lambda_i$ the $i$th coordinate of $\lambda$.

\begin{figure}[h]
\begin{center}
\noindent{\scalebox{0.6}{\includegraphics{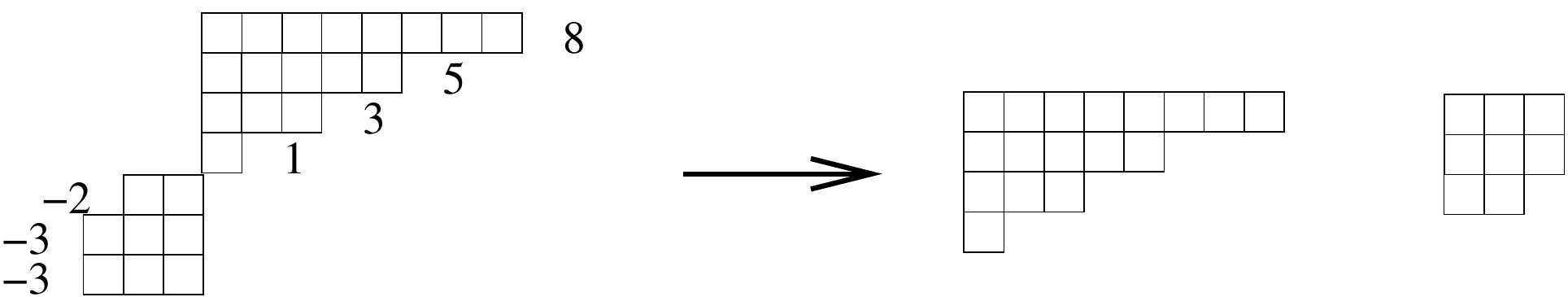}}} \caption{Signature $8\ge5\ge3\ge1\ge0\ge
-2\ge-3\ge-3$ and two corresponding Young diagrams.\label{Figure_signature} }
\end{center}
\end{figure}

Given a signature $\lambda$ we construct two Young diagrams $\lambda^+$ and $\lambda^-$ which are
called a \emph{positive diagram} and a \emph{negative diagram}, respectively. The row lengths of
the former one are all the positive $\lambda_i$s, while the row lengths of the latter are absolute
values of negative $\lambda_i$s. An example is shown in Figure \ref{Figure_signature}. Let
$|\lambda|$ denote the sum of coordinates of $\lambda$. Clearly,
$|\lambda|=|\lambda^+|-|\lambda^-|$.

Let us write $\lambda\le\mu$ for two signatures of the same size if $\lambda_i\le\mu_i$ for every
$i$. Clearly, $\lambda\le\mu$ is equivalent to $\lambda^+\subset\mu^+$ and $\mu^-\subset
\lambda^-$.

Let $\mathbb{GT}_N$ denote the set of all signatures of size $N$ and let $\mathbb{GT}_N^+$ denote
the set of al signatures of size $N$ with nonnegative coordinates. Every element of
$\mathbb{GT}_N^+$ can be identified with a Young diagram, however note, that a signature has an
additional information,  its size $N$.

\emph{A path} $\tau\in\mathcal T$ in $\mathbb{GT}$ is a sequence of signatures
$\{\lambda(n)\in\mathbb{GT}_n\}_{1,\dots}$ such that for every $n$, $\lambda(n)\prec\lambda(n+1)$.
In other words,
$$
 \lambda(n+1)_1\ge\lambda(n)_1\ge\lambda(n+1)_2\ge\dots\ge\lambda(n)_{n}\ge\lambda(n+1)_{n+1}.
$$
 Paths in $\mathbb{GT}$ are usually called \emph{Gelfand-Tsetlin schemes} in
representation-theoretic literature.

Let $\mathcal T_{N}$ denote the set of all paths $\tau(1)\prec\dots\prec\tau(N)$ of length $N$ with
$\tau(i)\in\mathbb{GT}_i$. The set $\mathcal T$ of all infinite paths is a projective limit of sets
$\mathcal T_N$. (Here projection is just a removal of the last step of a path.)

We define $\mathbb{GT}^+$ to be the part of $\mathbb{GT}$ consisting of signatures with nonnegative
coordinates. Let $\mathcal T^+$ and $\mathcal T_N^+$ denote the corresponding sets of paths.

A (semistandard Young) \emph{tableau} $T$ of shape $\lambda\in\mathbb{GT}_N^+$ is an assignment of
numbers $1,\dots,N$ to boxes of the Young diagram $\lambda$ in such a way that the numbers are
increasing along the columns and non-decreasing along the rows. Given a path $\tau\in\mathcal
T_{N}$ we construct two Young tableaux $T^+_{\tau}$ and $T^-_{\tau}$ as follows: Shape of
$T^+_{\tau}$ is $\tau(N)^+$, $T^+_{\tau}(i,j)=k$ if and only if
$(i,j)\in\tau(k)^+\setminus\tau(k-1)^+$, where we agree that $\tau(0)^+=\emptyset$. Similarly,
shape of $T^-_{\tau}$ is $\tau(N)^-$, $T^-_{\tau}(i,j)=k$ if and only if
$(i,j)\in\tau(k)^-\setminus\tau(k-1)^-$.

\begin{figure}[h]
\begin{center}
\noindent{\scalebox{0.9}{\includegraphics{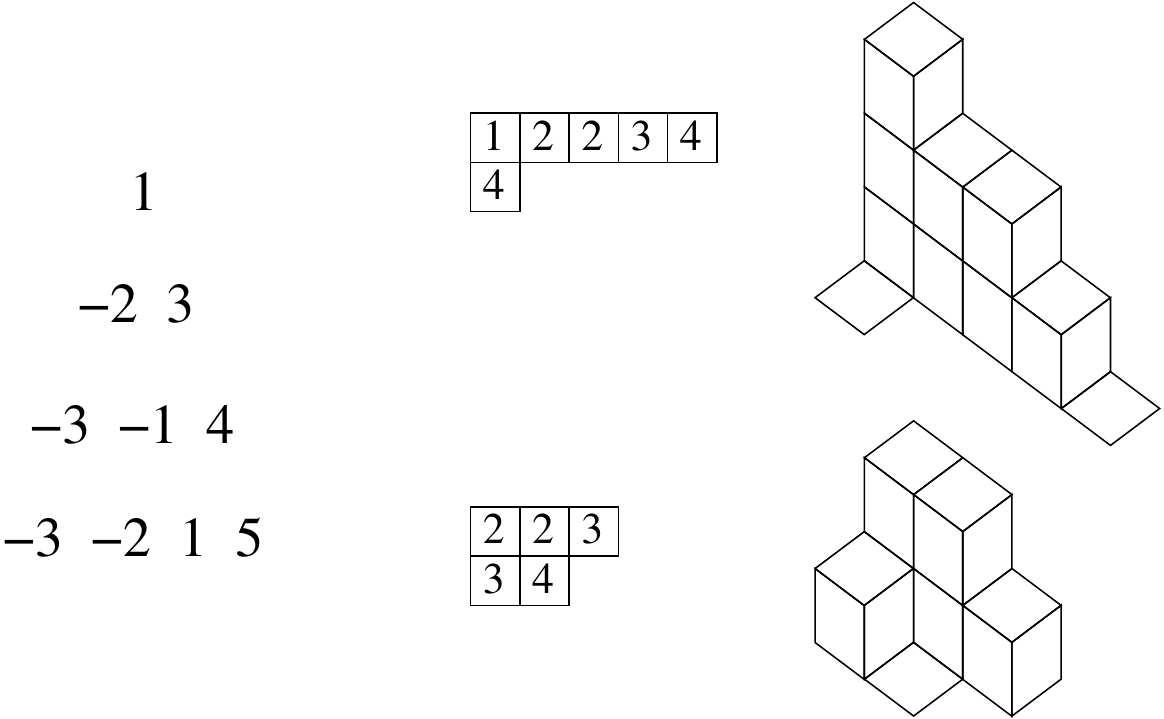}}} \caption{Finite path from ${\mathcal
T_4}$ , corresponding semistandard Young tableaux and 3D Young diagrams.\label{Figure_3d} }
\end{center}
\end{figure}

Every tableau $T$ of shape $\lambda\in\mathbb{GT}_N^{+}$ corresponds to a \emph{3D Young diagram}
in the following way: For every $k$ put $N-k$ unit cubes on all the boxes of $T$ with number $k$.
The union of all these unit cubes is the desired 3D Young diagram. An example of the above
procedure is shown in Figure \ref{Figure_3d}.

Let $V(T)$ denote the volume of 3D Young diagram corresponding to $T$.

\begin{lemma}
  The following formulas hold for an arbitrary $\tau\in\mathcal T_N$:
$$
 V(T^+_\tau)-V(T^-_\tau)=\sum_{i=1}^{N-1} |\tau(i)|=\sum\limits_{i=1}^{N}
(N-i)(|\tau(i)|-|\tau(i-1)|),
$$
where we agree that $|\tau(0)|=0$
\end{lemma}
We leave the proof to the reader.

\medskip

 For a finite path $\tau\in\mathcal T_N$ let $C_\tau$ be a corresponding cylinder set in $\mathcal
T$, i.e.\
$$
 C_\tau=\{\lambda(1)\prec\lambda(2)\prec\dots \in \mathcal T \mid
 \lambda(1)=\tau(1),\dots,\lambda(N)=\tau(N)\}.
$$
We equip $\mathcal T$ with a $\sigma$-algebra spanned by all cylinder sets. We are going to
consider various probability measures on $\mathcal T$. For any finite path $\tau$ we usually write
$P(\tau)$ instead of $P(C_\tau)$ where it leads to no confusion.

For any probability measure $P$ on $\mathcal T$ let $P_N$ be its projection on $\mathbb{GT}_N$,
i.e.\
$$
 P_N (\{\lambda\})=P(\{\tau\in\mathcal T:\quad \tau(N)=\lambda\}).
$$
To simplify the notation we usually write $P_N(\lambda)$ instead of $P_N(\{\lambda\}$

Recall that a probability measure $P$ on $\mathcal T$ is called \emph{$q$--central} if
probabilities of paths $\tau\in\mathcal T_N$ ending at the same signature $\lambda$ are
proportional to $q^{V(T^+_\tau)-V(T^-_{\tau})}$, i.e.\  for any $\tau\in\mathcal T_N$ such that
$\tau(N)=\lambda$ we have:
$$
 P(\tau)=P_N(\lambda)\dfrac {q^{V(T^+_\tau)-V(T^-_{\tau})}}{\sum\limits_{\theta\in\mathcal T_N:\,
\theta(N)=\lambda} q^{V(T^+_\theta)-V(T^-_{\theta})}}.
$$

The main goal of the present paper is to describe the convex set of all $q$--central probability
measures on $\mathcal T$.

\subsection*{Statistical mechanical interpretation}

Before starting the proofs of the main theorems let us show a way to interpret $q$--central
measures in the spirit of statistical mechanics. (The material of this subsection is not further
used throughout the paper.)

\begin{figure}[h]
\begin{center}
\noindent{\scalebox{0.25}{\includegraphics{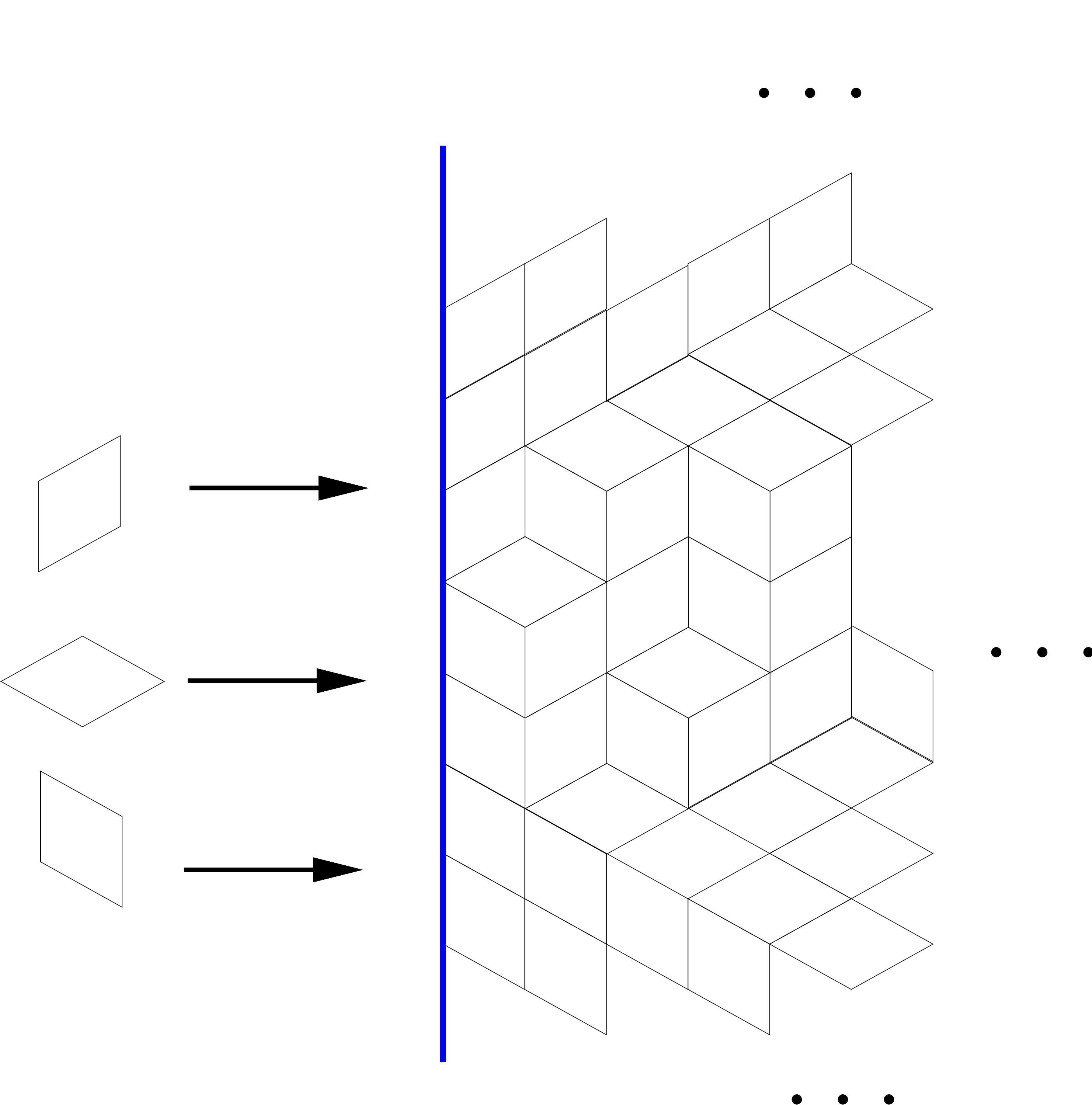}}} \caption{Lozenge tiling of
the halfplane \label{Figure_tiling} }
\end{center}
\end{figure}

Consider a tiling of the halfplane by rhombuses of 3 types (such rhombuses are usually called
\emph{lozenges}) as shown in Figure \ref{Figure_tiling}.

\begin{figure}[h]
\begin{center}
\noindent{\scalebox{0.25}{\includegraphics{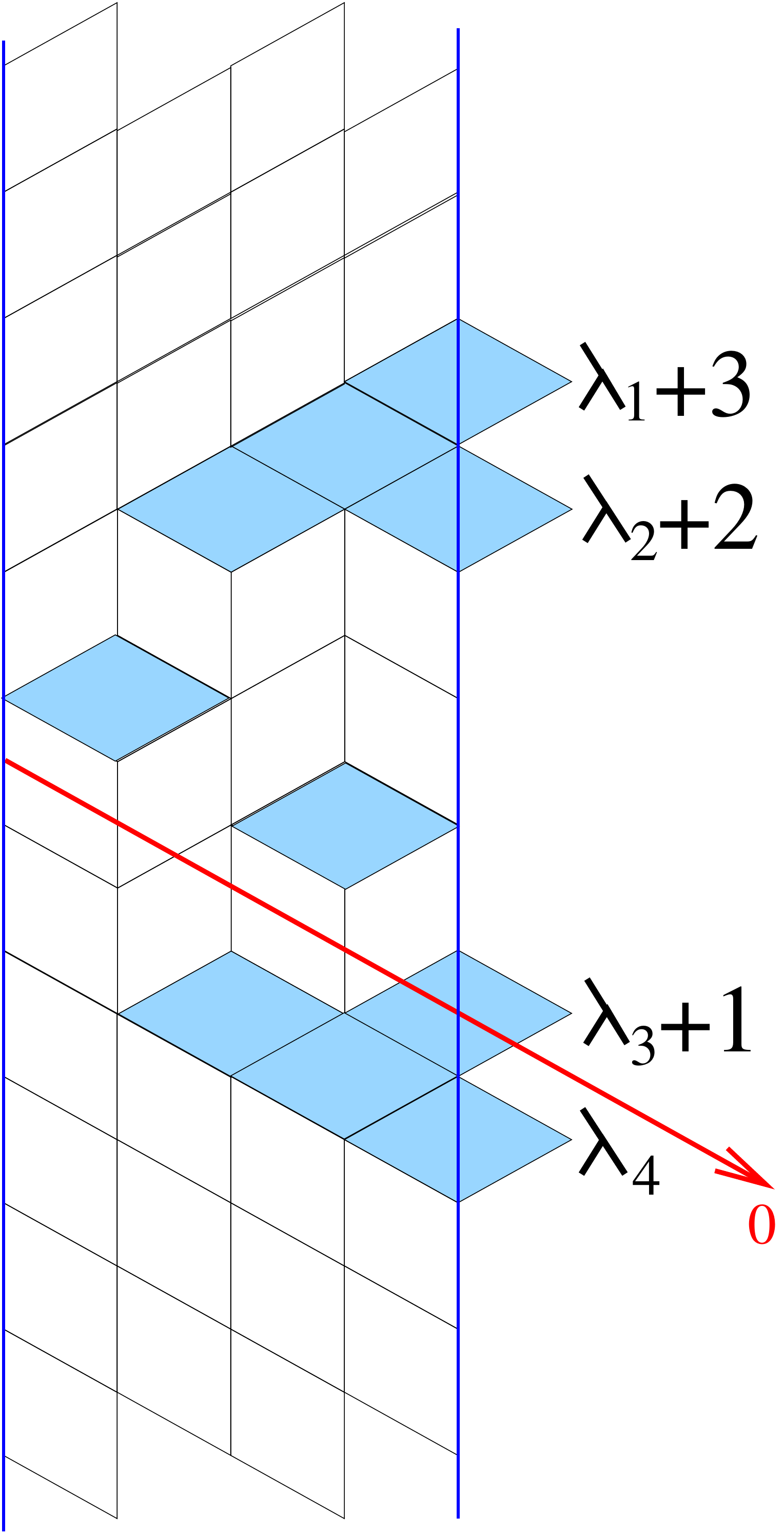}}} $\quad\quad$
\noindent{\scalebox{0.5}{\includegraphics{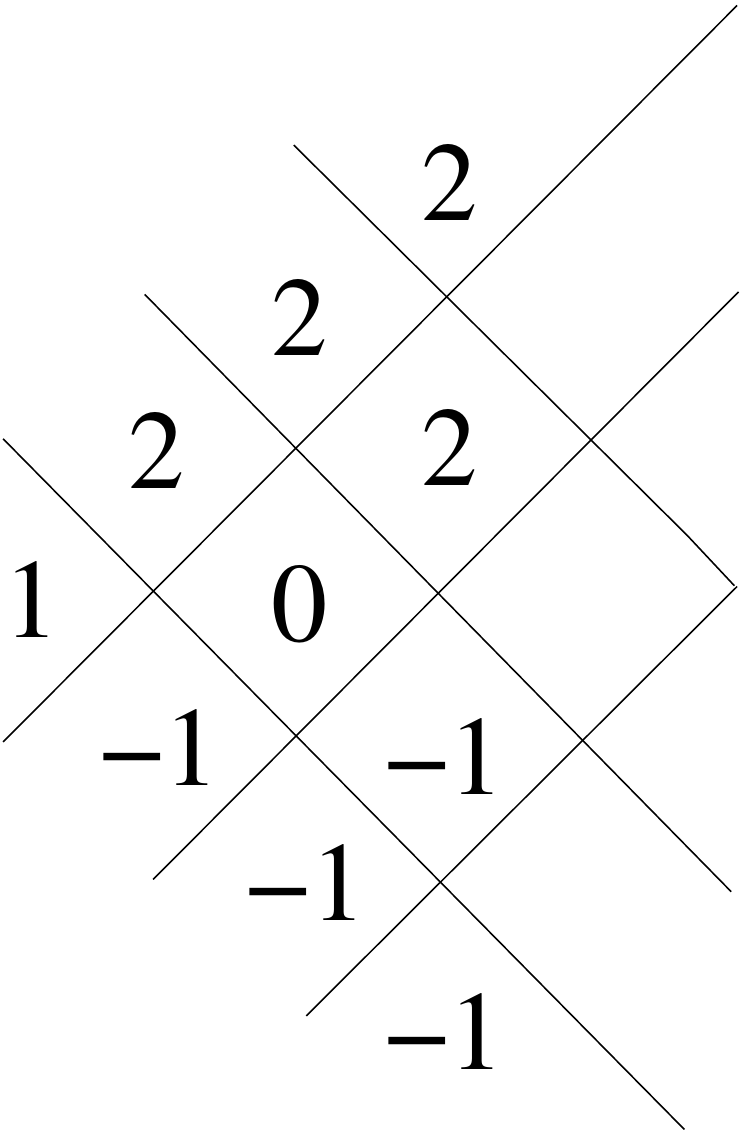}}} \caption{Part of the tiling corresponding
to the path in $\mathbb{GT}$ $(1)\prec (2\ge -1)\prec (2 \ge 0 \ge -1)\prec (2\ge 2\ge -1 \ge
-1)$\label{Figure_tiling_coords}}
\end{center}
\end{figure}

A tiling is uniquely defined by the positions of \emph{horizontal} lozenges (i.e.\ lozenges with
no sides parallel to the vertical line; such lozenge is the middle one in the left part of Figure
\ref{Figure_tiling}). Given a path $\lambda(1)\prec\lambda(2)\prec\dots$ in $\mathbb{GT}$ consider
a tiling with coordinates of horizontal lozenges $(N,\lambda(N)_{i}+N-i-1)$ as shown in Figure
\ref{Figure_tiling_coords}, left panel.

There is a well--known bijection between tilings and discrete stepped surfaces. To obtain such
surface we arrange coordinates of the signatures of a path as in Figure
\ref{Figure_tiling_coords}, right panel; we get some function, which is usually called \emph{the
height function}. Graph of this function is a stepped surface and projection of the surface in the
direction $(1,1,1)$ is precisely a lozenge tiling.

$q$--Centrality of probability measures on paths in $\mathbb{GT}$ turns into the following Gibbs
property of measures on tilings: given positions of horizontal lozenges on a vertical line,
tilings of the strip to the left from this line are distributed with weight $q^{vol}$, where $vol$
is the volume enclosed underneath the corresponding stepped surface. Note that the differences
between the volumes underneath surfaces corresponding to the distinct tilings of the strip are
finite. Thus, the weight $q^{vol}$ is well defined.

Therefore, the problem solved in the present paper can be restated as the classification of
certain Gibbs measures on the lozenge tilings of the halfplane.

\section{Symmetric polynomials}

In this section we introduce various symmetric functions and the algebras they belong to. These
functions play an important role in our study of $q$--central measures.

\subsection{Schur functions and factorial Schur functions}

 Recall that a (rational) Schur function $s_{\lambda}(x_1,\dots,x_N)$ parameterized by an arbitrary
 signature $\lambda_1\ge\lambda_2\ge\dots\lambda_N$, is a symmetric Laurent polynomial given by
 $$
  s_{\lambda}=\frac{\det\left[x_i^{\lambda_j+N-j}\right]_{i,j=1,\dots,N}}{\prod\limits_{i<j}(x_i-x_j)}.
 $$
 If $\lambda\in\mathbb{GT}^+$, then $s_{\lambda}$ is an ordinary symmetric polynomial.
 It is known that rational Schur functions form a linear basis in the space of all symmetric
 Laurent polynomials. We recommend \cite{Mac} as a general source of information about symmetric polynomials and Schur
 functions.

 \begin{proposition}[The branching rule for Schur functions]
 \label{proposition_Schur_branching_rule} For any $\lambda\in\mathbb{GT}_N$ we have:
 $$
 s_\lambda(x_1,\dots,x_N)=\sum_{\mu\prec\lambda} s_\mu(x_1,\dots,x_{N-1}) x_N^{|\lambda|-|\mu|},
 $$
 where $|\lambda|$ and $|\mu|$ stand for the sum of the coordinates of signatures $\lambda$ and
 $\mu$ correspondingly.
 \end{proposition}
 See e.g.\ \cite[Chapter 1, Section 5]{Mac} for the proof. Iterating the branching rule for Schur
functions we get the following:
 \begin{proposition}[The combinatorial formula]
  For any $\lambda\in\mathbb{GT}_k$ we have:
 $$
  s_{\lambda}(x_1,\dots,x_N)=\sum_{\tau(0)\prec\dots\prec\tau(N)} x_1^{|\tau(1)|-|\tau(0)|}\dots
 x_N^{|\tau(N)|-|\tau(N-1)|},
 $$
  where the sum is taken over all paths in $\mathbb{GT}$ $\tau(0)\prec\dots\prec\tau(N)$ such that
 $\tau(N)=\lambda$ and $|\tau(i)|$ stands for the sum of the coordinates of signature $\tau(i)$.
 \end{proposition}

From now on assume that $\lambda\in\mathbb{GT}_N^+$, i.e.\
$\lambda_1\ge\lambda_2\ge\dots\ge\lambda_N\ge 0$.

Let $\{a_n\}_{n\in \mathbb Z}$ be any sequence of numbers. For any $r\ge 0$ let
$$
 (x\mid a)^r=(x+a_1)\dots (x+a_r).
$$
 Factorial Schur function
$s_\lambda(x\mid a)(x_1,\dots,x_N)$ is a symmetric polynomial in variables $x_1,\dots x_N$ defined
through
$$
  s_{\lambda}(x\mid a)=\frac{\det\left[(x_i\mid a)^{\lambda_j+N-j}\right]_{i,j=1,\dots,N}}{\prod\limits_{i<j}(x_i-x_j)}.
$$
See \cite[6th Variation]{Mac2} for the properties of these polynomials.

\begin{proposition}[Combinatorial formula for factorial Schur polynomials]
 For any $\lambda\in\mathbb{GT}_N^+$ we have:
$$
 s_\lambda(x\mid a)=\sum_T\prod_{(i,j)\in\lambda}(x_{T(i,j)}+a_{T(i,j)+j-i}),
$$
where the sum is taken over all semistandard Young tableau $T(i,j)$ of shape $\lambda$ filled with
numbers $1, \dots, N$.
\end{proposition}

\subsection{$q$--interpolation Schur polynomials}

\emph{$q$--Interpolation Schur polynomials $s^*_{\lambda}(x;q)$} are factorial Schur polynomials
with $a_i=-q^{i-N}$:
 $$
  s_\lambda^*(x_1,\dots,x_N;q)=s_\lambda(x\mid a), \quad a=\{-q^{i-N}\}.
 $$

This polynomials are a particular case of Macdonald interpolation polynomials for the case $q=t$.
(See \cite{Kn}, \cite{S}, \cite{Ok_shifted_Macdonald}.)

Let us introduce some notations. Suppose $\mu\in\mathbb{GT}^+_N$, and recall that $\mu$ can be
identified with a Young diagram with not more than $N$ rows. Let ${\mu_1\ge\mu_2\ge\dots}$ be the
row lengths of $\mu$ and let $\mu'_1\ge\mu'_2\ge\dots$ be the row lengths of the transposed
diagram $\mu'$, or, equivalently, the column lengths of $\mu$. For any box $(i,j)\in\mu$  set
$c(i,j)=j-i$ and $h(i,j)=\mu_i-i+\mu'_j-j+1$. Denote $n(\mu)=\sum (i-1)\mu_i$.

\begin{proposition}[Interpolation property, \cite{Ok_shifted_Macdonald}]
\label{proposition_interpolation_property} The $q$--interpolation Schur polynomial
$s^*_{\mu}(x_1,\dots,x_N;q)$ is the unique symmetric polynomial in $x_1,\dots,x_N$ such that:
\begin{enumerate}
\item $\deg(s^*_{\mu}(x_1,\dots,x_N;q))=|\mu|$
\item $s^*_{\mu}(q^{\mu-\delta};q)=q^{n(\mu')-2n(\mu)}\prod\limits_{(i,j) \in \mu}(q^{h(i,j)}-1)$
\item $s^*_{\mu}(q^{\lambda-\delta};q)=0$ for all positive signatures (Young diagrams) $\lambda\in\mathbb{GT}_N^+$ such that $\mu \nsubseteq
\lambda$,
\end{enumerate}
where $q^{\lambda-\delta}$ is $(q^{\lambda_1},q^{\lambda_2-1},\dots,q^{\lambda_N-N+1})$
\end{proposition}

Rewriting the combinatorial formula for the factorial Schur polynomials we get:
\begin{proposition} We have
\label{Proposition_combinatoral_for_interp_Schur}
$$
 s^*_{\mu}(x_1,\dots,x_N;q)=\sum_{T} \prod\limits_{(i,j)\in\mu}
 \left(x_{T(i,j)}-q^{j-i+T(i,j)-N}\right),
$$
where the sum is taken over all semistandard Young tableau $T(i,j)$ of shape $\lambda$ filled with
numbers $1, \dots, N$.
\end{proposition}

In one-dimensional case interpolation Schur polynomials are enumerated by nonnegative integers and
$$s^*_{k}(x;q)=(x-1)\dots (x-q^{k-1})$$

Polynomials $s^*_{\mu}(x_1,x_2\cdot,\dots, x_N;q)$ form a linear basis in the space of symmetric
polynomials in $N$ variables. The leading term of $s^*_{\mu}(x_1,x_2\cdot ,\dots, x_N \cdot ;q)$
coincides with ordinary Schur function $s_{\mu}(x_1,x_2,\dots, x_N)$.

The basis $s^*_\mu$ is connected with the basis $s_\mu$ via the following formulas due to Okounkov
\cite{Ok_bynomial}

\begin{equation}
\label{eq_schur_through_interpolation}
 s_{\lambda}(x_1,x_2,\dots,x_k)
 =\sum_{\mu\in\mathbb{GT}_k^+}\frac{s^*_\mu(q^{\lambda-\delta};q)}{s^*_\mu(q^{\mu-\delta};q)}\frac{s_{\lambda}(1,q^{-1},\dots,
 q^{1-k})}{s_{\mu}(1,q^{-1},\dots,q^{1-k})}{s^*_{\mu}(x_1,\dots,x_k;q)},
\end{equation}

\begin{equation}
\label{eq_interpolation_through_schur}
 s_{\lambda}^*(x_1,x_2,\dots,x_k;q)
=\sum_{\mu\in\mathbb{GT}_k^+}\frac{s^*_\mu(q^{-(\lambda-\delta)};q^{-1})}{s^*_\mu(q^{-(\mu-\delta)});q^{-1})}\frac{s_{\lambda}^*(0,\dots,0;q
)}{s_{\mu}^*(0,\dots,0;q)}{s_{\mu}(x_1,\dots,x_k)}.
\end{equation}

Observe that in both formulas only diagrams $\mu$ such that $\mu\subset\lambda$ give a nonzero
contribution. Thus, both sums are actually finite. Also note that the formulas
\eqref{eq_schur_through_interpolation} and \eqref{eq_interpolation_through_schur} look similar, we
are going to further use this symmetry.

\subsection{Algebras of symmetric functions}
\label{subsection_algebras} Let $\Lambda_N$ be the graded algebra of symmetric polynomials in
$x_1,\dots,x_N$. Let $\Lambda$ be a projective limit of graded algebras $\Lambda_N$ with respect to
projections $\rho_N$:
$$
 \rho_N: \Lambda_{N+1}\to\Lambda_{N},\quad  \rho_N(f(x_1,\dots,x_{N+1}))=f(x_1,\dots,x_N,0).
$$
$\Lambda$ is usually called \emph{the algebra of symmetric functions}. There are 3 well-known
systems of algebraic generators of $\Lambda$. They are \emph{Newton power sums} $p_k$
$$
 p_k=\sum_{i=1}^{\infty} x_i^k,
$$
\emph{elementary symmetric functions} $e_k$
$$
 e_k=\sum_{\ell_1<\ell_2<\dots<\ell_k} x_{\ell_1}\cdots x_{\ell_k}
$$
and \emph{complete symmetric functions} $h_k$
$$
 h_k=\sum_{\ell_1\le \ell_2\le\dots\le\ell_k} x_{\ell_1}\cdots x_{\ell_k}.
$$
Let $\lambda=(\lambda_1,\dots,\lambda_k)\in\mathbb{GT}_k^+$ and let $\widehat
\lambda\in\mathbb{GT}_{k+1}^+$ be a signature obtained by adding one zero to $\lambda$, i.e.\
$\widehat\lambda=(\lambda_1,\dots,\lambda_k,0)$. Schur polynomials corresponding to non-negative
signatures are stable in the sense that
$$
 s_{\widehat \lambda}(x_1,\dots,x_k,0)=s_{\lambda}(x_1,\dots,x_k)
$$
for any $\lambda\in\mathbb{GT}_k^+$. Thus, the Schur polynomial corresponding to a non-negative
signature $\lambda$ defines an element of $\Lambda$ that we denote $s_\lambda$. Functions
$s_\lambda$ form a linear basis in algebra $\Lambda$.

\medskip

Let $\widehat \Lambda$ be a projective limit of filtered algebras $\widehat \Lambda_N$, where
$\widehat \Lambda_N$ is the ordinary algebra of symmetric polynomials in $N$ variables filtered by
the degree, and the projections $\widehat \rho_N$ are given by the following formula:
$$
 \widehat \rho_N: \widehat \Lambda_{N+1}\to\widehat \Lambda_{N},\quad \widehat
\rho_N(f(x_1,\dots,x_{N+1}))=f(x_1,\dots,x_N,q^{-N})
$$

Polynomials $s^*_{\mu}$ are stable in the sense that
$$
 s^*_{\widehat \mu}(x_1,\dots,x_N,q^{-N};q)=s^*_{\mu}(x_1,\dots,x_N;q)
$$
for any $\mu\in\mathbb{GT}_N^+$. This fact follows from Proposition
\ref{proposition_interpolation_property}. Thus, these polynomials can be viewed as elements of a
filtered algebra $\widehat \Lambda$. Furthermore, $s^*_{\mu}$ form a linear basis of $\widehat
\Lambda$.

\section{$q$--Central measures and probability generating functions}
In this section we state a number of propositions about $q$--central measures and probability
generating functions related to them. Some of the proofs are quite technical and we omit them for
the convenience of the reader. All the missing proofs are given in Section
\ref{Subsection_proofs_probability_gen_functions}.

For any $\mu\in\mathbb{GT}_N$ denote
$$
 {\rm Dim}_q(\mu)=\sum_{\tau\in \mathcal T_N,\, \tau(N)=\mu} q^{|\tau(1)|+\dots+|\tau(N-1)|}.
$$

\begin{lemma} We have
$$
 {\rm Dim}_q(\lambda)=s_{\lambda}(1,q,\dots,q^{N-1}).
$$
\end{lemma}
\begin{proof}
 Applying the combinatorial formula for Schur functions we obtain
 \begin{multline*}
  s_{\lambda}(1,\dots,q^{N-1})=s_{\lambda}(q^{N-1},\dots,1)\\ =\sum_{\tau(0)\prec\dots\prec\tau(N)}
  q^{(N-1)\cdot\left(|\tau(1)|-|\tau(0)|\right)+\dots+
  0\cdot\left(|\tau(N)|-|\tau(N-1)|\right)}\\=\sum_{\tau(0)\prec\dots\prec\tau(N)}
  q^{|\tau(1)|+|\tau(2)|+\dots+|\tau(N-1)|}={\rm Dim}_q(\lambda)
 \end{multline*}
\end{proof}

Let $\lambda\in\mathbb{GT}_{N+1}$, $\mu\in\mathbb{GT}_{N}$. We define \emph{cotransitional
probability $P(\lambda\to\mu)$}:
$$
P(\lambda\to\mu)=\begin{cases}q^{|\mu|}\frac{{\rm Dim}_q(\mu)}{{\rm Dim}_q(\lambda)},\quad
\mu\prec\lambda,\\0,\quad otherwise.\end{cases}
$$

\begin{proposition}
 For any $\lambda\in\mathbb{GT}_{N+1}$ we have
$$
 \sum_{\mu\in\mathbb{GT}_N} P(\lambda\to\mu)=1.
$$
\end{proposition}
\begin{proof}
 Proposition \ref{proposition_Schur_branching_rule} implies that
$$
 s_\lambda(1,\dots,q^N)=\sum_{\mu\prec\lambda} s_\mu(q,q^2,\dots,q^{N})=\sum_{\mu\prec\lambda}
q^{|\mu|} s_\mu(1,q,\dots,q^{N-1}).
$$
Dividing by $s_\lambda(1,\dots,q^N)$ we get the desired equality.
\end{proof}

Suppose that $P_N$ and $P_{N+1}$ are probability distributions on $\mathbb{GT}_N$ and
$\mathbb{GT}_{N+1}$, respectively. We call $P_N$ and $P_{N+1}$ \emph{$q$--coherent} if for any
$\mu\in\mathbb{GT}_N$:
$$
 P_{N}(\mu)=\sum_{\lambda\in\mathbb{GT}_{N+1}} P_{N+1}(\lambda) P(\lambda\to\mu).
$$

We call probability distributions $P_1,P_2,\dots$ on $\mathbb{GT}_1,\mathbb{GT}_2,\dots$
respectively \emph{a $q$--coherent system}, if $P_i$ and $P_{i+1}$ are $q$--coherent for every
$i=1,2,\dots$.

Note that for any probability distribution $P_N$ on $\mathbb{GT}_N$, there exist unique
distributions $P_1,\dots,P_{N-1}$ on $\mathbb{GT}_1,\dots,\mathbb{GT}_{N-1}$, respectively, such
that $P_1,\dots,P_N$ is a $q$--coherent system.

Let $P$ be an arbitrary probability measure on $\mathcal T_N$. Let $P_k$ denote  a projection of
$P$ on $\mathbb{GT}_k$. The following two propositions are proved in Section
\ref{Subsection_proofs_probability_gen_functions}:
\begin{proposition}
 \label{Proposition_coherence_from_centrality} If measure $P$ is such that
 \begin{equation}
 \label{eq_q_centrality} P(\tau(1)\prec\dots\tau(N))=\frac{q^{|\tau(1)|+\dots+|\tau(N-1)|}}{{\rm
Dim}_q(\tau(N))} P_N(\tau(N)),
 \end{equation}
 for any path $\tau\in\mathcal T_N$, then $P_1,P_2,\dots,P_N$ is a $q$--coherent system. In
particular, if $P$ is a $q$--central measure on $\mathcal T$, then $P_1,P_2,\dots$ is a
$q$--coherent system.
\end{proposition}

\begin{proposition}
\label{Proposition_from_coherent_system_to_measure} Let $P_1,P_2,\dots$ be a $q$--coherent system.
There exists a unique $q$--central measure $P$ such that
 $P_k$ is a projection of $P$ on $\mathbb{GT}_k$ for every $k=1,2,\dots$.
\end{proposition}

Next, we want to introduce a convenient tool for studying $q$--central measures and $q$--coherent
systems.

Suppose that $P$ is a probability measure on $\mathbb{GT}_N$. \emph{$q$--Schur generating function}
of $P$ is a symmetric function in $x_1,\dots,x_N$ given by:
\begin{equation}
\label{eq_Schur_gen_function} {\mathcal S}(x_1,\dots,x_N;P)=\sum_{\mu\in\mathbb{GT}_N} P(\mu)
\frac{s_\mu(x_1,\dots,x_N)}{s_\mu(1,q^{-1},\dots,q^{1-N})}.
\end{equation}
Note that when $N=1$, this definition turns into the usual definition of probability generating
function:
$$
 F(t)=\sum_\ell c_\ell t^\ell.
$$

For every $N$ we define
$$
 T_N=\{(x_1,\dots,x_N)\in\mathbb C^N\mid |x_i|=q^{1-i}\}
$$
and
$$
 D_N=\{(x_1,\dots,x_N)\in\mathbb C^N\mid |x_i|\le q^{1-i}\}.
$$

\begin{proposition}
 \label{proposition_Schur_gen_function_converges} The series \eqref{eq_Schur_gen_function}
 converges uniformly on $T_N$. If ${\rm supp}(P)\subset \mathbb{GT}_N^+$, then the series
 \eqref{eq_Schur_gen_function} converges uniformly on $D_N$.
\end{proposition}
\begin{proof}
 The proposition follows from the fact that
 $$
  \left|\frac{s_\mu(x_1,\dots,x_N)}{s_\mu(1,q^{-1},\dots,q^{1-N})}\right|\le 1,
 $$
 for all $(x_1,\dots,x_N)\in T_N$, and if $\mu\in\mathbb{GT}_N^+$, then
$$
  \left|\frac{s_\mu(x_1,\dots,x_N)}{s_\mu(1,q^{-1},\dots,q^{1-N})}\right|\le 1,
 $$
 for all $(x_1,\dots,x_N)\in D_N$.

  These inequalities, in turn, easily follow from the combinatorial formula for rational Schur
functions.
\end{proof}
 The following proposition is proved in Section \ref{Subsection_proofs_probability_gen_functions}:
\begin{proposition}
 \label{Proposition_coherency_in_terms_of_Schur_polynomials} Suppose that $P_N$ and $P_{N+1}$ are
two probability measures on $\mathbb{GT}_N$ and $\mathbb{GT}_{N+1}$ respectively. $P_N$ and
$P_{N+1}$ are $q$--coherent if and only if
$$
 {\mathcal S}(x_1,\dots,x_N;P_{N})={\mathcal S}(x_1,\dots,x_N,q^{-N};P_{N+1}).
$$
\end{proposition}

If ${\rm supp}(P)\subset \mathbb{GT}^+_N$  (i.e.\  $P$ is a probability measure on
$\mathbb{GT}^+_N$) then we define \emph{$q$--interpolation Schur generating function} of $P$
through
\begin{equation}
\label{eq_interpolation_Schur_gen_function} {\mathcal
S^*}(x_1,\dots,x_N;P)=\sum_{\mu\in\mathbb{GT}^+_N} P(\mu)
\frac{s^*_\mu(q^{N-1}x_1,\dots,q^{N-1}x_N;q^{-1})}{s^*_\mu(0,\dots,0;q^{-1})}.
\end{equation}

\begin{proposition}
\label{proposition_convergence_in_def_of_class_F_star} Series
\eqref{eq_interpolation_Schur_gen_function} converges uniformly on compact subsets of $\mathbb
C^N$.
\end{proposition}
\begin{proof}
 Using the combinatorial formula for $q$--interpolation polynomials we get:

 \begin{multline*}
s^*_\mu(q^{N-1}x_1,q^{N-1}x_2,\dots,q^{N-1}x_N;q^{-1})\\=q^{(N-1)|\mu|}\sum_{T}\prod_{(i,j)\in\mu}\left(x_{T(i,j)}-q^{i-j-1+T(i,j)}\right)
  \\=q^{(N-1)|\mu|}(-1)^{|\mu|}\left(\prod_{(i,j)\in\mu}q^{i-j-1}\right)\\ \times \sum_{T}\left(\prod_{(i,j)\in\mu}q^{T(i,j)}\right)
  \prod_{(i,j)\in\mu}\left(1-x_{T(i,j)}q^{j-i-T(i,j)+1}\right).
 \end{multline*}

 Let $M$ be a constant such that $|x_i|<M$ for all $i$. We have
\begin{multline*}
\left|\frac{s^*_\mu(q^{N-1}x_1,q^{N-1}x_2,\dots,q^{N-1}x_N;q^{-1})}{s^*_\mu(0,\dots,0;q^{-1})}\right|\\=\frac{\left|\sum_{T}\left(\prod_{(i,j)\in\mu}q^{T(i,j)}\right)\prod_{(i,j)\in\mu}\left(1-x_{T(i,j)}q^{j-i-T(i,j)+1}\right)\right|}{\sum_{T}\left(\prod_{(i,j)\in\mu}q^{T(i,j)}\right)}
 \\ \le
 \frac{\sum_{T}\left(\prod_{(i,j)\in\mu}q^{T(i,j)}\right)\prod_{(i,j)\in\mu}\left(1+Mq^{j-i-T(i,j)+1}\right)}{\sum_{T}\left(\prod_{(i,j)\in\mu}q^{T(i,j)}\right)}\\
 \le \max_T \prod_{(i,j)\in\mu} \left(1+Mq^{j-i-T(i,j)+1}\right)
  \\ \le \max_T \exp\left(\sum_{(i,j)\in\mu} Mq^{j-i-T(i,j)+1}\right)
 \\ \le  \exp\left(M\sum_{(i,j)\in\mu} q^{j-i-N}\right).
\end{multline*}
 Since $\mu$ has at most $N$ rows,
$$\sum_{(i,j)\in\mu} q^{j-i-N}< N (q^{1-2N}+q^{2-2N}+q^{3-2N}+\dots)=\frac{Nq^{1-2N}}{1-q}.
$$

We conclude that if $(x_1, \dots, x_N) $ is such that $|x_i|<M$ for every $i$, then
$$
\left|\frac{s^*_\mu(q^{N-1}x_1,q^{N-1}x_2,\dots,q^{N-1}x_N;q^{-1})}{s^*_\mu(0,\dots,0;q^{-1})}\right|<A(M).
$$

Consequently,
$$
\sum\limits_{\mu\in\mathbb{GT}^+_N}c_\mu
\left|\frac{s^*_\mu(q^{N-1}x_1,q^{N-1}x_2,\dots,q^{N-1}x_N;q^{-1})}{s^*_\mu(0,\dots,0;q^{-1})}\right|
< \infty,
$$
\eqref{eq_interpolation_Schur_gen_function} absolutely converges and this convergence is uniform on
compact subsets of $\mathbb C^N$.
\end{proof}

The following proposition is proved in Section \ref{Subsection_proofs_probability_gen_functions}:
\begin{proposition}
\label{Proposition_coherency_in_terms_of_interpolation_polynomials}
 Suppose that $P_N$ and $P_{N+1}$ are two probability measures on $\mathbb{GT}^+_N$ and
$\mathbb{GT}_{N+1}^+$ respectively. $P_N$ and $P_{N+1}$ are $q$--coherent if and only if
$$
 {\mathcal S^*}(x_1,\dots,x_N;P_{N})={\mathcal S^*}(x_1,\dots,x_N,0;P_{N+1}).
$$
\end{proposition}

Next, we want to show that convergence of $q$--central probability measures is equivalent to
uniform convergence of their $q$--Schur generating functions.

Let $P_N$ and  $P^i_N,\, i=1,2,\dots$ be probability measures on $\mathbb{GT}_N$. We say that
$P^i_N$ \emph{weakly converges} to $P^N$ as $i\to\infty$ if $\lim_{i\to\infty} P^i_N(\mu)=P_N(\mu)$
for every $\mu\in\mathbb{GT}_N$.

Let $P$ and $P^i$ be probability measures on $\mathcal T$. We say that $P^i$ \emph{weakly
converges} to $P$ as $i\to\infty$ if $\lim_{i\to\infty} P^i(C_\tau) = P(C_\tau)$ for any cylinder
set $C_{\tau}$.

Suppose that $P$ and $P^i$ are $q$--central probability measures on $\mathcal T$ and let $P_N$,
$P^i_N$ be the corresponding $q$--coherent systems (i.e.\  projections of the measures on
$\mathbb{GT}_N$).

\begin{proposition}
\label{proposition_convergence_of_central_measures_and_coherent_systems}
 Measures $P^i$ weakly converge to $P$ if and only if $P_N^i$ weakly converge to $P_N$ for every
$N$.
\end{proposition}
\begin{proof}
 This proposition follows from the correspondence between $q$--central measures and $q$--coherent
systems. (See Proposition \ref{Proposition_coherence_from_centrality} and Proposition
\ref{Proposition_from_coherent_system_to_measure}).
\end{proof}

In what follows the sign $\rightrightarrows$ stays for the uniform convergence of functions on
various sets.

\begin{proposition}
\label{Proposition_convergene_of_measures_implies_uniform_convergence_of_function}
 Let $P_N$, $P^i_N$ be probability measures on $\mathbb{GT}_N$. If measures $P_N^i$ weakly converge to
$P_N$, then
$$
 {\mathcal S}(x_1,\dots,x_N;P_N^i)\rightrightarrows {\mathcal S}(x_1,\dots,x_N;P_N)
$$
uniformly on $T_N$. If  $P_N$ and $P^i_N$ are supported on $\mathbb{GT}_N^+$, then the convergence
is uniform on $D_N$.
\end{proposition}

\begin{proposition}
\label{proposition_Convergence_of_S_implies_weak_convergence} Let $P^i_N$ be probability measures
on $\mathbb{GT}_N$. If
$$
 {\mathcal S}(x_1,\dots,x_N;P_N^i)
$$
 converge uniformly on $T_N$ to a function $S(x_1,\dots,x_N)$, then there exists a probability
measure $P_N$ such that
$$
 S(x_1,\dots,x_N)={\mathcal S}(x_1,\dots,x_N;P_N)
$$
and $P_N^i$ weakly converge to $P_N$.
\end{proposition}

\begin{proposition}
\label{Proposition_convergence_of_measures_implies_uniform_convergence_of_interpolation_generating_functions}
 Let $P_N$, $P^i_N$ be probability measures on $\mathbb{GT}^+_N$. $P_N^i$ weakly converge to $P_N$
if and only if
$$
 {\mathcal S^*}(x_1,\dots,x_N;P_N^i)\rightrightarrows {\mathcal S^*}(x_1,\dots,x_N;P_N)
$$
uniformly on compact subsets of $\mathbb C^N$.
\end{proposition}

The proofs are quite involved and we give them in Section
\ref{Subsection_proofs_probability_gen_functions}.

\section{The boundary of the $q$--Gelfand--Tsetlin graph }

Suppose that $\lambda\in\mathbb{GT}_N$. Let $P_i^\lambda$, $i=1,\dots,N$ be probability measures on
$\mathbb{GT}_i$ such that $P_1^\lambda,\dots,P_N^\lambda$ is a $q$--coherent system and
$P_N^\lambda$ is delta-measure on $\lambda$ (i.e.\  $P_N^{\lambda}(\lambda)=1$). Clearly, such
$q$--coherent system exists and is unique. We call $\{P_i^{\lambda}\}$ \emph{a primitive
$q$--coherent system} corresponding to $\lambda$.

%


Let $N(i)$ be an increasing sequence of positive integers and  $\lambda(i)\in\mathbb{GT}_{N(i)}$.
We call $\lambda(i)$ a \emph{regular} sequence if for every $k$ probability measures
$P^{\lambda(i)}_k$ weakly converge to a certain probability measure $P_k$. When we vary $k$, we see
that the limit measures $P_k$ form a $q$--coherent system. The set of all possible $q$--coherent
system that can be obtained in such a way is called \emph{the Martin boundary} of the graph. In
Sections \ref{subsection_Neces_cond_for_convergence}-\ref{subsection_properties_of_limit_measures}
we describe the Martin boundary of the $q$--Gelfand--Tsetlin graph. And in Section
\ref{subsection_proof_of_main_th} we prove that the minimal boundary of the $q$--Gelfand--Tsetlin
graph coincides with the Martin boundary of the $q$--Gelfand--Tsetlin graph.

\smallskip

In Section \ref{subsection_algebras} we defined the algebra of symmetric function $\Lambda$.
\emph{A specialization} ${\rm Spec}$ of $\Lambda$ is an arbitrary homomorphism from algebra
$\Lambda$ to $\mathbb C$. $H$--generating function of the specialization is the following power
series
$$
 H(t)=\sum_{i=0}^{\infty} {\rm Spec}(h_k) t^k.
$$
Since complete symmetric functions $h_k$ generate $\Lambda$, $H$--generating function uniquely
defines the specialization.

Recall that $\mathcal N$ is the set of all non-decreasing sequences of integers
$\nu={(\nu_1\le\nu_2\le\dots)}$. We prove the following theorem in Sections
\ref{subsection_Neces_cond_for_convergence}-\ref{subsection_computation_of_the_limits}.

\begin{theorem}
\label{Theorem_limits_of_measures} The Martin boundary of the $q$--Gelfand--Tsetlin graph
$\{\mathcal E^k_\nu\}$ is parameterized by points of $\mathcal N$.

 If a sequence of signatures $\lambda(i)$ is regular, then there exists $\nu\in\mathcal N$ such
that the last coordinates of $\lambda(i)$ stabilize to $\nu$, i.e.\  for any $j$
$$
 \lim_{i\to\infty} \lambda(i)_{N(i)+1-j}=\nu_j.
$$

If the last coordinates of $\lambda(i)$ stabilize to $\nu$, then
 the measures $P^{\lambda(i)}_k$ weakly tend to $\mathcal E^\nu_k$.

The $q$--Schur generating function of the measure $\mathcal E^k_\nu$ can be uniquely extended to a
function analytic everywhere in $(\mathbb C^*)^k$.

If $\nu'=\nu+\ell$, i.e.\  for every $j$ we have $\nu'_j=\nu_j+\ell$, then
$$
 {\mathcal S}(x_1,\dots,x_k;\mathcal E^k_{\nu'})=\frac{x_1^\ell\cdots
x_k^\ell}{1^{\ell}q^{-\ell}\cdots q^{(1-k)\ell}} {\mathcal S}(x_1,\dots,x_k;\mathcal E^k_{\nu'})
$$

If $\nu_1\ge 0$, then
$$ {\mathcal S}(x_1,\dots,x_k;\mathcal E^k_\nu)=\sum_{\mu\in\mathbb{GT}_k^+}(-1)^{|\mu|}q^{n(\mu)-n(\mu')}{\rm
Spec}_{\nu}(s_\lambda){s^*_{\mu}(x_1,\dots,x_k;q)},
$$
the series converges everywhere in $\mathbb C^k$ and ${\rm Spec}_{\nu}$ is a specialization of
algebra $\Lambda$ with $H$--generating function
$$
H^\nu(t)= \frac{\prod\limits_{i\ge 0}(1-q^it)}{\prod\limits_{j=1}^{\infty}(1-q^{\nu_j+j-1}t)}
$$
\end{theorem}

 It is worth stressing that Propositions \ref{Proposition_convergene_of_measures_implies_uniform_convergence_of_function} and
\ref{proposition_Convergence_of_S_implies_weak_convergence} imply that $\lambda(i)$ is regular if
and only if for any $k$ functions
$$
\frac{s_{\lambda(i)}(x_1,x_2,\dots,x_k,q^{-k},q^{-k-1},\dots,q^{1-N(i)})}{s_\lambda(i)(1,q^{-1},\dots,q^{1-N(i)})}
$$
converge uniformly on $T_k$. Thus, Theorems \ref{theorem_limit_of_polynomials} and
\ref{Theorem_limits_of_measures} are equivalent.

\subsection{Simple necessary conditions for convergence}
\label{subsection_Neces_cond_for_convergence}
  Recall that $A_\ell$ is an automorphisms of graph $\mathbb{GT}$ acting on an arbitrary signature
$\lambda$ by
$$
 A_\ell(\lambda_1\ge\dots\ge\lambda_N)=(\lambda_1+\ell)\ge\dots\ge(\lambda_N+\ell).
$$
We call $A_\ell$ an \emph{$\ell$--shift}.

\begin{lemma} If $P_1,P_2,\dots$ is a $q$--coherent system, then $A_\ell(P_1),A_\ell(P_2),\dots$ is a
$q$--coherent system too.
\end{lemma}
\begin{proof}
 This follows form the fact that $P(\lambda\to\mu)=P(A_\ell(\lambda)\to A_\ell(\mu))$, which is
straightforward.
\end{proof}

\begin{proposition}
 \label{Proposition_bound_on_P_1^lambda} For any $N>0$ and any
$\lambda=(\lambda_1\ge\dots\ge\lambda_N)\in\mathbb{GT}_N$ we have
 $$
  P_1^\lambda(\lambda_N)\ge\prod_{i=1}^{\infty}(1-q^i).
 $$
 Here $\lambda_N$ stands for a signature of $\mathbb{GT}_1$ with coordinate $\lambda_N$.
\end{proposition}
\begin{proof}
 Note that
 $$
   P_1^\lambda(\lambda_N)=P_1^{A_\ell(\lambda)}(\lambda_N+\ell).
 $$
 Thus, without loss of generality we may assume that $\lambda_N=0$.

 Since $P_N^\lambda$ is concentrated on the signature $\lambda$, we have
 $$
  {\mathcal S} (x_1,\dots,x_N; P_N^{\lambda})
=\frac{s_\lambda(x_1,\dots,x_N)}{s_\lambda(1,\dots,q^{1-N})}.
 $$
 Using Proposition \ref{Proposition_coherency_in_terms_of_Schur_polynomials} we conclude that
 $$
  {\mathcal S} (x;
P_1^{\lambda})=\frac{s_\lambda(x,q^{-1},\dots,q^{1-N})}{s_\lambda(1,\dots,q^{1-N})}.
 $$
 Recall that
 $$
  {\mathcal S} (x; P_1^{\lambda})=\sum_\ell P_1^{\lambda}(\ell) x^\ell.
 $$
 Since ${\mathcal S} (x_1,\dots,x_N; P_N^{\lambda})$ is a polynomial, so is  ${\mathcal S} (x;
P_1^{\lambda})$. It follows that
$$
P_1^{\lambda}(\lambda_N)=P_1^\lambda(0)= {\mathcal S} (0;
P_1^{\lambda})=\frac{s_\lambda(0,q^{-1},\dots,q^{1-N})}{s_\lambda(1,\dots,q^{1-N})}=\frac{s_\lambda(q^{-1},\dots,q^{1-N})}{s_\lambda(1,\dots,q^{1-N})}
$$

 The value $s_\lambda(1,q^{-1},\dots,q^{1-N})$ can be computed using the definition of the
 rational Schur function (see e.g.\ \cite[Example 3.1]{Mac}). Recall that $n(\lambda)=\sum
 (i-1)\lambda_i$. We have
\begin{equation}
\label{eq_explicit_formula_for_specialized_Schur_function}
 s_\lambda(1,q,\dots,q^{1-N})=q^{-n(\lambda)}\prod_{1\le i<j \le N}
 \frac{1-q^{\lambda_j-\lambda_i+i-j}}{1-q^{i-j}}
\end{equation}
It follows that
\begin{multline*}
 \frac{s_\lambda(q^{-1},\dots,q^{1-N})}{s_\lambda(1,\dots,q^{1-N})}=\dfrac{q^{-|\lambda|}\prod\limits_{1\le
 i<j \le N-1} \dfrac{1-q^{-\lambda_i+\lambda_j+i-j}}{1-q^{i-j}}}{\prod\limits_{1\le i<j \le N}
 \dfrac{1-q^{-\lambda_i+\lambda_j+i-j}}{1-q^{i-j}}}
 =\dfrac{q^{-|\lambda|}}{\prod\limits_{i=1}^{N-1}\dfrac{1-q^{-\lambda_i+i-N}}{1-q^{i-N}}}\\=
 \prod\limits_{i=1}^{N-1}\dfrac{q^{-\lambda_i}-q^{-\lambda_i+i-N}}{1-q^{-\lambda_i+i-N}}
 =\prod_{i=1}^{N-1}\dfrac{1-q^{N-i}}{1-q^{\lambda_i-i+N}}\ge\prod_{i=1}^{N-1}(1-q^{N-i})\ge\prod_{i=1}^{\infty}{(1-q^i)}.
\end{multline*}
 Thus,
$$
 P_1^{\lambda}(\lambda_N)\ge\prod_{i=1}^{\infty}(1-q^i).
$$
\end{proof}

\begin{corollary}
 \label{Corollary_regular_sequence_is_bounded} If a sequence of measures $P_1^{\lambda(i)}$ with
 $\lambda(i)\in\mathbb{GT}_{N_i}$, $N_i\to\infty$, weakly converges to a certain probability
 measure $P_1$, then $\lambda(i)_{N_i}$ is  bounded from below.
\end{corollary}
\begin{proof}
 Let $k$ be an integer such that $$P_1(\{\lambda\in\mathbb{GT}_1:\,
-k<\lambda<k\})>1-\frac12\prod_{j=1}^{\infty}(1-q^j).$$ Since $P_1^{\lambda(i)}$ weakly converges
to $P_1$ and the set $\{\lambda\in\mathbb{GT}_1:\, -k<\lambda<k\}$ is finite, we conclude that
$$P_1^{\lambda(i)}(\{\lambda\in\mathbb{GT}_1:\, -k<\lambda<k\})>1-\prod_{j=1}^{\infty}(1-q^j)$$ for
$i>i_0$. Thus, if $\lambda(i)_{N(i)}<{-k}$, then
$P_1^{\lambda(i)}(\lambda_{N(i)})<\prod_{j=1}^{\infty}(1-q^j)$. This is a contradiction with
Proposition \ref{Proposition_bound_on_P_1^lambda}, consequently, $\lambda(i)_{N(i)}\ge (-k)$ for
$i\ge i_0$.
\end{proof}

Corollary \ref{Corollary_regular_sequence_is_bounded} implies that for the full description of the
Martin boundary of the $q$--Gelfand--Tsetlin graph it is enough to study only measures concentrated
on positive signatures and their $\ell$--shifts.

\begin{proposition}
 Let $k<N$, $\lambda\in\mathbb{GT}_N$ and $\mu\in\mathbb{GT}_k$. If $P_k^\lambda(\mu)>0$, then
$\mu\ge(\lambda_{N-k+1},\lambda_{N-k+2},\dots,\lambda_N)$.
\end{proposition}
\begin{proof}
 From the definitions of measures $P_k^\lambda$ and $q$--coherent systems it follows that if
$P_k^\lambda(\mu)>0$, then there exists a sequence $\tau(k)\prec\tau(k+1)\prec\dots\prec\tau(N)$
such that $\tau(k)=\mu$ and $\tau(N)=\lambda$. Consequently, for $i=0,1,\dots,k-1$ we have
$$
 \mu_{k-i}=\tau(k)_{k-i}\le\tau(k+1)_{k+1-i}\le\dots\le \tau(N)_{N-i}=\lambda_{N-i}.
$$
\end{proof}

\begin{corollary}
\label{Corollary_regular_sequence_partially_bounded_from_above}
 If $\lambda(i)\in\mathbb{GT}_{N(i)}$ is a sequence of signatures such that the measures
$P_k^{\lambda(i)}$ weakly converge, then the sequence of integers $\lambda(i)_{N(i)-m}$ is bounded
from above for any $0\le m\le k-1$.
\end{corollary}

\subsection{Proof of Theorems \ref{theorem_limit_of_polynomials} and \ref{Theorem_limits_of_measures}}

\label{subsection_computation_of_the_limits}

We start with the following compactness result

\begin{proposition}
\label{proposition_tightness} Let $\lambda(i)$ be a sequence of signatures stabilizing to $\nu$.
The family of functions
$$
g_i(x_1,\dots,x_k)=\frac{s_{\lambda(i)}(x_1,x_2,\dots,x_k,q^{-k},q^{-k-1},\dots,q^{1-N(i)})}{s_\lambda(i)(1,q,\dots,q^{1-N(i)})}
$$
is a relatively compact subset of the set of continuous functions on $k$--dimensional torus $T_k$
with uniform convergence topology.
\end{proposition}
The proof is quite technical and we present it in Section \ref{Subsection_proofs_tightness}. The
main idea is to find a uniform estimate for the derivatives of functions $g_i(x_1,\dots,x_k)$.

\smallskip

In order to identify all the possible limits of the functions $g_i$, we want to decompose the
functions in $q$--interpolation polynomials series. We need the following proposition that is
proved in Section \ref{Subsection_proofs_spaces_of_functions}. Recall that
$q^{\lambda-\delta}=(q^{\lambda_1},q^{\lambda_2-1},\dots,q^{\lambda_N-N+1})$.

\begin{proposition}
\label{Proposition_all_properties_of_interpolation_decomposition}
 Let $P_N$ be a probability measure on $\mathbb{GT}_N$ such that ${\rm supp}(P_N)\subset \mathbb{
GT}_N^+$, in other words  the probability of any signature with at least one negative coordinate
is zero. The $q$--Schur generating function of $P_N$ is well-defined for all $(x_1,\dots,x_N)\in
D_N$ and it can be uniquely represented as
\begin{equation}
\label{eq_x9}
 {\mathcal S}(x_1,\dots,x_N; P_N)=\sum_{\mu\in\mathbb{GT}_N^+} c_\mu s^*_\mu(x_1,\dots,x_N; q).
\end{equation}
 The series converges uniformly on any ball $B(0,r)$ (in the usual Euclidian metric)  with radius $0<r<1$ and in every point
$q^{\lambda-\delta}$.

Furthermore, suppose that $P_N^i$ and $P_N$ are probability measures on $\mathbb{GT}_N$ such that
${\rm supp}(P_N^i)\subset \mathbb{GT}_N^+$ and ${\rm supp}(P_N)\subset \mathbb{GT}_N^+$. Let
$c_\mu^i$ and $c_\mu$ be the coefficients of the decomposition \eqref{eq_x9} for the $q$--Schur
generating functions of $P_N^i$ and $P_N$, respectively. If $P_N^i$ weakly converge to $P_N$ as
$i\to\infty$ (in other words, if the $q$--Schur generating functions of the measures uniformly
converge), then for every $\mu$
$$
 \lim_{i\to\infty} c^i_\mu=c_\mu.
$$
\end{proposition}
{\bf Remark.} Note that here we use functions $s^*_\mu(\cdot; q)$ while in the definition of a
$q$--interpolation Schur generating function we use $s^*_\mu(\cdot; q^{-1})$.

Let $\nu=(0\le\nu_1\le\nu_2\le \dots)$ be an arbitrary nondecreasing sequence of non-negative
integers.

Recall that $H^\nu(t)$ is the following function:
$$
H^\nu(t)=\frac{\prod\limits_{i\ge 0}(1-q^it)}{\prod\limits_{j=1}^{\infty}(1-q^{\nu_j+j-1}t)},
$$
$H^\nu(t)$ is an analytic function in $\mathbb{C}$.

There is a one-to-one correspondence $X$ between $\nu$s and subsets of $\mathbb Z_{\ge 0}$ with
infinite complement given by
$$
 X(\nu)=\mathbb Z_{\ge 0}\setminus \{\nu_i+i-1\}_{i=1,2,\dots}.
$$
Observe that
$$
H^{\nu}(t)=\prod_{x\in X(\nu)} (1-q^xt).
$$

The main part of the proof of Theorem \ref{theorem_limit_of_polynomials} is contained in the
following proposition.

\begin{proposition}
\label{proposition_necessary_conditions}
 Let $N(i)$ be an increasing sequence of positive integers and $\lambda(i)\in\mathbb{GT}^+_{N(i)}$.
Suppose that $P_k^{\lambda(i)}$ weakly converges as $i\to\infty$ to a certain probability measure
$P_k$ for $k=1,2,\dots$. Then the last coordinates of $\lambda(i)$ stabilize, i.e.\  there exists a
nondecreasing sequence of integers $\nu_j$ such that for every $j=1,2,\dots$
$$
 \lambda(i)_{N(i)+1-j}\to \nu_j.
$$
The $q$--Schur generating function of the limit measure $P_k=\mathcal E^{\nu}_k$ has the following
decomposition:

$$
 {\mathcal S}(x_1,\dots,x_k; \mathcal
E^{\nu}_k)=\sum_{\mu\in\mathbb{GT}_k^+}(-1)^{|\mu|}q^{n(\mu)-n(\mu')}{\rm
Spec}_{\nu}(s_\lambda){s^*_{\mu}(x_1,\dots,x_k;q)},
$$
where ${\rm Spec}_{\nu}$ is a specialization of algebra $\Lambda$ with $H$--generating function
$H^\nu(t)$
\end{proposition}
\begin{proof}
The probability measure $P^{\lambda(i)}_{N(i)}$ is concentrated on a single signature $\lambda(i)$,
thus
$$
 {\mathcal S}\left(x_1,\dots,x_{N(i)};
P^{\lambda(i)}_{N(i)}\right)=\frac{s_{\lambda(i)}(x_1,\dots,x_{N(i)})}{s_{\lambda(i)}(1,\dots,q^{-N(i)})}.
$$
Measures $P_1^{\lambda(i)},\dots,P_{N(i)}^{\lambda(i)}$ form a $q$--coherent system. Consequently,
Proposition \ref{Proposition_coherency_in_terms_of_Schur_polynomials} yields
$$
{\mathcal S}\left(x_1,\dots,x_k;P^{\lambda(i)}_k\right)=
\frac{s_{\lambda(i)}(x_1,x_2,\dots,x_k,q^{-k},q^{-k-1},\dots,q^{1-N(i)})}{s_\lambda(i)(1,q^{-1},\dots,q^{1-N(i)})},
$$

Let us expand the right side of the last formula into the sum of $q$-interpolation Schur
polynomials. Coefficients of this expansion are given by the formula
\eqref{eq_schur_through_interpolation}:

\begin{multline}
 \frac{s_{\lambda(i)}(x_1,x_2,\dots,x_k,q^{-k},\dots, q^{1-N(i)})} {s_{\lambda(i)}(1,q^{-1},\dots,
q^{1-N(i)})}\\=\sum_{\mu\in\mathbb{GT}^+_{N(i)}}\frac{q^{(N(i)-1)|\mu|}s^*_\mu(q^{\lambda(i)-\delta};q)}{s^*_\mu(q^{\mu-\delta};q)}
\frac{s^*_{\mu}(x_1,\dots,x_k,q^{-k},\dots,q^{1-N(i)})}{s_{\mu}(1,q^{-1},\dots,q^{N(i)-1})}
\end{multline}
The combinatorial formula for interpolation Schur polynomials (see Proposition
\ref{Proposition_combinatoral_for_interp_Schur}) implies that if $\mu_{k+1}\ne 0$, then
$$
s^*_{\mu}(x_1,\dots,x_k,q^{-k},\dots,q^{1-N(i)})=0.
$$
Furthermore, polynomials $s^*_\mu$ are stable, i.e.\  if $\mu_{k+1}=\mu_{k+2}=\dots=0$, then
$$
s^*_{\mu}(x_1,\dots,x_k,q^{-k},\dots,q^{1-N(i)})=s^*_{(\mu_1,\dots,\mu_k)}(x_1,\dots,x_k).
$$
It follows that
\begin{multline}
\label{eq_decomposition_interp_Schur_prelimit}
 \frac{s_{\lambda(i)}(x_1,x_2,\dots,x_k,q^{-k},\dots, q^{1-N(i)})} {s_{\lambda(i)}(1,q^{-1},\dots,
q^{1-N(i)})}\\=\sum_{\mu\in\mathbb{GT}^+_k}\frac{q^{(N(i)-1)|\mu|}s^*_\mu(q^{\lambda(i)-\delta};q)}{s^*_\mu(q^{\mu-\delta};q)}
\frac{s^*_{\mu}(x_1,\dots,x_k)}{s_{\mu}(1,q,\dots,q^{N(i)-1})},
\end{multline}

Using Proposition \ref{Proposition_convergene_of_measures_implies_uniform_convergence_of_function}
and Proposition \ref{Proposition_all_properties_of_interpolation_decomposition} we conclude that
weak convergence of measures $P_k^{\lambda(i)}$ implies that for any $\mu$
$$
\frac{q^{(N(i)-1)|\mu|}s^*_\mu(q^{\lambda(i)-\delta};q)}{s^*_\mu(q^{\mu-\delta};q){s_{\mu}(1,q,\dots,q^{N(i)-1})}}
$$
converges as $i\to\infty$.

Using again \cite[Example 3.1]{Mac} we obtain
$$
\lim_{N\to\infty} s_\mu(1,q,\dots,q^N)=\lim_{N\to\infty}
\frac{q^{n(\mu)}\prod\limits_{(i,j)\in\mu}(1-q^{N+j-i})}{\prod\limits_{(i,j)\in\mu}(1-q^{h(i,j)})}
= \frac{q^{n(\mu)}}{\prod\limits_{(i,j)\in\mu}(1-q^{h(i,j)})},
$$
Proposition \ref{proposition_interpolation_property} yields
$$
 s^*_\mu(q^{\mu};q)=q^{n(\mu')-2n(\mu)}\prod\limits_{(i,j) \in \mu}\left(q^{h(i,j)}-1\right).
$$
Thus,
$$
\frac1{s^*_\mu(q^{\mu};q)s_{\mu}(1,q,\dots)}\to (-1)^{|\mu|}q^{n(\mu)-n(\mu')}
$$
Consequently,
\begin{equation}
\label{eq_prelimit_s_star}
 q^{(N(i)-1)|\mu|}s^*_\mu(q^{\lambda(i)-\delta};q)
\end{equation}
should tend to a limit as $i\to\infty$.

Recall that functions $s^*_\mu$ form a linear basis of filtered algebra $\widehat \Lambda$.
Therefore, \eqref{eq_prelimit_s_star} has a limit for any $\mu$ if and only if
 $$
 \lim_{i\to\infty} q^{(N(i)-1)k}f(q^{\lambda(i)-\delta})
 $$
 for any $f\in\widehat \Lambda$ of the degree $k$.

 Let us introduce an analogue of Newton power sums in algebra $\widehat \Lambda$:
 $$
  p^*_k=\sum_{j\ge 1}\left( x_j^k-(q^{1-j})^k\right)
 $$
 (with $p^*_0=1$). Note that functions
 $$
  p^*_{k_1}\cdots p^*_{k_l}
 $$
 also form a linear basis of $\widehat \Lambda$,

 We conclude that \eqref{eq_prelimit_s_star} has a limit for every $\mu$ if and only if
 $$
  q^{(N(i)-1)(k_1+\dots+k_l)}p^*_{k_1}(q^{\lambda(i)-\delta})\cdots
p^*_{k_1}(q^{\lambda(i)-\delta})
 $$
 converges for any positive integers $k_1,\dots,k_l$. The last limits exist if and only if
 $$
  q^{(N(i)-1)k}p^*_k(q^{\lambda(i)-\delta})
 $$
  tends to a finite limit as $i\to\infty$ (for any $k$). We claim that then the last coordinates of $
  \lambda(i)$ should stabilize.

 Indeed, we have
 $$
  q^{(N(i)-1)k}p^*_k (q^{\lambda(i)-\delta})=\sum_{j=1}^{N(i)}
(q^{k\lambda(i)_j}-1)q^{k(N(i)-i)}=\sum_{j=0}^{N(i)}q^{kj}(q^{k\lambda(i)_{N(i)-j}}-1).
 $$

Denote $f(i,j)=\lambda(i)_{N(i)-j}+j$. Note that $\sum_{j=0}^{N(i)} q^{kj}(-1)$ converges for any
$k$. Thus,
$$
 \sum_{j=0}^{N(i)}q^{kj}q^{k\lambda(i)_{N(i)-j}}=\sum_{j=0}^{N(i)}q^{kf_{i,j}}
$$
should also converge. But since $f(i,j)$ is increasing in $j$, we have
$$
\left|\sum_{j=0}^{N(i)}q^{kf(i,j)}-q^{kf(i,0)}\right|<q^{k(f(i,0)+1)}(1+q+q^2+\dots)=q^{kf(i,0)}\frac{q^k}{1-q}.
$$
Informally speaking, if $k$ is large enough, then
$$
 \sum_{j=0}^{N(i)}q^{kf(i,j)}\approx q^{kf(i,0)}.
$$
Recall that $f(i,0)$ is an integer. Therefore, if $ \sum_{j=0}^{N(i)}q^{kf(i,j)}$ converges, then
either $f(i,0)\to+\infty$ or $f(i,0)$ stabilize as $i\to\infty$. Repeating the same argument for
$\sum_{j=w}^{N(i)}q^{kf_{i,j}}$, $w=1,2,\dots$ we conclude that for every $w$, either
$f(i,w)\to+\infty$ or $f(i,w)$ stabilize as $i\to\infty$. In the former case $\lambda(i)_{N(i)-w}$
is unbounded which is a contradiction with Corollary
\ref{Corollary_regular_sequence_partially_bounded_from_above}. Thus, $f(i,w)$ stabilize, in other
words there exists a sequence $0\le\nu_1\le\nu_2\le\dots$ such that
$$
 \lambda(i)_{N(i)+1-w}\to \nu_w.
$$
In this case
$$
\lim_{i\to\infty}  q^{(N(i)-1)k}p^*_k (q^{\lambda(i)-\delta};q)=\sum_{i\ge 1}
(q^{k\nu_i}-1)q^{k(i-1)}.
$$

For any Young diagram (i.e.\  positive signature) let $s_\lambda$ denote the element of $\Lambda$
corresponding to Schur function $s_\lambda(x_1,\dots,x_N)$. Newton power sums $p_k$ are
algebraically independent generators of $\Lambda$. Thus, for any $\lambda$ there exists a unique
polynomial $R_\lambda$, such that
$$
 s_\lambda=R_\lambda(p_1,\dots,p_m).
$$
(Here $m$ also depends on $\lambda$, but we omit this dependence to simplify the notations.)

Now consider the following element of $\widehat \Lambda$:
$$
 r^*_\lambda=s^*_\lambda-R_{\lambda}(p_1^*,\dots,p_m^*).
$$
Note that, if we work with finite sets of variables (i.e.\  in algebras $\Lambda_N$ and $\widehat
\Lambda_N$) then the highest homogenous component of $s^*_\lambda$ is exactly $s_\lambda$, and the
same is true for $p^*_k$ and $p_k$. It follows that the degree of $r^*_\lambda$ is less than
$|\lambda|$. Consequently,
$$
\lim_{i\to\infty} q^{(N(i)-1)|\lambda|} r^*_\lambda(q^{\lambda(i)-\delta})=0.
$$

Let ${\rm Spec}_{\nu}$ be a specialization of $\Lambda$ defined on generators $p_k$ through
$$
 {\rm Spec}_{\nu}(p_k)=\sum_{i\ge 1} (q^{k\nu_i}-1)q^{k(i-1)}.
$$

The above arguments show that
\begin{multline*}
 \lim_{i\to\infty} q^{(N(i)-1)|\lambda|}s^*_\lambda=\lim_{i\to\infty}
q^{(N(i)-1)|\lambda|}R_{\lambda}(p_1^*,\dots,p_m^*)\\=R_{\lambda}({\rm Spec}_{\nu}(p_1),\dots,{\rm
Spec}_{\nu}(p_m))={\rm Spec}_{\nu}(s_\lambda).
\end{multline*}

Thus, we have proved that if $P_k^{\lambda(i)}$ weakly converges as $i\to\infty$ to a certain
probability measure $P_k$ for any $k=1,2,\dots$, then  $\lambda(i)_{N(i)+1-j}\to \nu_j$ and
\begin{multline}
\label{eq_Limit_fucntion}
 {\mathcal S}(x_1,\dots,x_k;P_k)=\lim_{i\to\infty} {\mathcal
S}\left(x_1,\dots,x_k;P_k^{\lambda(i)}\right)\\=
\sum_{\mu\in\mathbb{GT}_k^+}(-1)^{|\mu|}q^{n(\mu)-n(\mu')}{\rm
Spec}_{\nu}(s_\mu){s^*_{\mu}(x_1,\dots,x_k)}
\end{multline}

It remains to prove that $H$--generating function of ${\rm Spec}_{\nu}$ is $H^\nu(t)$.

$P$--\emph{generating function} of ${\rm Spec}_{\nu}$ is the following power series
$$
 P(t)=\sum_{k=1}^{\infty} {\rm Spec}_{\nu}(p_k)  t^{k-1}
$$

The following equality relates  $P$--generating function to $H$--generating function:
$P(t)=H'(t)/H(t)$. See e.g. \cite[2.10]{Mac} for the proof. The following computation completes the
proof:

\begin{multline*}
 P(t)=\sum_{k=1}^{\infty} {\rm Spec}_\nu (p_k) t^{k-1}=\sum_{k=1}^{\infty}\sum_{i=1}^{\infty}
 (q^{k\nu_i}-1)q^{k(i-1)}t^{k-1}
\\=\sum_{i=1}^{\infty}
 \frac{q^{\nu_i+i-1}}{1-q^{\nu_i+i-1}t}-\sum_{i=1}^{\infty}\frac{q^{i-1}}{1-q^{i-1}t},
\end{multline*}

\begin{multline}
\label{eq_H_gen_1formula} H(t)=\sum_{k=0}^{\infty} {\rm Spec}_\nu (h_k) t^{k}=\exp\left(\int {\rm
Spec}_\nu(P(t))\right)\\=\exp\left(\sum_{i\ge 1} -\ln(1-q^{\nu_i+i-1}t)+\sum_{i\ge
1}\ln(1-q^{i-1}t)\right)\\=\frac{\prod\limits_{i\ge
0}(1-q^it)}{\prod\limits_{j=1}^{\infty}(1-q^{\nu_j+j-1}t)}=H^\nu(t).
\end{multline}

\end{proof}

We also need yet another technical proposition that will be proved in Section
\ref{Subsection_proofs_analyticity}:
\begin{proposition}
\label{Proposition_analyticity} For every $\nu$ the series $$
\sum_{\mu\in\mathbb{GT}_k^+}(-1)^{|\mu|}q^{n(\mu)-n(\mu')}{\rm
Spec}_{\nu}(s_\mu){s^*_{\mu}(x_1,\dots,x_k)}$$ converges for all $x_1,\dots x_k$ and defines an
entire function in $\mathbb C^k$.
\end{proposition}

Now we are ready to finish the proof of Theorems \ref{theorem_limit_of_polynomials} and
\ref{Theorem_limits_of_measures}.

\begin{proof}[Proof of Theorem \ref{theorem_limit_of_polynomials}]
Let $\lambda(i)\in\mathbb{GT}_{N(i)}$ be a regular sequence of signatures. It means that for any
$k$ probability measures $P_k^{\lambda(i)}$ weakly converge to a certain probability measure
$P_k$. Corollary \ref{Corollary_regular_sequence_is_bounded} implies that there exists $\ell$ such
that $\lambda(i)_j\ge\ell$ for every $i,j$. Let $\mu(i)=A_{-\ell}(\lambda(i))$, then
$\mu(i)\in\mathbb{GT}_{N(i)}^+$ and for any $k$ probability measures $P_k^{\mu(i)}$ weakly
converge to a measure $\widetilde P_k$. Applying Proposition
\ref{proposition_necessary_conditions} we conclude that the last coordinates of $\mu(i)$ stabilize
to some $\nu_j$. Thus, the last coordinates of $\lambda(i)$ also stabilize
$$
 \lim_{i\to\infty} \lambda(i)_{N(i)+1-j}\to \nu_j+\ell.
$$
If all the coordinates of $\lambda(i)$ are non-negative starting from some $i$ (equivalently, if
the last coordinates stabilize to non-negative numbers), then Proposition
\ref{proposition_necessary_conditions} and Proposition \ref{Proposition_analyticity} imply that
$q$--Schur generating function of $P_k$ is an analytic function with desired interpolation
polynomials series decomposition. Again using the correspondence between weak convergence of
probability measures and convergence of their $q$--Schur generating functions (Proposition
\ref{proposition_Convergence_of_S_implies_weak_convergence} and Proposition
\ref{Proposition_convergene_of_measures_implies_uniform_convergence_of_function})  we conclude
that polynomials of Theorem \ref{theorem_limit_of_polynomials} converge to the desired analytic
function. For general $\lambda(i)$ the limit function is a product of the analytic $q$--Schur
generating function corresponding to measure $\mathcal E_k^\nu$ and polynomial
$$
 \frac{x_1^\ell\cdots x_k^\ell}{1^\ell q^{-\ell}\cdots q^{-(k-1)\ell}},
$$
consequently the limit function is analytic.

\smallskip
Now suppose that a sequence of signatures $\lambda(i)\in\mathbb{GT}_{N(i)}$, $N(i)\to\infty$ is
such that
$$
 \lim_{i\to\infty} \lambda(i)_{N(i)+1-j}\to \nu_j
$$
for some sequence $\nu_1\le\nu_2,\dots$. Proposition \ref{proposition_tightness} yields that the
sequence of functions
$$
g_i(x_1,\dots,x_k)=\frac{s_{\lambda(i)}(x_1,x_2,\dots,x_k,q^{-k},q^{-k-1},\dots,q^{1-N(i)})}{s_\lambda(i)(1,q^{-1},\dots,q^{1-N(i)})}
$$
has converging subsequences. But in the first part of the theorem we have identified all the
possible subsequential limits and their are the same. Thus, $g_i(x_1,\dots,x_k)$ converges
uniformly on $T_k$.
\end{proof}

\subsection{Limit measures}

The only formula we have for $q$--Schur generating functions of $\mathcal E^\nu_k$ is an infinite
series expansion. The situation is different if we turn to $q$--interpolation Schur generating
functions.

\begin{proposition} We have
$$
{\mathcal S^*}(x_1,\dots,x_k; \mathcal E^\nu_k)=H^\nu(x_1)\cdots H^\nu(x_k).
$$
\label{Proposition_interpolation_gen_function_for_limit_measure}
\end{proposition}

Recall that $\mathcal S^* (x_1,\dots,x_k;P_k)$ is an entire function for any measure $P_k$ with
${\rm supp}(P_k)\subset \mathbb{GT}_k^+$. We need the following lemma:
\begin{lemma}
\label{lemma_G=G'} Let $P_k$ be a probability measure of $\mathbb{GT}_k^+$
 Suppose that $\mathcal S^* (x_1,\dots,x_k;P_k)$ has the following Taylor series expansion
 $$
  \mathcal S^* (x_1,\dots,x_k;P_k)=\sum_{\mu\in\mathbb {GT}^+_k}
  a_{\mu}\frac{s_\mu(x_1,\dots,x_k)}{s_\mu(1,q^{-1},\dots,q^{1-k})}
 $$
 Define
 \begin{equation}
 \label{eq_definition_of_G_extended} F=\sum_{\mu\in\mathbb {GT}^+_k}
 a_{\mu}\frac{s^*_\mu(x_1,x_2,\dots,x_k;q)}{s^*_\mu(0,\dots,0;q)}.
 \end{equation}
 Then  the series on the right side of \eqref{eq_definition_of_G_extended} converges uniformly on
 any ball $B(0,r)$ (in the usual Euclidian metric) with radius $0<r<1$ and in every point $q^{\lambda-\delta}$. Furthermore,
$$F=\mathcal S(x_1,\dots,x_k;P_k).$$
\end{lemma}
If the support of the measure $P_k$ is finite, then both $\mathcal S^* (x_1,\dots,x_k;P_k)$ and
$F$ are polynomials and Lemma \ref{lemma_G=G'} immediately follows from the symmetry between
formulas \eqref{eq_interpolation_through_schur} and \eqref{eq_schur_through_interpolation}. The
complete proof of the lemma is given in Section \ref{Subsection_proofs_spaces_of_functions}.

\begin{proof}[Proof of Proposition \ref{Proposition_interpolation_gen_function_for_limit_measure}]
Let us expand $\mathcal S^* (x_1,\dots,x_k;\mathcal E^\nu_k)$ in Taylor series:
$$
\mathcal S^* (x_1,\dots,x_k;\mathcal E^\nu_k)=\sum_{\mu\in\mathbb {GT}^+_k}
  a_{\mu}\frac{s_\mu(x_1,\dots,x_k)}{s_\mu(1,q^{-1},\dots,q^{1-k})}.
$$
Using Lemma \ref{lemma_G=G'} we conclude that
$$
\mathcal S (x_1,\dots,x_k;\mathcal E^\nu_k)=\sum_{\mu\in\mathbb {GT}^+_k}
  a_{\mu}\frac{s^*_\mu(x_1,x_2,\dots,x_k;q)}{s^*_\mu(0,\dots,0;q)}.
$$
Comparing with Theorem \ref{Theorem_limits_of_measures} we see that
$$
 a_\mu=(-1)^{|\mu|}q^{n(\mu)-n(\mu')}{\rm Spec}_{\mu}(s_\mu){s^*_\mu(0,\dots,0;q)}.
$$
The combinatorial formulas for polynomials $s_\lambda$ and $s^*_\lambda$ imply that
$$
(-1)^{|\mu|}q^{n(\mu)-n(\mu')}\frac{s^*_\mu(0,\dots,0;q)}{s_\mu(1,q^{-1},\dots,q^{1-k})}=1.
$$
Thus,
\begin{multline*}
\mathcal S^* (x_1,\dots,x_k;\mathcal E^\nu_k)\\=\sum_{\mu\in\mathbb {GT}^+_k}
  (-1)^{|\mu|}q^{n(\mu)-n(\mu')}{\rm
Spec}_{\mu}(s_\mu){s^*_\mu(0,\dots,0;q)}\frac{s_\mu(x_1,\dots,x_k)}{s_\mu(1,q^{-1},\dots,q^{1-k})}\\=
 \sum_{\mu\in\mathbb {GT}^+_k}{\rm Spec}_{\mu}(s_\mu)s_\mu(x_1,\dots,x_k).
\end{multline*}
It remains to prove that
\begin{equation}
\label{eq_x11}
 H^\nu(x_1)\cdots H^\nu(x_k)=\sum_{\mu\in\mathbb {GT}^+_k}{\rm
Spec}_{\nu}(s_\mu)s_\mu(x_1,\dots,x_k),
\end{equation}
 where
$$
 H^\nu(x)=\sum_{j=0}^{\infty} {\rm Spec}_{\nu}(h_j)x^j.
$$
The formula \eqref{eq_x11} is a particular case of the well-known Cauchy identity for symmetric
functions (see e.g. \cite[Chapter 1, Section 4]{Mac}) and we leave its proof to the reader.
\end{proof}

\subsection{Some properties of measures ${\mathcal E^\nu}$}
\label{subsection_properties_of_limit_measures} In this section we discuss various properties of
the measures $\mathcal E^\nu_k$.

 For any nondecreasing infinite sequence of integers $\nu=(\nu_1,\nu_2,\dots)$ set
$$
 A_\ell(\nu)=(\nu_1+\ell,\nu_2+\ell,\dots)
$$
\begin{proposition} For any nondecreasing sequence of nonnegative integers $\nu=(\nu_1,\nu_2,\dots)$
any any $\ell=1,2,\dots$ we have
$$
 A_\ell(\mathcal E^\nu)=\mathcal E^{A_\ell(\nu)}.
$$
\end{proposition}
\begin{proof}
 There exists a sequence of signatures $\lambda(i)$, such that $\mathcal E^\nu_k=\lim
P^{\lambda(i)}_k$ for any $k$. Using Theorem \ref{Theorem_limits_of_measures} we conclude that
$$
 A_\ell(\mathcal E^\nu_k)=\lim A_\ell P^{\lambda(i)}_k =\lim P^{A_\ell(\lambda_i)}_k=\mathcal
E^{A_\ell(\nu)}_k.
$$
\end{proof}

Let us introduce a partial order on $\mathcal N$. We write $\nu\le\nu'$ if $\nu_i\le\nu'_i$ for
every $i$.

\begin{proposition}
\label{Proposition_support_and_monotonicity_of_extremal_measures}
 \begin{enumerate}
\item For any signature $\mu\in\mathbb{GT}_k$, $\mathcal E^{\nu}_k(\mu)=0$ unless
$\mu\ge(\nu_k,\nu_{k-1},\dots,\nu_1)$. Furthermore, $\mathcal E^{\nu}_k( (\nu_k,\dots,\nu_1))>0$.
\item If $\nu'>\nu$, then
$$
 \mathcal E^{\nu}_k((\nu_k,\dots,\nu_1))> \mathcal E^{\nu'}_k((\nu_k,\dots,\nu_1)).
$$
\end{enumerate}
\end{proposition}
\begin{proof}
 Without loss of generality assume that $\nu_1\ge 0$. For $\mu$ not belonging to $\mathbb{GT}_k^+$,
we have $\mathcal E^{\nu}_k(\mu)=0$. For $\mu\in\mathbb{GT}_k^+$, $\mathcal E^{\nu}_k(\mu)$ is
found from the representation:
\begin{equation}
\label{eq_x2}
 H^{\nu}(x_1)\cdots H^{\nu}(x_k)=\sum_{\mu\in\mathbb{GT}_k^+} \mathcal E^{\nu}_k(\mu)
\frac{s_\mu^*(q^{k-1}x_1,\dots,q^{k-1}x_k;q^{-1})}{s_\mu^*(0,\dots,0;q^{-1})}.
\end{equation}
To find $\mathcal E^{\nu}_k(\mu)$ we substitute $x=q^{-(k-1+\lambda-\delta)}$ into \eqref{eq_x2}
and use interpolation property of polynomials $s^*$ (which is stated in Proposition
\ref{proposition_interpolation_property}). We start from the empty diagram $\lambda=\emptyset$ and
then add boxes to it and find the numbers   $\mathcal E^{\nu}_k(\mu)$ inductively. (We discuss this
procedure in more details later. See Proposition
\ref{proposition_how_to_obtain_interpolation_coefficients}.) It follows from the definition of
function $H^{\nu}(t)$ that
$$
 H^{\nu}(q^{-(\lambda_1+k-1)})\cdots H^{\nu}(q^{-(\lambda_k)})=0,
$$
unless $(\lambda_1,\dots,\lambda_k)\ge(\nu_k,\dots,\nu_1)$. Thus, $\mathcal E^{\nu}_k(\mu)=0$
unless ${\mu\ge(\nu_k,\nu_{k-1},\dots,\nu_1)}$.

Next, substitute $x=(q^{-(\nu_k+k-1)},\dots,q^{-\nu_1})$ in \eqref{eq_x2}. We obtain
\begin{multline*}
 H^{\nu}(q^{-(\nu_k+k-1)})\cdots
H^{\nu}(q^{-(\nu_1)})\\=\mathcal
E^{(\nu_k,\dots,\nu_1)}_k((\nu_k,\dots,\nu_1))\frac{s_{(\nu_k,\dots,\nu_1)}^*(q^{-\nu_k},\dots,q^{-\nu_1+k-1};q^{-1})}{s_{(\nu_k,\dots,\nu_1)}^*(0,\dots,0;q^{-1})}.
\end{multline*}
Observe that for $i=1,\dots,k$ we have $H^{\nu}(q^{-(\nu_i+i-1)})>0$. It follows that
$$
  \mathcal E^{\nu}_k((\nu_k,\dots,\nu_1))>0.
$$
Now, let us prove the second part of Proposition
\ref{Proposition_support_and_monotonicity_of_extremal_measures}. If for some $i\in\{1,\dots,k\}$,
$\nu'_i>\nu_i$, then
$$
 \mathcal E^{\nu'}_k((\nu_k,\dots,\nu_1))=0.
$$
Otherwise note that for $i=1,\dots,k$ we have
$H^{\nu}(q^{-(\nu_i+i-1)})>H^{\nu'}(q^{-(\nu_i+i-1)})$. In both cases
$$
  \mathcal E^{\nu}_k((\nu_k,\dots,\nu_1))> \mathcal E^{\nu'}_k((\nu_k,\dots,\nu_1)).
$$
\end{proof}

Actually we can say more, i.e.\  $\mathcal E^{\nu}_k( (\nu_k,\dots,\nu_1))>c>0$ for some constant
$c$ depending solely on $k$. Let us prove this fact for $k=1$.
\begin{lemma} We have
\label{lemma_value_of_measure_at_minimal_point}
$$
 \mathcal E^{\nu}_1( \nu_1)\ge\prod_{i=1}^{\infty}(1-q^i).
$$
\end{lemma}
\begin{proof}
 Assume without loss of generality that $\nu_1\ge0$. Repeating the argument of Proposition
\ref{Proposition_support_and_monotonicity_of_extremal_measures} we conclude that
\begin{multline*}
 \mathcal E^{\nu}_1(
\nu_1)=\dfrac{H^\nu(q^{-\nu_1})}{(1-q^{-\nu_1})\cdots(1-q^{-1})}\\=\dfrac{\prod\limits_{j\ge 0,\,
j\ne \nu_1}(1-q^{-\nu_1}q^{j})}{ \prod\limits_{j=2}^{\infty}(1-q^{-\nu_1}q^{\nu_j+j-1}) }
\dfrac1{(1-q^{-\nu_1})\cdots(1-q^{-1})}\\=
\dfrac{\prod\limits_{j=\nu_1+1}^{\infty}(1-q^{-\nu_1}q^{j})}{
\prod\limits_{j=2}^{\infty}(1-q^{-\nu_1}q^{\nu_j+j-1}) }\ge
\prod\limits_{j=\nu_1+1}^{\infty}(1-q^{-\nu_1}q^{j})=\prod_{i=1}^{\infty}(1-q^i)
\end{multline*}
\end{proof}

The set $\mathcal N\subset \mathbb{Z}^\infty$ has a natural topology as a subset of the direct
product of discrete spaces $\mathbb{Z}$. A sequence $\theta(i)$ converges to $\nu$ in this topology
if and only if for every $k$ there exist $i_0$ such that $\theta(i)_k=\nu_k$ for $i>i_0$.

\begin{proposition}
\label{proposition_topology_on_extreme_measures}
 Sequence of probability measures $\mathcal E^{\theta(i)}$ weakly converges to a probability measure $P$ if
and only if $P=\mathcal E^\nu$ for some $\nu$, and $\theta(i)\to\nu$.
\end{proposition}
\begin{proof}
 Suppose that $\theta(i)\to\nu$. Without loss of generality assume that $\nu_1\ge0$. Then for
$i>i_i$, $\theta(i)_1\ge 0$. Thus, for $i>i_1$ we may use Proposition
\ref{Proposition_interpolation_gen_function_for_limit_measure}. One proves that
$$
 H^{\theta(i)}(t)\rightrightarrows H^{\nu}(t)
$$
uniformly on compact subsets of $\mathbb C$. It follows that for any $k$
\begin{multline*}
 {\mathcal S^*}(x_1,\dots,x_k;\mathcal E^{\theta(i)})=H^{\theta(i)}(x_1)\cdots
H^{\theta(i)}(x_k)\rightrightarrows H^{\nu}(x_1)\cdots H^{\nu}(x_k)\\={\mathcal
S^*}(x_1,\dots,x_k;\mathcal E^{\nu}).
\end{multline*}
Using Proposition
\ref{Proposition_convergence_of_measures_implies_uniform_convergence_of_interpolation_generating_functions}
and Proposition \ref{proposition_convergence_of_central_measures_and_coherent_systems} we conclude
that measures $\mathcal E^{\theta(i)}$ weakly converge to $\mathcal E^{\nu}$.

\smallskip

Now suppose that $\mathcal E^{\theta(i)}$ is a weakly convergent sequence. Lemma
\ref{lemma_value_of_measure_at_minimal_point} implies that $\mathcal
E^{\theta(i)}_1(\theta(i)_1)>c$ for some constant $c>0$. Thus, $\theta(i)_1$ should be bounded
from below. (Otherwise, measures $\mathcal E^{\theta(i)}_1$ ``escape to infinity''.) Choose
$l\ge\max(-\theta(i)_1)$. Applying, if necessary, $A_\ell$ we may assume without loss of
generality that $\theta(i)_1\ge 0$ for all $i$. Then Proposition
\ref{Proposition_convergence_of_measures_implies_uniform_convergence_of_interpolation_generating_functions}
yields that $q$--interpolation Schur generating functions of measures $\mathcal E^{\theta(i)}_1$
converge as $i\to\infty$ to $q$--interpolation Schur generating function of measure $P$ uniformly
on compact subsets of $\mathbb C$. Recall that
$$
 {\mathcal S^*}(t;\mathcal E^{\theta(i)}_1)=H^{\theta(i)}(t)=\prod_{x\in X(\theta(i))}(1-q^x t).
$$
Due to Rouch\'{e}'s theorem, uniform convergence of analytic functions implies the convergence of
positions of their simple zeros. Since zeros of $H^{\theta(i)}(t)$ are precisely $q^{-x},\, x\in
X(\theta(i))$, we conclude that the sets $X(\theta(i))$ converge (in a sense that for any $n$,
$X(\theta(i))\cap[0,\dots,n]$ stabilizes as $i\to\infty$) to a set $\hat X$. It remains to prove
that $\hat X$ has infinite complement. Indeed, otherwise there exists $k$ such that
$\theta(i)_k\to\infty$. But $\mathcal E^{\theta(i)}_k$ is concentrated on signatures $\lambda$
such that $\lambda_i\ge\theta(i)_{k+1-i}$, thus, measures $\mathcal E^{\theta(i)}_k$ can not
converge to a probability measure on $\mathbb{GT}_k$.

Consequently, $\hat X= X(\nu)$ for a certain $\nu\in\mathcal N$. Therefore, $\mathcal
E^{\theta(i)}\to \mathcal E^{\nu}$.
\end{proof}

Next we prove that measures $E^{\nu}$ are linearly independent in the following sense:

\begin{proposition}
\label{Proposition_extremality}
 If for some $\theta\in\mathcal N$ and for some probability measure $\pi$ defined on
 $\sigma$--algebra of Borel sets in $\mathcal N$ we have
$$
 \mathcal E^\theta=\int_{\mathcal N} \mathcal E^\nu d\pi,
$$
i.e.\  for any cylinder set $C_\tau$
$$
 \mathcal E^\theta(C_\tau)=\int_{\mathcal N} \mathcal E^\nu(C_\tau) d\pi,
$$
then $\pi$ is a delta measure on $\theta$, in other words $\pi(\theta)=1$.
\end{proposition}
\begin{proof}
 Let us prove that $\pi(\{\nu:\, \nu\ge \theta\})=1$. Assume the opposite. Then there exist $k$ and
$\tilde\nu_1,\dots,\tilde\nu_k$, such that $\pi(\{\nu:\, \nu_i=\tilde\nu_i\text{ for
}i=1,\dots,k\})>0$ and $\nu_i<\theta_i$ for some $i\in\{1,\dots,k\}$. But we have
$$
 \mathcal E^\theta_k=\int_{\mathcal N} \mathcal E^\nu_k d\pi,
$$
in particular
\begin{equation}
\label{eq_x3}
 \mathcal E^\theta_k((\tilde\nu_1,\dots,\tilde\nu_k))=\int_{\mathcal N} \mathcal E^\nu_k
((\tilde\nu_k,\dots,\tilde\nu_1)) d\pi,
\end{equation}
 The first part of Proposition \ref{Proposition_support_and_monotonicity_of_extremal_measures}
implies that the left side of \eqref{eq_x3} vanishes, while the right side is positive. This
contradiction proves that $\pi$--almost surely $\nu\ge\theta$.

On the other hand
\begin{equation}
\label{eq_x4}
 \mathcal E^\theta_k((\theta_k,\dots,\theta_1))=\int_{\mathcal N} \mathcal E^\nu_k ((\theta_k,\dots,\theta_1)) d\pi,
\end{equation}
Using the second part of Proposition
\ref{Proposition_support_and_monotonicity_of_extremal_measures} we conclude that if $\pi(\{\nu:\,
\nu>\theta\})>0$ then the right side of \eqref{eq_x4} should be strictly less than the left side.
Thus, $\pi(\theta)=1$.
\end{proof}

\subsection{Proof of Theorem \ref{theorem_Main}}

\label{subsection_proof_of_main_th} We start with 3 propositions which hold not only for the for
the $q$--Gelfand--Tsetlin graph, but in a much larger  generality.

Let $\Omega_q$ denote the convex set of $q$--central probability measures on $\mathcal T$. There is
a natural topology on $\Omega_q$, i.e.\  a minimal topology such that for any cylinder set $C_\tau$
the map
$$
 O_{\tau}: \Omega_q\to \mathbb{R}, \quad O_\tau(P)=P(C_\tau)
$$
is continuous. Convergence in this topology coincides with weak convergence of probability
measures.

Recall that the minimal boundary of the $q$--Gelfand-Tsetlin graph ${\rm Ex\,}\Omega_q$  is the set
of extremal points of $\Omega_q$. We denote elements of ${\rm Ex\,}\Omega_q$ by $\omega$.
\begin{proposition}
\label{Proposition_Omega_simplex}
 $\Omega_q$ is a simplex, i.e.\  for any $P\in\Omega_q$ there is a unique measure $\pi$ on ${\rm Ex
\,}\Omega_q$ such that
$$
 P=\int_{{\rm Ex\,}\Omega_q} \omega d\pi.
$$
\end{proposition}

\begin{proposition}
\label{Proposition_minimal_martin}
 The minimal boundary is a subset of the Martin boundary of the $q$--Gelfand--Tsetlin graph. More
precisely, if $P\in{\rm Ex \,}\Omega_q$ and $P_k$ is a $q$--coherent system corresponding to $P$,
then $P_k$ belongs to the Martin boundary of the $q$--Gelfand--Tsetlin graph.
\end{proposition}

{\bf Remark.} For the most non-degenerate examples the Martin boundary of the graph coincides with
the minimal boundary. However, there exist graphs for which the minimal boundary is strictly less
than the Martin boundary.

\begin{proposition}
\label{Proposition_boundary_regular}
  Let $P\in{\rm Ex \,}\Omega_q$ and let $\tau\in\mathcal T$. $P$-almost surely the sequence of
signatures $\tau(N)$ is regular and
$$
 P^{\tau(N)}_k\to P_k
$$
 for every $k$.
\end{proposition}

For the proofs see \cite[Theorem 9.2]{Olsh}, \cite[Section 6]{OkOlsh} and \cite[Proposition
10.8]{Olsh} respectively. Similar propositions were proved by Diaconis and Freedman in the
framework of \emph{partial exchangeability}. See \cite[Theorem 1.1]{DF}.

\begin{theorem}
 The set $\Omega_q$ of all $q$--central probability measures on $\mathcal T$ is a simplex with
extreme points $\mathcal E^\nu$, i.e.\  for any $q$--central probability measure $P\in\Omega_q$
there exists a unique probability measure $\pi$ on $\mathcal N$ such that
$$
 P=\int_{\mathcal N} \mathcal E^\nu d\pi.
$$
\end{theorem}
\begin{proof}
Proposition \ref{Proposition_Omega_simplex} implies that $\Omega_q$ is a simplex. Proposition
\ref{Proposition_minimal_martin} and Theorem \ref{Theorem_limits_of_measures} imply that
$$
{\rm Ex \,}\Omega_q\subset\{\mathcal E^\nu\}_{\nu\in\mathcal N},
$$
where $\mathcal E^\nu$ is a $q$--central measure corresponding to $q$--coherent system $\mathcal
E^\nu_k$. It follows from Proposition \ref{proposition_topology_on_extreme_measures} that
$\{\mathcal E^\nu\}_{\nu\in\mathcal N}$ with topology induced from $\Omega_q$ is isomorphic to
$\mathcal N$. Finally, Proposition \ref{Proposition_extremality} implies that  ${\rm Ex
\,}\Omega_q=\{\mathcal E^\nu\}_{\nu\in\mathcal N}$. Indeed, if $Q\in \mathcal E^\nu$, then
$$
  Q=\int_{{\rm Ex\,}\Omega_q} \omega d\pi=\int_{\{\mathcal E^\nu\}_{\nu\in\mathcal N}} \omega
d\pi',
$$
where $\pi'$ is a probability measure on $\{\mathcal E^\nu\}_{\nu\in\mathcal N}$ such that
$$
 \pi'({\rm Ex\,}\Omega_q)=1.
$$
 But then $\pi'$ is $\delta$--measure on a single element $Q=\mathcal E^\nu$. Thus, $\mathcal E^\nu\in{\rm
Ex\,}\Omega_q$.
\end{proof}
The proved theorem is readily seen to be equivalent to Theorem \ref{theorem_Main}. Now Proposition
\ref{Proposition_description_of_regular_point} is a straightforward corollary of Theorem
\ref{theorem_Main} and Proposition \ref{Proposition_boundary_regular}

\section{Some proofs}

\subsection{Relations between Schur and interpolation Schur functions}

\label{Subsection_proofs_spaces_of_functions}

In this section we aim to prove Proposition
\ref{Proposition_all_properties_of_interpolation_decomposition} and Lemma \ref{lemma_G=G'}. To do
that we need some preparations.

Denote by $\mathcal F_N$ the class of $q$--Schur generating functions of probability measures
supported on $\mathbb{GT}_N^+$. In other words, $F(x_1,\dots,x_N)\in \mathcal F_N$ if $F$ is a
symmetric analytic functions $F(x_1,\dots,x_N)$ on $D_N$ such that
\begin{equation} \label{eq_Schur_decomposition} F(x_1,x_2,\dots,x_N)=\sum_{\mu\in\mathbb{GT}^+_N}
\frac{s_{\mu}(x_1,\dots ,x_N)}{s_{\mu}(1,q^{-1},\dots,q^{1-N})} c_\mu,
\end{equation}
for a certain sequence of numbers, $c_\mu$, $\mu\in\mathbb{GT}^+_N$, such that $c_\mu\ge 0$ and
$\sum\limits_{\mu\in\mathbb{GT}^+_N}c_\mu=1$.

 Clearly, \eqref{eq_Schur_decomposition} is essentially just a Taylor series decomposition, thus,
if such decomposition exists, then it is unique.

We also want to consider decompositions of symmetric functions into interpolation Schur
polynomials:
\begin{equation}
 F(x_1,\dots,x_N)=\sum_\mu a_\mu s^*_{\mu}(x_1,x_2\cdot,\dots, x_N;q)
\end{equation}

\begin{proposition}
 \label{proposition_how_to_obtain_interpolation_coefficients} There exist coefficients
 $K_\mu^\lambda$, ($\mu,\lambda\in\mathbb{GT}^+_N$), such that for every function
 $F(x_1,\dots,x_N)$ defined in points $q^{\lambda-\delta}$ for every $\lambda\in\mathbb{GT}^+_N$,
 there is a unique decomposition
\begin{equation}
 \label{eq_simple_interp_polin_decomposition} F(x_1,\dots,x_N)=\sum_\mu a_\mu
 s^*_{\mu}(x_1,x_2\cdot,\dots, x_N;q),
\end{equation}
 where the series converges in points $x=q^{\lambda-\delta}$ for every $\lambda\in\mathbb{GT}^+_N$.
 We have
 \begin{equation}
  \label{eq_defining_interp_coeffs} a_\mu=\sum_{\lambda\in\mathbb{GT}^+_{N}} K_\mu^{\lambda}
  F(q^{\lambda-\delta}).
 \end{equation}
 Furthermore, $K_\mu^\lambda=0$ unless $\lambda\subset\mu$, thus, all sums in
 \eqref{eq_defining_interp_coeffs} are finite.
\end{proposition}
{\bf Remark.} Note that this proposition is valid not only for $0<q<1$ but also for $q>1$.
\begin{proof}[Proof of Proposition \ref{proposition_how_to_obtain_interpolation_coefficients}]
 Substitute $x=q^{\lambda-\delta}$ for every $\lambda\in\mathbb{GT}^+_N$ in
 \eqref{eq_simple_interp_polin_decomposition}. We obtain a system of linear equations
 \begin{equation}
 \label{eq_system_of_interpolation_equations} F(q^{\lambda-\delta})=\sum_{\mu\in\mathbb{GT}^+_N}
 a_\mu L^\mu_\lambda,
 \end{equation}
  where
 $$
  L^\mu_\lambda=s^*_\mu(q^{\lambda-\delta}).
 $$
 Proposition \ref{proposition_interpolation_property} implies that $L^\mu_\lambda=0$, unless
 $\mu\subset\lambda$, and $L^\mu_\mu\neq 0$. Thus, the matrix $L^\mu_\lambda$ has a triangular
 structure (with respect to the partial order on signatures defined in Section
 \ref{section_probabilistic_setup}) and all the sums in
 \eqref{eq_system_of_interpolation_equations} are finite. Consequently, one may solve the system
 \eqref{eq_system_of_interpolation_equations} inductively, starting from $c_\mu$ with $|\mu|=0$,
 then proceeding to $c_\mu$ with $|\mu|=1$ and so on. In other words, there exists a matrix
 $K_\mu^\lambda$ (which also has a triangular structure), such that
 $$
  \sum_{\lambda\in\mathbb{GT}_N^+} K^{\lambda}_{\mu_1} L^{\mu_2}_{\lambda}=\begin{cases} 1,\quad
  \mu_1=\mu_2,\\ 0.\quad otherwise.\end{cases}
 $$
 The coefficients of the matrix $K_\mu^\lambda$ are the desired ones.
\end{proof}

For a general function $F$ we can not claim that the series
\eqref{eq_simple_interp_polin_decomposition} converges in any points other that
$q^{\lambda-\delta}$. However, for functions from $\mathcal F_n$ the following lemma holds:

\begin{lemma}
 \label{lemma_interpolation_schur_decomposition} Suppose that function $F$ belong to $\mathcal
 F_N$. There exists a unique decomposition
 $$
  F(x_1,\dots,x_N)=\sum_\mu a_\mu s^*_{\mu}(x_1,x_2,\dots, x_N ;q),
 $$
 such that the series converges uniformly on any ball $B(0,r)$ with radius $0<r<1$ and in every point
 $q^{\lambda-\delta}$.
\end{lemma}
\begin{proof}
 To simplify the notations let us consider the case $N=1$. The proof in the general case follows the same steps.

 Equality \eqref{eq_schur_through_interpolation} yields
 \begin{equation}
 \label{eq_q_Newton_interpolation} x^l=\sum_{m\le l} \left[\frac{(q^l-1)\dots
 (q^l-q^{m-1})}{(q^m-1)\dots(q^m-q^{m-1})}  \right] (x-1)\dots(x-q^{m-1})
 \end{equation}
 By the definition of $\mathcal F_1$ we have
 $$
  F(x)=\sum_{l\ge 0} c_l x^l,
 $$
 where $c_l\ge 0 $ and $\sum c_l=1$.

 Substituting \eqref{eq_q_Newton_interpolation} we get
 \begin{equation}
 \label{eq_double_sum_schur_interp_1dim} F(x)=\sum_{l\ge 0} c_l \sum_{m\le l}
 \left[\frac{(q^l-1)\dots (q^l-q^{m-1})}{(q^m-1)\dots(q^m-q^{m-1})}  \right] (x-1)\dots(x-q^{m-1})
 \end{equation}
 Observe that the coefficients $\frac{(q^l-1)\dots (q^l-q^{m-1})}{(q^m-1)\dots(q^m-q^{m-1})}$ are
 uniformly bounded in $l,m$. If $x=q^k$, then in \eqref{eq_q_Newton_interpolation} only first $k+1$
 terms are nonzero and if $|x|<r<1$, then the series in \eqref{eq_q_Newton_interpolation} converges
 exponentially fast. In both cases \eqref{eq_double_sum_schur_interp_1dim} absolutely converges and
 we may change the order of summation.

 We obtain
 \begin{equation}
  F(x)=\sum_{m\ge 0} (x-1)\dots(x-q^{m-1}) \cdot \sum_{l\ge m} c_l \left[\frac{(q^l-1)\dots
  (q^l-q^{m-1})}{(q^m-1)\dots(q^m-q^{m-1})}  \right],
 \end{equation}
which is the required decomposition.

The uniqueness of the decomposition follows from Proposition
\ref{proposition_how_to_obtain_interpolation_coefficients}.
\end{proof}

\begin{proposition}
\label{Proposition_from_un_conv_to_interpolaion_coeff_convergence} Suppose that functions $F^n$
and $F$ belong to $\mathcal F_N$ , $F^n \rightrightarrows F$ on $D_N$. Then the coefficients
$a^n_{\mu}$ of the interpolation Schur polynomials expansion
\eqref{eq_simple_interp_polin_decomposition} for the functions $F_n$ converge to the corresponding
coefficients $a_\mu$ of the function $F$.
\end{proposition}
\begin{proof}
Uniform convergence on $D_N$ implies that $F^n(q^{\lambda-\delta})\to F(q^{\lambda-\delta})$.
Applying Proposition \ref{proposition_how_to_obtain_interpolation_coefficients} we conclude that
$a^n_{\mu}\to a_\mu$.
\end{proof}

Combining Lemma \ref{lemma_interpolation_schur_decomposition} with Proposition
\ref{Proposition_from_un_conv_to_interpolaion_coeff_convergence} we arrive at Proposition
\ref{Proposition_all_properties_of_interpolation_decomposition}.

Now let us turn to $q$--interpolation Schur generating functions. Let $\mathcal F^*_N$ denote the
class of $q$--interpolation Schur generating functions of probability measures supported on
$\mathbb{GT}_N^+$. In other words, $F(x_1,\dots,x_N)\in \mathcal F^*_N$ if $F$ is a symmetric
analytic functions $F(x_1,\dots,x_N)$ on $\mathbb C^N$ such that
\begin{equation}
 \label{eq_definition_class_F_star} F(x_1,\dots,x_N)=\sum\limits_{\mu\in\mathbb{GT}^+_N}c_\mu
 \frac{s^*_\mu(q^{N-1}x_1,q^{N-1}x_2,\dots,q^{N-1}x_N;q^{-1})}{s^*_\mu(0,\dots,0;q^{-1})}
\end{equation}
for a certain sequence of numbers, $c_\mu$, $\mu\in\mathbb{GT}^+_N$, such that $c_\mu\ge 0$ and
${\sum_{\mu\in\mathbb{GT}^+_N}c_\mu=1}$.

 Note that if decomposition \eqref{eq_definition_class_F_star} exists, then it is unique. To prove
this fact we repeat the argument of Proposition
\ref{proposition_how_to_obtain_interpolation_coefficients}.

\smallskip

Next, we want to study the relation between $q$--Schur generating function and $q$--interpolation
Schur generating function of the same probability measure.

Let ${\rm Sym(N)}$ denote the space of symmetric polynomials in $N$ variables $x_1,\dots, x_N$.
Consider a linear map $G: {\rm Sym(N)} \to {\rm Sym(N)} $ defined on Schur polynomials' basis
through
\begin{equation}
\label{eq_map_schur_to_interpolation} G\left(
\frac{s_\mu(x_1,\dots,x_N)}{s_\mu(1,q^{-1},\dots,q^{1-N})}\right)=\frac{s^*_\mu(x_1,x_2,\dots,x_N;q)}{s^*_\mu(0,\dots,0;q)},
\end{equation}
or, equivalently,
$$
 G(s_{\mu}(x_1,\dots,x_N))=(-1)^{|\mu|}q^{n(\mu)-n(\mu')} {s^*_{\mu}(x_1,\dots,x_N ; q )}
$$

Note that $s^*_\lambda(x_1,\dots,x_N;1)=s_\lambda(x_1-1,\dots,x_N-1)$. Thus for $q=1$ the map $G$
becomes a simple change of variables:
$$
 G_{q=1}f(x_1,\dots, x_N)=f(1-x_1,\dots,1-x_N).
$$

Let us also consider another map  $G': {\rm Sym(N)} \to {\rm Sym(N)} $ defined on
$q$--interpolation polynomials through:

\begin{equation}
\label{eq_map_interpolation_to_schur}
G'\left(\frac{s^*_\mu(q^{k-1}x_1,q^{k-1}x_2,\dots,q^{k-1}x_k;q^{-1})}{s^*_\mu(0,\dots,0;q^{-1})}\right)=
\frac{s_\mu(x_1,\dots,x_k)}{s_\mu(1,q^{-1},\dots,q^{1-k})}
\end{equation}

\begin{lemma}
\label{lemma_equivalence_G_def}  $G=G'$ on all finite degree polynomials.
\end{lemma}
{\bf Remark.} In one-dimensional case we have:
$$
 G(x^k)=(1-x)(1-xq^{-1})\dots(1-xq^{1-k}),
$$
$$
 G'( (1-x)(1-xq)\dots(1-xq^{k-1}))=x^k
$$
The fact that $G=G'$ follows from the $q$--binomial theorem.
\begin{proof}[Proof of Lemma \ref{lemma_equivalence_G_def}]
Take the equality \eqref{eq_interpolation_through_schur} with $q$ replaced by $q^{-1}$ and with $x$
replaced by $q^{k-1}x$. We get:
\begin{multline*}
 s_{\lambda}^*(q^{k-1}x_1,\dots,q^{k-1}x_N;q^{-1})
\\=\sum_{\mu}\frac{s^*_\mu(q^{\lambda-\delta};q)}{s^*_\mu(q^{\mu-\delta});q)}\frac{s_{\lambda}^*(0,\dots,0;q^{-1}
)}{s_{\mu}^*(0,\dots,0;q^{-1})}{q^{(N-1)|\mu|}s_{\mu}(x_1,\dots,x_N)}.
\end{multline*}

 Apply $G$ in the sense of formula \eqref{eq_map_schur_to_interpolation} to the righthand-side and
 $G'$ in the sense of formula \eqref{eq_map_interpolation_to_schur} to the lefthand-side. We get
 the equality

\begin{multline}
 \label{eq_apply_G}
 s_\lambda(x_1,\dots,x_N)\frac{s^*_\lambda(0,\dots,0;q^{-1})}{s_\lambda(1,q^{-1},\dots,q^{1-N})}
\stackrel{?}{=}\sum_{\mu}\frac{s^*_\mu(q^{\lambda-\delta};q)}{s^*_\mu(q^{\mu-\delta});q)}\frac{s_{\lambda}^*(0,\dots,0;q^{-1}
)}{s_{\mu}^*(0,\dots,0;q^{-1})}\\ \times\frac{s_\mu(1,q^{-1},\dots,q^{1-N})}{s^*_\mu(0,\dots,0;q)}
{q^{(N-1)|\mu|}s^*_\mu(x_1,x_2,\dots,x_Nk;q)}.
\end{multline}

Since $$\frac{q^{(N-1)|\mu|}
s_\mu(1,q^{-1},\dots,q^{1-N})}{s^*_\mu(0,\dots,0;q^{-1}){s_{\mu}^*(0,\dots,0;q^{-1})}}=\frac{1}{s_\mu(1,q^{-1},\dots,q^{1-N})},
$$

 \eqref{eq_apply_G} is equivalent to
$$
 s_\lambda(x_1,\dots,x_N)
\stackrel{?}{=}\sum_{\mu}\frac{s^*_\mu(q^{\lambda-\delta};q)}{s^*_\mu(q^{\mu-\delta};q)}
\frac{s_\lambda(1,q^{-1},\dots,q^{1-N})}{s_\mu(1,q^{-1},\dots,q^{1-N})}
s^*_\mu(x_1,x_2,\dots,x_N;q),
$$

which is exactly \eqref{eq_schur_through_interpolation}. Hence, the last equality is true and $G$
coincides with $G'$ on all $q$--interpolation Schur polynomials. Consequently, they coincide on
all polynomials by linearity.
\end{proof}

Next, we want to extend the domain of definition of the maps $G$ and $G'$.

For any function $f\in\mathcal F^*_N$ we can define $G'(f)$ as follows:
$$
 f(x_1,\dots,x_N)= \sum\limits_{\mu\in\mathbb{GT}^+_N}c_\mu
 \frac{s^*_\mu(q^{N-1}x_1,q^{N-1}x_2,\dots,q^{N-1}x_N;q^{-1})}{s^*_\mu(0,\dots,0;q^{-1})},
$$
$$
 G'(f)= \sum\limits_{\mu\in\mathbb{GT}^+_N}c_\mu
 \frac{s_\mu(x_1,\dots,x_N)}{s_\mu(1,q^{-1},\dots,q^{1-N})}.
$$
It is clear that $G'$ is a bijection between $\mathcal F^*_N$ and $\mathcal F_N$. More precisely,

$$
 G'(\mathcal S^*(x_1,\dots,x_N;P))=\mathcal S(x_1,\dots,x_N;P)
$$
for any probability measure $P$ on $\mathbb{GT}_N$ with ${\rm supp}(P)\subset\mathbb{GT}_N^+$.

\begin{proof}[Proof of Lemma \ref{lemma_G=G'}] Recall that we want to prove the following statement:
Suppose that $F\in\mathcal F^*_N$ has the following Taylor series expansion
 $$
  F(x_1,\dots,x_N)=\sum_{\mu\in\mathbb {GT}^+_N}
  a_{\mu}\frac{s_\mu(x_1,\dots,x_N)}{s_\mu(1,q^{-1},\dots,q^{1-N})}
 $$
 Define
 \begin{equation}
 \label{eq_definition_of_G_extended2} G(f)=\sum_{\mu\in\mathbb {GT}^+_N}
 a_{\mu}\frac{s^*_\mu(x_1,x_2,\dots,x_N;q)}{s^*_\mu(0,\dots,0;q)}.
 \end{equation}
 Then  the series on the right side of \eqref{eq_definition_of_G_extended2} converges uniformly on
 any ball $B(0,r)$ with radius $0<r<1$ and in every point $q^{\lambda-\delta}$. Furthermore,
 $G(f)=G'(f)$.

\smallskip

Let us start the proof. Lemma \ref{lemma_interpolation_schur_decomposition} yields that $G'(f)$ can
be represented as a linear combination of $q$--interpolation Schur polynomials:
$$
   G'(f)=\sum_{\mu\in\mathbb {GT}^+_N}
   b_{\mu}\frac{s^*_\mu(x_1,x_2,\dots,x_N;q)}{s^*_\mu(0,\dots,0;q)}.
$$
 This series converges uniformly on any ball $B(0,r)$ with radius $0<r<1$ and in every point
 $q^{\lambda-\delta}$. To prove the proposition we should check that $a_\mu=b_\mu$.

 Suppose that
 $$
  F(x_1,\dots,x_N)=\sum\limits_{\mu\in\mathbb{GT}^+_N}c_\mu
  \frac{s^*_\mu(q^{N-1}x_1,q^{N-1}x_2,\dots,q^{N-1}x_N;q^{-1})}{s^*_\mu(0,\dots,0;q^{-1})}.
 $$
 Set
 $$
   f^m(x_1,\dots,x_N)=\sum\limits_{|\mu|\le m}c_\mu
   \frac{s^*_\mu(q^{N-1}x_1,q^{N-1}x_2,\dots,q^{N-1}x_N;q^{-1})}{s^*_\mu(0,\dots,0;q^{-1})}
 $$
 and $\widehat f^m=F-f^m$. Let $a^m_\mu$ be the corresponding coefficient of the Taylor series
 expansion of $f^m$ and let $\widehat a^m_\mu$ be the corresponding coefficient of the Taylor
 series expansion of $\widehat f^m$. Clearly, $a_\mu=a^m_\mu+\widehat a^m_\mu$. Represent $G'(f^m)$
 and $G'(\widehat f^m)$ as linear combinations of $q$--interpolation Schur polynomials and let
 $b_\mu^m$, $\widehat b_\mu^m$ be the corresponding coefficients:
 $$
   G'(f^m)=\sum_{\mu\in\mathbb {GT}^+_k}
   b^m_{\mu}\frac{s^*_\mu(x_1,x_2,\dots,x_k;q)}{s^*_\mu(0,\dots,0;q)},
 $$
 $$
   G'(\widehat f^m)=\sum_{\mu\in\mathbb {GT}^+_k} \widehat
   b^m_{\mu}\frac{s^*_\mu(x_1,x_2,\dots,x_k;q)}{s^*_\mu(0,\dots,0;q)},
 $$
 It is clear that $b_\mu=b_\mu^m+\widehat b_\mu^m$.

 Lemma \ref{lemma_equivalence_G_def} implies that $b_\mu^m=a_\mu^m$. It remains to prove that both
 $\widehat a_\mu^m$ and $\widehat b_\mu^m$ tend to zero as $m$ tends to infinity.

 It follows from Proposition \ref{proposition_convergence_in_def_of_class_F_star} that $\widehat
 f^m \rightrightarrows 0$ uniformly on compact sets as $m\to\infty$. Uniform convergence of
 analytical functions implies convergence of their Taylor expansion coefficients. Thus, $\widehat
 a_\mu^m \to 0$ as $m\to\infty$. By the definition
$$
 G'(f^m)= \sum\limits_{|\mu|>m}c_\mu \frac{s_\mu(x_1,\dots,x_N)}{s_\mu(1,q^{-1},\dots,q^{1-N})}.
$$

Proposition \ref{proposition_Schur_gen_function_converges} implies that $G'(f^m) \rightrightarrows
0$ uniformly on $D_N$. Repeating the argument in the proof of Proposition
\ref{Proposition_from_un_conv_to_interpolaion_coeff_convergence} we conclude that $\widehat
b_\mu^m \to 0$ as $m\to\infty$.
\end{proof}

\subsection{$q$--Central measures and probability generating functions}

\label{Subsection_proofs_probability_gen_functions}

\begin{proof}[Proof of Proposition \ref{Proposition_coherence_from_centrality}]
We want to prove the following statement: if measure $P$ on $\mathcal T_N$  is such that
 $$ P(\tau(1)\prec\dots\tau(N))=\frac{q^{|\tau(1)|+\dots+|\tau(N-1)|}}{{\rm
 Dim}_q(\tau(N))} P_N(\tau(N)), $$
 for any path $\tau\in\mathcal T_N$, then $P_1,P_2,\dots,P_N$ is
 a $q$--coherent system.

We proceed by induction in $N$. Case $N=1$ is trivial. For general $N$ we have

 \begin{multline*}
  P_{N-1}(\mu)=\sum_{\tau\in\mathcal T_N,\, \tau(N-1)=\mu}
  P_N(\tau(N))\frac{q^{|\tau(1)|+\dots+|\tau(N-1)|}}{{\rm Dim}_q(\tau(N))}\\=\sum_\lambda
  \frac{P_N(\lambda)q^{|\mu|}}{{\rm Dim}_q (\lambda)}\sum_{\tau\in\mathcal T_N,\,
  \tau(N)=\lambda,\tau(N-1)=\mu} q^{|\tau(1)|+\dots+|\tau(N-2)|}\\=P_N(\lambda)q^{|\mu|}\frac{{\rm
  Dim}_q(\mu)}{{\rm Dim}_q(\lambda)}=\sum_\lambda P_N(\lambda)P(\lambda\to\mu)
\end{multline*}
 Thus, $P_{N-1}$ and $P_N$ are $q$--coherent.

Next, projection of measure $P$ on $\mathcal T_{N-1}$ defines a probability measure that we denote
by $\tilde P$. If $\tau'\in\mathcal T_{N-1}$, then
$$
\tilde P(\tau')= \sum_\lambda P(\tau'(1)\prec\dots\prec\tau'(N-1)\prec\lambda)
=q^{|\tau(1)|+\dots+|\tau(N-1)|} \sum_{\lambda\mid \tau'(N-1)\prec\lambda}\frac{P_N(\lambda)}{{\rm
Dim}_q(\lambda)},
$$
It follows that if $\tau^1,\tau^2\in\mathcal T_{N-1}$ and $\tau^1(N-1)=\tau^2(N-1)$, then
$$
\frac{\tilde P(\tau^1)}{\tilde
P(\tau^2)}=\frac{q^{|\tau^1(1)|+\dots+|\tau^1(N-2)|}}{q^{|\tau^2(1)|+\dots+|\tau^2(N-2)|}},
$$
Thus,
\begin{multline*}
 \tilde P(\tau')=\frac{q^{|\tau'(1)|+\dots+|\tau'(N-2)|}}{\sum\limits_{\tau\in\mathcal T_{N-1},\,
 \tau(N-1)=\tau'(N-1)} q^{|\tau(1)|+\dots+|\tau(N-2)|}} \sum_{\tau\in\mathcal T_{N-1},\,
 \tau(N-1)=\tau'(N-1)} \tilde P(\tau)\\=\frac{q^{|\tau'(1)|+\dots+|\tau'(N-2)|}}{{\rm
 Dim}_q(\tau'(N-1))} \tilde P_{N-1}(\tau'(N-1))
\end{multline*}
Consequently, it follows by induction that $\tilde P_1,\dots,\tilde P_{N-1}$ form a $q$--coherent
system. Since $P_i=\tilde P_i$ for $i=1,\dots,N-1$ and measures $P_N$ and $P_{N-1}$ are
$q$--coherent, $P_1,\dots P_N$ is a $q$--coherent system.
\end{proof}

\begin{proof}[Proof of Proposition \ref{Proposition_from_coherent_system_to_measure}]
 We want to prove that for any $q$--coherent system $P_1,P_2,\dots$, there exists a unique
$q$--central measure $P$ such that $P_k$ is a projection of $P$ on $\mathbb{GT}_k$ for every
$k=1,2,\dots$.

 For any $N$ let $P^{(N)}$ be a probability measure on $\mathcal T_N$ defined through
 $$
  P^{(N)}(\tau)=\frac{q^{|\tau(1)|+\dots+|\tau(N-1)|}}{{\rm Dim}_q(\tau(N))} P_N(\tau(N)).
 $$
 Proposition \ref{Proposition_coherence_from_centrality} yields that projections $P^{(N)}_k$ of
 $P^{(N)}$ on $\mathbb{GT}_k$ form a $q$--coherent system. Since $P^{(N)}_N=P_N$ and $P_1,\dots,
 P_N$ also form a $q$--coherent system, thus, $P^{(N)}_k=P_k$.

Repeating the argument of Proposition \ref{Proposition_coherence_from_centrality} we conclude that
projection of $P^{(N)}$ on $\mathcal T_{N-1}$ coincides with $P^{(N-1)}$.

Let $P$ be projective limit of the measures $P^{(N)}$ as $N\to\infty$. Clearly, $P$ is a
$q$--central measure on $\mathcal T$ and projections of $P$ on $\mathbb{GT}_k$ are exactly $P_k$.
\end{proof}

\begin{proof}[Proof of Proposition \ref{Proposition_coherency_in_terms_of_Schur_polynomials}]
  We want to prove that two probability measures $P_N$ and $P_{N+1}$ on $\mathbb{GT}_N$ and
$\mathbb{GT}_{N+1}$, respectively, are $q$--coherent if and only if
$$
 {\mathcal S}(x_1,\dots,x_N;P_{N})={\mathcal S}(x_1,\dots,x_N,q^{-N};P_{N+1}).
$$

 Using Proposition \eqref{proposition_Schur_branching_rule} we obtain
 \begin{multline*}
  \frac{s_\lambda(x_1,\dots, x_N,q^{-N})}{s_\lambda(1,q^{-1},\dots,q^{-N})}=
  \dfrac{\sum\limits_{\mu\prec\lambda} s_\mu(x_1,\dots, x_N)
  q^{-N(|\lambda|-|\mu|)}}{q^{-N|\lambda|}s_\lambda(1,q,\dots,q^{N})}\\= \sum_{\mu\prec\lambda}
  \frac{s_\mu(x_1,\dots, x_N)}{s_\mu(1,q^{-1},\dots, q^{1-N})} q^{|\mu|}\frac{s_\mu(1,q,\dots,
  q^{N-1})}{s_\lambda(1,q,\dots,q^{N})}\\= \sum_{\mu\in\mathbb{GT}_N} \frac{s_\mu(x_1,\dots,
  x_N)}{s_\mu(1,q^{-1},\dots, q^{1-N})} P(\lambda\to\mu).
 \end{multline*}
 Thus,
 \begin{multline*}
{\mathcal S}(x_1,\dots,x_N,q^{-N};P_{N+1}) =\sum_{\lambda\in\mathbb{GT}_{N+1}}
\frac{s_\lambda(x_1,\dots,x_N,q^{-N})}{s_\lambda(1,q^{-1},\dots,q^{-N})}P_{N+1}(\lambda)\\=
\sum_{\mu\in\mathbb{GT}_N} \frac{s_\mu(x_1,\dots, x_N)}{s_\mu(1,q^{-1},\dots, q^{1-N})}
\sum_{\lambda\in\mathbb{GT}_{N+1}} P(\lambda\to\mu) P_{N+1}(\lambda),
  \end{multline*}
 On the other hand,
$$
 {\mathcal S}(x_1,\dots,x_N;P_{N}) = \sum_{\mu\in\mathbb{GT}_N} \frac{s_\mu(x_1,\dots,
 x_N)}{s_\mu(1,q^{-1},\dots, q^{1-N})} P_{N}(\mu).
$$
Comparing the last two expressions we conclude that
$${\mathcal
S}(x_1,\dots,x_N,q^{-N};P_{N+1})={\mathcal S}(x_1,\dots,x_N;P_{N})$$ if and only if
 $$P_{N}(\mu)=\sum_{\lambda\in\mathbb{GT}_{N+1}} P(\lambda\to\mu)
P_{N+1}(\lambda)$$ for all $\mu\in\mathbb{GT}_N$.
\end{proof}


\begin{proof}[Proof of Proposition \ref{Proposition_coherency_in_terms_of_interpolation_polynomials}]

 We want to prove that two probability measures $P_N$ and $P_{N+1}$ on $\mathbb{GT}^+_N$ and
$\mathbb{GT}_{N+1}$, respectively,  are $q$--coherent if and only if
$$
 {\mathcal S^*}(x_1,\dots,x_N;P_{N})={\mathcal S^*}(x_1,\dots,x_N,0;P_{N+1}).
$$

To prove this proposition we need a lemma.
\begin{lemma}[The branching rule for $q$--interpolation Schur polynomials]
For any $\lambda\in\mathbb{GT}^+_{N+1}$ we have
\begin{multline*}
\frac{s^*_\lambda(q^{N}x_1,q^{N}x_2,\dots,q^{N}x_N,0;q^{-1})}{s^*_\lambda(0,\dots,0;q^{-1})}\\=
\sum_{\mu\in\mathbb{GT}_N}
\frac{s^*_\mu(q^{N-1}x_1,q^{N-1}x_2,\dots,q^{N-1}x_N;q^{-1})}{s^*_\mu(0,\dots,0;q^{-1})}
P(\lambda\to\mu).
\end{multline*}
\end{lemma}
\begin{proof}
 Using the combinatorial formula (see Proposition \ref{Proposition_combinatoral_for_interp_Schur})
 we obtain
$$
 s^*_{\lambda}(q^Nx_1,\dots,q^Nx_N,q^Nx_{N+1};q^{-1})=\sum_{T} \prod\limits_{(i,j)\in\lambda}
 \left(q^N x_{T(i,j)}-q^{i-j-T(i,j)+N+1}\right),
$$
where the sum is taken over all semistandard Young tableaux $T(i,j)$ of shape $\lambda$ filled with
numbers $1, \dots, N+1$. Note that the part of $T$ filled with numbers $1,\dots, N$ is a Young
diagram of shape $\mu$ such that $\mu\prec\lambda$. Consequently, substituting $x_{N+1}=0$ we get
\begin{multline*}
 s^*_{\lambda}(q^Nx_1,\dots,q^Nx_N,0;q^{-1})\\=\sum_{\mu\prec\lambda}
 \left(\prod\limits_{(i,j)\in\lambda\setminus\mu} \left(0-q^{i-j}\right)\right) \sum_{T}
 \prod\limits_{(i,j)\in\mu} \left(q^N x_{T(i,j)}-q^{i-j-T(i,j)+N+1}\right),
\end{multline*}
where the sums are taken over semistandard Young tableaux $T(i,j)$ of shape $\mu$ filled with
numbers $1, \dots, N$. Consequently,
\begin{multline*}
 s^*_{\mu}(q^Nx_1,\dots,q^Nx_N,0;q^{-1}) \\=\sum_{\mu\prec\lambda}\left(
 \prod\limits_{(i,j)\in\lambda\setminus\mu} \left(-q^{i-j}\right)\right) q^{|\mu|}
 s^*_{\mu}(q^{N-1}x_1,\dots,q^{N-1}x_N;q^{-1}).
\end{multline*}
It remains to prove that
\begin{equation}
\label{eq_equation_in_int_schur_coherency_lemma} \prod\limits_{(i,j)\in\lambda\setminus\mu}
\left(-q^{i-j}\right) q^{|\mu|}=\frac{s^*_\lambda(0,\dots,0;q^{-1})}{s^*_\mu(0,\dots,0;q^{-1})}
P(\lambda\to\mu).
\end{equation}
We claim that for any $\lambda\in\mathbb{GT}_N$
\begin{multline*}
 s^*_\lambda(0,\dots,0;q^{-1})=s_\lambda(1,q^{-1},\dots,q^{1-N}) q^{(N-1)|\lambda|}
\prod_{(i,j)\in\lambda}
  \left(-q^{i-j}\right) \\=s_\lambda(1,q,\dots,q^{N-1}) \prod_{(i,j)\in\lambda}
 \left(-q^{i-j}\right)
\end{multline*}
The claim follows from the combinatorial formulas for Schur polynomials and for $q$--interpolation
Schur polynomials. Using the last relation one immediately proves
\eqref{eq_equation_in_int_schur_coherency_lemma}.
\end{proof}

 To finish the proof of Proposition \ref{Proposition_coherency_in_terms_of_interpolation_polynomials}
 we use the last lemma and repeat the argument of Proposition
\ref{Proposition_coherency_in_terms_of_Schur_polynomials}.
\end{proof}

\begin{proof}[Proof of Proposition \ref{Proposition_convergene_of_measures_implies_uniform_convergence_of_function}]
We want to prove that if a sequence of probability measure $P^i_N$ on $\mathbb{GT}_N$.  weakly
converge to $P_N$, then
$$
 {\mathcal S}(x_1,\dots,x_N;P_N^i)\rightrightarrows {\mathcal S}(x_1,\dots,x_N;P_N)
$$
uniformly on $T_N$. Furthermore, if  $P_N$, $P^i_N$ are supported on $\mathbb{GT}_N^+$, then the
convergence is uniform on $D_N$.

Since $P_N$ is a probability measure on $\mathbb{GT}_N$, there exists $k$ such that $P_N(\{\mu:\,
-k<\mu_i<k\})>1-\varepsilon$. Observe that $\{\mu:\, -k<\mu_i<k\}$ is a finite set with less than
$(2k)^N$ elements. Thus, weak convergence of measures $P_N^i$ implies that for all $i>i_0$ and all
$\lambda\in\mathbb{GT}_N$ such that $-k<\lambda_i<k$ we have
$$
 \left| P_N^i(\lambda)-P_N(\lambda)\right| <\frac{\varepsilon}{(2k)^N}
$$
 and
$$
 P_N^i(\{\mu:\, |\mu|<k\})>1-2\varepsilon.
$$

Consequently, since
$$
 \left| \frac{s_\mu(x_1,\dots,x_N)}{s_\mu(1,\dots,q^{1-N})}\right| \le 1
$$
on $T_N$, we have for $i>i_0$ and $(x_1,\dots,x_N)\in T_N$
\begin{multline*}
 \left| {\mathcal S}(x_1,\dots,x_N;P_N^i)- {\mathcal S}(x_1,\dots,x_N;P_N) \right|\\ \le
 \sum_{|\mu|<k} |P_N^i(\mu)- P_N(\mu)| \left|\frac{s_\mu(x_1,\dots,x_N)}{s_\mu(1,\dots,q^{1-N})}
 \right| \\+\sum_{|\mu|\ge k} (P_N^i(\mu)+ P_N(\mu))
 \left|\frac{s_\mu(x_1,\dots,x_N)}{s_\mu(1,\dots,q^{1-N})}\right| \le
 k^N\frac\varepsilon{k^N}+\varepsilon+2\varepsilon=4\varepsilon
\end{multline*}
It follows, that ${\mathcal S}(x_1,\dots,x_N;P_N^i)$ converges uniformly on $T_N$.

If  $P_N$, $P^i_N$ are supported on $\mathbb{GT}_N^+$, then the argument is the same with $T_N$
replaced by $D_N$.
 \end{proof}

\begin{proof}[Proof of Proposition
\ref{proposition_Convergence_of_S_implies_weak_convergence}] We want to prove the following
statement: if $P^i_N$ is a sequence of probability measures on $\mathbb{GT}_N$ such that the
functions
$$
 {\mathcal S}(x_1,\dots,x_N;P_N^i)
$$
 converge uniformly on $T_N$ to a function $S(x_1,\dots,x_N)$, then there exists a probability
 measure $P_N$ such that
$$
 S(x_1,\dots,x_N)={\mathcal S}(x_1,\dots,x_N;P_N)
$$
and $P_N^i$ weakly converge to $P_N$.

\smallskip

Recall that (by its definition) a rational Schur function is the ratio of the alternating sum of
monomials and the Vandermonde determinant:
$$
 s_\lambda(x_1,\dots,x_N)=\frac{{\rm Alt}(x^{\lambda_1+N-1}\cdots
x^{\lambda_N})}{V(x_1,\dots,x_N)},
$$
where
$$
 V(x_1,\dots,x_N)=\prod_{1\le i<j \le N} (x_i-x_j)
$$
and
$$
 {\rm Alt}(x_1^{k(1)}\dots x_N^{k(N)})=\sum_{\sigma} (-1)^{\sigma} x_1^{k(\sigma(1))}\cdots
x_N^{k(\sigma(N))},
$$
sum is taken over all permutations of length $N$ and $(-1)^{\sigma}$ is the sign of the permutation
$\sigma$.

The following estimates are useful in the sequel:
\begin{lemma}
\label{lemma_Alternating_sum_estimate} There exist three constants $C_1(N)$, $C_2(N)$, $C_3(N)$
such that for any $(x_1,\dots,x_N)\in T_N$ and any $\lambda\in\mathbb{GT}_N$ we have
$$
 C_1(N)<|V(x_1,\dots,x_N)|<C_2(N)
$$
and absolute values of the coefficients of the monomials in the alternating sum
$$
 \frac{{\rm Alt}(x^{\lambda_1+N-1}\cdots x^{\lambda_N})}{s_\lambda(1,q^{-1},\dots,q^{1-N})}
$$
are bounded by $C_3(N)$.
\end{lemma}
We leave the proof of this lemma to the reader.

Denote
$$
R^i(x_1,\dots,x_N)=V(x_1,\dots,x_N) {\mathcal S}(x_1,\dots,x_N;P_N^i)
$$
and
$$
R(x_1,\dots,x_N)=V(x_1,\dots,x_N) S(x_1,\dots,x_N).
$$
Clearly, $R^i(x_1,\dots,x_N)$ converges to $R(x_1,\dots,x_N)$ uniformly on $T_N$.

We have
\begin{equation}
\label{eq_x12}
 R^i(x_1,\dots,x_N)=\sum_{\lambda\in\mathbb{GT}_N}P_N^i(\lambda)\frac{{\rm
Alt}(x^{\lambda_1+N-1}\cdots x^{\lambda_N})}{s_\lambda(1,q^{-1},\dots,q^{1-N})}.
\end{equation}
Note that if we expand \eqref{eq_x12} in a single sum of monomials, then we get a Fourier series
expansion of $R^i(x_1,\dots,x_N)$ in the conventional sense. Furthermore, estimates of Lemma
\ref{lemma_Alternating_sum_estimate} guarantee that this Fourier series uniformly converges on
$T_N$. It is well known that uniform convergence of the continuous functions implies the
convergence of their Fourier coefficients. (This fact follows from the integral formula for the
Fourier coefficients.) Consequently, $P_N^i(\lambda)$ converges to a certain number $P_N(\lambda)$
for every $\lambda$. Since $P_N^i(\lambda)\ge 0$ and $\sum_\lambda P_N^i(\lambda)=1$, we
 have
\begin{equation}
\label{eq_x13}
 P_N(\lambda)\ge 0, \quad \sum_{\lambda\in{\mathbb GT}_N} P_N(\lambda) \le 1.
\end{equation}
If the last sum equals to $1$, then the numbers $P_N(\lambda)$ define a probability measure $P_N$.
Measures $P_N^i$ weakly converge to $P_N$ and we are done. However, the proof of the fact that $
\sum_\lambda P_N(\lambda) =1$ needs an additional argument. Define
$$
 \widehat R(x_1,\dots x_N)=\sum_{\lambda\in\mathbb{GT}_N}P_N(\lambda)\frac{{\rm
Alt}(x^{\lambda_1+N-1}\cdots x^{\lambda_N})}{s_\lambda(1,q^{-1},\dots,q^{1-N})}.
$$
Inequalities \eqref{eq_x13} and estimates of Lemma \ref{lemma_Alternating_sum_estimate} guarantee
that $\widehat R(x_1,\dots x_N)$ is well-defined on $T_N$. Now observe that by the construction all
the Fourier coefficients of $R$ and $\widehat R$ coincide. Thus, all Fourier coefficients of
$R-\widehat R$ are equal to zero. Consequently, (since Laurent polynomials are dense in the space
of continuous functions of $T_N$) $R=\widehat R$. Therefore
\begin{multline*}
 \sum_{\lambda\in{\mathbb GT}_N} P_N(\lambda)=\frac{\widehat
R(1,q^{-1},\dots,q^{-N})}{V(1,q^{-1},\dots,q^{-N})}=\frac{R(1,q^{-1},\dots,q^{-N})}{V(1,q^{-1},\dots,q^{-N})}\\=
\lim_{i\to\infty}\frac{R^i(1,q^{-1},\dots,q^{-N})}{V(1,q^{-1},\dots,q^{-N})}=\lim_{i\to\infty}
1=1.
\end{multline*}

\end{proof}

\begin{proof}[Proof of Proposition \ref{Proposition_convergence_of_measures_implies_uniform_convergence_of_interpolation_generating_functions}]

We want to prove that a sequence of probability measures $P_N^i$ on $\mathbb{GT}^+_N$ weakly
converges to $P_N$ if and only if
$$
 {\mathcal S^*}(x_1,\dots,x_N;P_N^i)\rightrightarrows {\mathcal S^*}(x_1,\dots,x_N;P_N)
$$
uniformly on compact subsets of $\mathbb C^N$.

\smallskip

 If functions ${\mathcal S^*}(x_1,\dots,x_N;P_N^i)$ converge, then using the argument of
 Proposition \ref{proposition_how_to_obtain_interpolation_coefficients} we conclude that the
 coefficients of their $q$--interpolation Schur expansions converge. Thus, measures $P^i_N$ weakly
 converge.

The second part of the proof repeats the argument of Proposition
\ref{Proposition_convergene_of_measures_implies_uniform_convergence_of_function} with inequality
$$
 \left| \frac{s_\mu(x_1,\dots,x_N)}{s_\mu(1,\dots,q^{1-N})}\right| \le 1
$$
replaced by the bound
$$
\left|\frac{s^*_\mu(q^{N-1}x_1,q^{N-1}x_2,\dots,q^{N-1}x_N;q^{-1})}{s^*_\mu(0,\dots,0;q^{-1})}\right|<A(M)
$$
valid for $|x_i|\le M$, as was explained in the proof of Proposition
\ref{proposition_convergence_in_def_of_class_F_star}.
\end{proof}

\subsection{Tightness (Proof of Proposition \ref{proposition_tightness})}


\label{Subsection_proofs_tightness}

We want to prove that if $\lambda(i)$ is a sequence of signatures stabilizing to $\nu$, then the
family of functions
$$
g_i(x_1,\dots,x_k)=\frac{s_{\lambda(i)}(x_1,x_2,\dots,x_k,q^{-k},q^{-k-1},\dots,q^{1-N(i)})}{s_\lambda(i)(1,q,\dots,q^{1-N(i)})}
$$
is a relatively compact subset of the set of continuous functions on $k$--dimensional torus $T_k$
with uniform convergence topology.

\smallskip
We need the following lemma.
\begin{lemma}
\label{lemma_derivative_estimate} Let $\lambda(i)$ be a sequence of positive signatures stabilizing
to $\nu$ with $\nu_1=0$. For any integer $k\ge 0$ set
$$
V_{k,i}(\phi)=
\frac{s_{\lambda(i)}(q^{-k}e^{-i\phi},q^{-1},\dots,q^{1-N(i)})}{s_\lambda(i)(1,q^{-1},\dots,q^{1-N(i)})}.
$$
There exists a function $c(k,\nu)$ and a number $i_0$ such that
$$
 \left| \frac{\partial V_{k,i}(\phi)}{\partial \phi}\right|\le c(k,\nu)
$$
for $i>i_0$.
\end{lemma}
\begin{proof}
 Assume without lost of generality that $\lambda(i)_{N(i)}=0$. Denote $x=q^{-k}e^{-i\phi}$. The
 branching rule for Schur polynomials yields
\begin{multline*}
V_{k,i}(\phi)=\sum_{\mu\prec\lambda(i)} x^{|\lambda|-|\mu|}
\frac{s_{\mu}(q^{-1},\dots,q^{1-N(i)})}{s_{\lambda(i)}(1,q^{-1},\dots,q^{1-N(i)})}\\=
\frac{s_{\lambda(i)}(q^{-1},\dots,q^{1-N(i)})}{s_{\lambda(i)}(1,q^{-1},\dots,q^{1-N(i)})}
\sum_{m\ge 0}x^{m} \sum_{\mu\prec\lambda(i),\,|\lambda(i)|-|\mu|=m}
\frac{s_{\mu}(q^{-1},\dots,q^{1-N(i)})}{s_{\lambda(i)}(q^{-1},\dots,q^{1-N(i)})} .
\end{multline*}
(Here we used the notation $\lambda(i)$ both for $(\lambda(i)_1,\dots,\lambda(i)_{N(i)})$ and
$(\lambda(i)_1,\dots,\lambda(i)_{N(i)-1})$.) Thus,
\begin{multline}
\label{eq_x7} \left| \frac{\partial V_{k,i}(\phi)}{\partial \phi}\right|\le
\frac{s_{\lambda(i)}(q^{-1},\dots,q^{1-N(i)})}{s_{\lambda(i)}(1,q^{-1},\dots,q^{1-N(i)})}
\\ \times \sum_{m\ge 0}m|x|^m \sum_{\mu\prec\lambda(i),\,|\lambda(i)|-|\mu|=m}
\frac{s_{\mu}(q^{-1},\dots,q^{1-N(i)})}{s_{\lambda(i)}(q^{-1},\dots,q^{1-N(i)})} .
\end{multline}
The branching rule for Schur polynomials applied to $s_{\lambda(i)}(1,q^{-1},\dots,q^{1-N(i)})$
yields that the first fraction in the last formula is not greater than $1$. To estimate the double
sum, we first consider the terms with $m=1$. Using \cite[Example 3.1]{Mac} we obtain
\begin{multline*}
\frac{s_{\mu}(q^{-1},\dots,q^{1-N(i)})}{s_{\lambda(i)}(q^{-1},\dots,q^{1-N(i)})}=q^{|\mu|(1-N(i))-|\lambda(i)|(1-N(i))}
\frac{s_{\mu}(1,q,\dots,q^{N(i)-2})}{s_{\lambda(i)}(1,q,\dots,q^{N(i)-2})}\\=q^{1-N(i)}q^{n(\mu)-n(\lambda(i))}\prod_{p<q}\frac{1-q^{\mu_p-\mu_q-p+q}}{1-q^{\lambda(i)_p-\lambda(i)_q-p+q}}
\end{multline*}
Suppose that $\mu$ differs from $\lambda(i)$ in row $N(i)-j$, then
$$
q^{1-N(i)} q^{n(\mu)-n(\lambda(i))}=q^j
$$
and
$$
 \prod_{p<q}\frac{1-q^{\mu_p-\mu_q-p+q}}{1-q^{\lambda(i)_p-\lambda(i)_q-p+q}}\le
 \prod_{\ell=1}^{\infty}(1-q^{\ell})^{-1}.
$$
It follows that the terms corresponding to $m=1$ in the double sum in \eqref{eq_x7} are bounded
from above by
$$
 \sum_{j=1}^{\infty} q^j \prod_{\ell=1}^{\infty}(1-q^{\ell})^{-1}.
$$
For general $m$ let $\mu_n=\lambda(i)_n-f_{N(i)-n}$, $n=1,2,\dots,N(i)-1$. Similarly to $m=1$ case
we get
$$
 \frac{s_{\mu}(q^{-1},\dots,q^{1-N(i)})}{s_{\lambda(i)}(q^{-1},\dots,q^{1-N(i)})} \le
 q^{f_1+2f_2+\dots+(N(i)-1)f_{N(i)-1}} \prod_{\ell=1}^{\infty}(1-q^{\ell})^{-(f_1+f_2+\dots+f_n)}.
$$
Consider the generating function
$$
 a_i(t)=\sum_{f_1,f_2,\dots} q^{f_1+2f_2+3f_3\dots}\cdot t^{f_1+f_2+f_3+\dots},
$$
where the sum is taken over all finite collections of integers $\{f_n\}$ satisfying $0\le f_n\le
\lambda_{N(i)-n}-\lambda_{N(i)-n+1}$. We will use $a_i'(t)$ for $\frac{\partial}{\partial t}
a_i(t)$. The above arguments prove that
$$
\left| \frac{\partial V_{k,i}(\phi)}{\partial \phi}\right|\le
a_i'\left(q^{-k}\prod_{\ell=1}^{\infty}(1-q^{\ell})^{-1}\right).
$$
Now choose $r$ such that $q^{-k}\prod_{\ell=1}^{\infty}(1-q^{\ell})^{-1}<q^{-r}$. Suppose that $i$
is large enough, so that $\lambda(i)_{N(i)-n+1}=\nu_n$ for $n=1,2,\dots,r+1$. Let
$$
b(t)=\prod_{j=1}^{r}\left(\sum_{u=0}^{\nu_{j+1}-\nu_j}\left(tq^{j}\right)^u\right) \cdot
\prod_{j=r+1}^{\infty}\left(\sum_{u=0}^{\infty}\left(tq^{j}\right)^u\right).
$$
Clearly, $a_i(t)\le b(t)$, furthermore, $a_i'(t)\le b'(t)$. But
$$
 b(t)=\prod_{j=1}^{r}\left(\sum_{u=0}^{\nu_{j+1}-\nu_j}\left(tq^{j}\right)^u\right) \cdot
 \prod_{j=r+1}^{\infty}\left(1-\left(tq^{j}\right)\right)^{-1}.
$$
Consequently, $b(t)$ is an analytic function in ball $\{t: t\le q^{-r}\}$. Thus,
$$a_i'\left(q^{-k}\prod_{\ell=1}^{\infty}(1-q^{\ell})^{-1}\right)<b'\left(q^{-k}\prod_{\ell=1}^{\infty}(1-q^{\ell})^{-1}\right)<\infty.$$
\end{proof}

\begin{proof}[Proof of Proposition \ref{proposition_tightness}]
First, let $\nu_1=0$. Observe that $g_i(x_1,\dots,x_k)$ is a symmetric polynomial with positive
coefficients (for large enough $i$). Thus, for any $(x_1,\dots,x_k)\in T_k$  and any $1\le \ell
\le k$ we have
\begin{multline*}
\left| \frac{\partial}{\partial \phi}
g_i(\dots,x_{\ell-1},e^{i\phi}x_\ell,x_{\ell+1},\dots)\right|\\ \le \left|
\frac{\partial}{\partial \phi}
g_i(1,q^{-1},\dots,q^{2-\ell},q^{1-\ell}e^{i\phi},q^{-\ell},\dots,q^{1-k})\right|\\= \left|
\frac{\partial}{\partial \phi}
g_i(q^{1-\ell}e^{i\phi},1,q^{-1},\dots,q^{2-\ell},q^{-\ell},\dots,q^{1-k})\right|\\ \le
\left|\frac{\partial}{\partial \phi}
g_i(q^{1-\ell}e^{i\phi},q^{-1},q^{-2},\dots,q^{1-k})\right|<const,
\end{multline*}
where the last inequality is Lemma \ref{lemma_derivative_estimate}. The above uniform estimate for
the derivatives yields that the family of functions $g_i$ is equicontinuous on $T_n$. We also have
$|g_i|\le 1$ on $T_n$. Thus, by the Arzel\`{a}--Ascoli theorem the set of functions $\{g_i\}$ is
relatively compact.

\smallskip

 For general $\nu$ we note that if $\lambda(i)$ stabilizes to $\nu$ then $A_{-\nu_1}(\lambda(i))$
 stabilizes to $\nu'=A_{-\nu_1}(\nu)$ with $\nu'_1=0$. Since for any $\lambda$ and any $\ell$
$$
 s_{A_{\ell}(\lambda)}(x_1,\dots,x_N)=(x_1\cdots x_N)^{\ell}s_{\lambda}(x_1,\dots,x_N),
$$
the case of general $\nu$ reduces to the case $\nu_1=0$.
\end{proof}

\subsection{Analyticity (proof of Proposition 5.10)}

\label{Subsection_proofs_analyticity}

We want to prove that for every $\nu$ the series $$
\sum_{\mu\in\mathbb{GT}_k^+}(-1)^{|\mu|}q^{n(\mu)-n(\mu')}{\rm
Spec}_{\nu}(s_\mu){s^*_{\mu}(x_1,\dots,x_k)}$$ converges for all $x_1,\dots x_k$ and defines an
entire function.

\smallskip

The combinatorial formula for $s^*_\mu$ (Proposition
\ref{Proposition_combinatoral_for_interp_Schur}) implies that for every $M$ there exists a constant
$C(M)$ such that for every $x_1,\dots,x_k$ with $|x_i|<M$, we have
$$
 \left|s^*_{\mu}(x_1,\dots,x_k)\right|<C(M)^{|\mu|}<\left(C(M)^k\right)^{\mu_1}.
$$

Let us fix $\nu$ and estimate $q^{n(\mu)-n(\mu')}{\rm Spec}_{\nu}(s_\lambda)$.

\begin{lemma} There exists a constant $W$ such that
 $$
  |{\rm Spec}_{\nu}(h_k)|\le W^k q^{k^2/2}
 $$
 for any $k$. Furthermore, for any $V>0$ there exists $k_0$ such that for every $k>k_0$ we have
$$
 |{\rm Spec}_{\nu}(h_k)|\le V^k q^{k^2/2}.
$$
\label{lemma_Spectialization_estimate}
\end{lemma}
\begin{proof}
 Recall that
$$
 \sum_{k=0}^{\infty}{\rm Spec}_{\nu}(h_k)t^k=H^{\nu}(t)=\prod_{x\in X(\nu)} (1-q^xt).
$$
(See Section \ref{subsection_computation_of_the_limits} for the definition of $X(\nu)$.) If
$X(\nu)$ is finite, then the statement of the lemma is obvious. Otherwise,
$$
 {\rm Spec}_{\nu}(h_k)=(-1)^k\sum_{1\le i_1<i_2<\dots<i_k} q^{X(\nu,i_1)+\dots+X(\nu,i_k)},
$$
where $X(\nu,i)$ stays for the $i$th element in $X(\nu)$. Clearly, $i_\ell\ge \ell$, consequently
$$
 \left|{\rm Spec}_{\nu}(h_k)\right|\le \sum_{i_1\ge 1,\, i_2\ge 2, \dots,i_k\ge k}
 q^{X(\nu,i_1)+\dots+X(\nu,i_k)}=q^{X(\nu,1)+\dots+X(\nu,k)}(1-q)^{-k}.
$$
Since $X(\nu,i)\ge i-1$, the first part of Lemma \ref{lemma_Spectialization_estimate} is proved.
Now recall that the complement of the set $X(\nu)$ in $\mathbb Z_{\ge 0}$ is infinite. Thus, for
any $r$ there exists $k_0$ such that for every $k>k_0$,
$$
X(\nu,1)+\dots+X(\nu,k)>(0+1+\dots+k-1)+rk.
$$
Choosing $r$ large enough we prove the second part of Lemma \ref{lemma_Spectialization_estimate}
\end{proof}
The Jacobi--Trudy formula for Schur polynomials (see e.g.\ \cite[3.4]{Mac}) implies that for
$\lambda\in{\mathbb{GT}}^+_k$, we have
$$
 {\rm Spec}_{\nu}(s_\lambda)={\rm det}\left[{\rm Spec}_{\nu}(h_{\lambda_i-i+j})\right]_{i,j=1,\dots,k},
$$
where we agree that ${\rm Spec}_{\nu}(h_{-1})={\rm Spec}_{\nu}(h_{-2})=\dots=0$. Since $k$ is
fixed, writing determinant as an alternating sum of products and applying Lemma
\ref{lemma_Spectialization_estimate} we conclude that for any $V>0$ there exists $n_0$ such that
for every $\lambda\in{\mathbb {GT}}^+_k$ with $\lambda_1>n_0$ we have
$$
 |{\rm Spec}_{\nu}(s_\lambda)|\le V^{\lambda_1} \prod_{i=1}^k q^{\lambda_i^2/2}.
$$
For any $\lambda\in{\mathbb{GT}}^+_k$ a simple computation proves that $n(\lambda')=\sum_{i=1}^k
\lambda_i(\lambda_i-1)/2$. Thus for $x_1,\dots,x_k$ with $|x_i|<M$ we have
\begin{multline*}
 \left|\sum_{\mu\in\mathbb{GT}_k^+,\, \mu_1\ge n_0}(-1)^{|\mu|}q^{n(\mu)-n(\mu')}{\rm
 Spec}_{\nu}(s_\mu){s^*_{\mu}(x_1,\dots,x_k)}\right|\\\le  \sum_{\mu\in\mathbb{GT}_k^+,\,
 \mu_1\ge n_0} V^{\mu_1} |s^*_{\mu}(x_1,\dots,x_k)|\le \sum_{\mu_1=n+0}^{\infty} V^{\mu_1}
 \left((\mu_1+1)^{k-1}\right) \left(C(M)^k\right)^{\mu_1}.
\end{multline*}
(Here we used the rough estimate that the number of signatures $\mu\in {\mathbb GT}_k$ with
$\mu_1=n$ is at most $(n+1)^{k-1}$.) If $V$ is small enough, then the above series converges.
Consequently,
$$ \sum_{\mu\in\mathbb{GT}_k^+}(-1)^{|\mu|}q^{n(\mu)-n(\mu')}{\rm
Spec}_{\nu}(s_\lambda){s^*_{\mu}(x_1,\dots,x_k)}$$ converges and defines an entire function.

\section{$q$--Toeplitz matrices (Proof of Proposition \ref{Proposition_q_toeplitz_non_neg})}


In this section we prove the following theorem:
\begin{theorem}
\label{theorem_q_toeplitz} Let $c_l$, $l\ge0$ be a sequence of non-negative numbers and suppose
that $\sum_l c_l q^{-l(l-1)}=1$. Let $d_{i,j}$, $i>0$, $j>0$ be a unique $q$--Toeplitz matrix such
that
 $$
  d[i,1]={c_{i-1}}, \quad i=1,2,\dots.
 $$
 Set
 $$
 \phi(t)=\sum_{\ell=0}^{\infty} c_\ell \prod_{i=0}^{\ell-1}(q^{-i}-t)
 $$
 an define coefficients $c_\lambda$ ($\lambda\in\mathbb{GT}^+_N$) as the coefficients of the
 expansion
\begin{equation}
\label{eq_x10}
 \phi(x_1)\cdots\phi(x_N)=\sum_{\lambda\in\mathbb{GT}^+_N} c_\lambda
 (-1)^{|\lambda|}{s^*_\lambda(q^{N-1}x_1,\dots,q^{N-1}x_N;q^{-1})}.
\end{equation}
Then
$$
 c_\lambda=q^{(N-1)|\lambda|} \det\biggl[d[\lambda_{N-i+1}+i,j] \biggr]_{i,j=1,\dots,N}.
$$
\end{theorem}
{\bf Remark.} If $\phi(t)$ is a polynomial, then the series \eqref{eq_x10} is finite. For general
$\phi(t)$ Proposition \ref{proposition_how_to_obtain_interpolation_coefficients} yields that the
series \eqref{eq_x10} converges at least in points $q^{1-N-\lambda+\delta}$, and such convergence
is enough for our purposes.

Before proving Theorem \ref{theorem_q_toeplitz} let us obtain Proposition
\ref{Proposition_q_toeplitz_non_neg} as its corollary.

\begin{proof}[Proof of Proposition \ref{Proposition_q_toeplitz_non_neg}]
Every initial column minor of the matrix $d^{\nu}[i,j]$ is
$$
 \det\biggl[d[\lambda_{N-i+1}+i,j] \biggr]_{i,j=1,\dots,N}
$$
for a certain $N$ and $\lambda\in\mathbb{GT}_N^+$. Theorem \ref{theorem_q_toeplitz} and Theorem
\ref{theorem_Main} imply that
\begin{multline*}
 \det\biggl[d_{\lambda_{N-i+1}+i,j} \biggr]_{i,j=1,\dots,N}=q^{-(N-1)|\lambda|}c_\lambda=
q^{-(N-1)|\lambda|}(-1)^{|\lambda} \frac{\mathcal
E^\nu_N(\lambda)}{s^*_{\lambda}(0,\dots,0;q^{-1})}.
\end{multline*}
Observe that
$$
 \frac{(-1)^{|\lambda}}{s^*_{\lambda}(0,\dots,0;q^{-1})}>0,
$$
thus, all initial column minors of $d^{\nu}[i,j]$ are non-negative.

To finish the proof note that since $d^{\nu}[i,j]$ is a triangular matrix all its initial row
minors (i.e.\ \,minors corresponding to the first $N$ rows and arbitrary $N$ columns) are either
zero or equal to certain initial column minors.
\end{proof}

To prove Theorem \ref{theorem_q_toeplitz} we need two general formulas involving factorial Schur
polynomials.

\begin{proposition} We have
\label{Proposition_dual_Cauchy}
$$
\prod_{i=1}^{N}\prod_{j=1}^m (y_j-x_i)=\sum_{\lambda\subset m^N} (-1)^{|\lambda |}
s_\lambda(x_1,\dots,x_N\mid a) s_{\hat \lambda'}(y_1,\dots,y_n\mid a),
$$
where $\lambda'$ is the transpose diagram and $\widehat\lambda'$ is the complement of $\lambda'$ in
diagram $N^m$
\end{proposition}
\begin{proof}
See \cite[Proof of (6.17)]{Mac2}.
\end{proof}

Denote
$$e_k(y_1,\dots,y_m\mid a)=\begin{cases}
  s_{1^k}(y\mid a),\quad 0<k\le m,\\
  1, \quad k=0, \\
  0,\quad k<0 \text{ or } k>m. \end{cases}
$$

\begin{proposition}
\label{Proposition_determinantal_formula_for_factorial_Schur}
Factorial Schur polynomials admit a determinantal formula:
$$
 s_{\lambda}(y_1,\dots,y_m\mid a)=\det[e_{\lambda'_i-i+j}(y_1,\dots,y_m\mid
\tau^{j-1}a)]_{i,j=1,\dots,m},
$$
where $m$ is an arbitrary integer greater or equal than $\lambda_1$, $\lambda'_i$ are row lengths
of the transposed diagram (equivalently, they are column lengths of $\lambda$), and $\tau^{j-1}a$
stands for the sequence $(\tau^{j-1}a)_i=a_{j-1+i}$.
\end{proposition}
\begin{proof}
See \cite[(6.7)]{Mac2}.
\end{proof}

Now let $\widehat a$ be the following sequence:
$$
 \widehat a_j=-q^{1-j},\quad j\in\mathbb{Z}.
$$
Now fix $m$ arbitrary numbers $y_1,\dots,y_m$ and let
$d[i,j]_{i,j=1,2,\dots}$ be an infinite matrix given by
$$
d[i,j]=e_{m-i+j}(y_1,\dots,y_m\mid \tau^{j-1}\widehat a).
$$
\begin{lemma}
\label{lemma_is_a_qToeplitz}
 $d[i,j]$ is a $q$--Toeplitz matrix, i.e.\
$$
 d[i,j+1]=d[i-1,j]+(q^{1-j}-q^{1-i})d[i,j].
$$
(Here we agree that $d[i,j]=0$ is either $i<1$ or $j<1$.)
\end{lemma}
\begin{proof}
Setting $N=1$ in Proposition \ref{Proposition_dual_Cauchy}, we get
$$
 \prod_{i=1}^m (y_i-t)= \sum_{k=0}^m (-1)^{k} h_k(t\mid a)
 e_{m-k}(y_1,\dots,y_m\mid a),
$$
where
$$
 h_k(t\mid a)=s_{k}(t\mid a)=(t+a_1)(t+a_2)\dots (t+a_k).
$$

We see that the numbers $e_{m-k}(y\mid a)$ are the coefficients of
the decomposition of the function $\prod_{i=1}^m (y_i-t)$ into the
sum of polynomials $(-1)^k h_k(t\mid a)$. We can use any sequence
$a$. In particular, if we set $a_i=-q^{2-i-\ell}$ then we get
\begin{equation}
\label{eq_defining_matrix} \prod_{i=1}^m (y_i-t)=\sum_{k=0}^m
d[k+\ell,\ell](q^{1-\ell}-t)(q^{1-\ell-1}-t)\dots(q^{2-k-\ell}-t).
\end{equation}

Observe that the left side of \eqref{eq_defining_matrix} does not
depend on $\ell$. Comparing \eqref{eq_defining_matrix} for $\ell$
and $\ell+1$ and using the fact that
\begin{multline*}
 (q^{1-\ell}-t)(q^{1-\ell-1}-t)\dots(q^{2-k-\ell}-t)\\=(q^{1-(\ell+1)}-t)(q^{1-(\ell+1)-1}-t)\dots(q^{2-k-(\ell+1)}-t)\\+
 (q^{1-\ell}-q^{2-k-(\ell+1)} ) (q^{1-(\ell+1)}-t)(q^{1-(\ell+1)-1}-t)\dots(q^{2-k-\ell}-t),
\end{multline*}
we get
$$
d[k+\ell+1,\ell+1]=d[k+\ell,\ell]+(q^{1-\ell}-q^{-k-\ell})d[k+\ell+1,\ell].
$$
To finish the proof substitute $j=\ell$, $i=k+\ell+1$.
\end{proof}

\begin{proposition}
\label{Proposition_q_Toeplitz_for_polynomials} Let $H(t)$ be a
polynomial of degree $m$ with $H(0)=1$. Let the coefficients
$c_\ell$, $i=0,\dots,m$ be defined from the expansion
 $$
 H(t)=\sum_{\ell=0}^{m} c_\ell \prod_{i=0}^{\ell-1}(q^{-i}-t).
 $$
 More generally, define coefficients $c_\lambda$
 ($\lambda\in\mathbb{GT}^+_N$) as the coefficients of the expansion
\begin{equation}
 \label{eq_x5} H(x_1)\cdots H(x_N)=\sum_{\lambda\in\mathbb{GT}^+}
 c_\lambda
 (-1)^{|\lambda|}{s^*_\lambda(q^{N-1}x_1,\dots,q^{N-1}x_N;q^{-1})}.
\end{equation}
 Let $\widetilde d[i,j]$ denote a unique $q$--Toeplitz matrix such that
$$
 \widetilde d[i,1]=c_{i-1}.
$$

Then we have
$$
 c_\lambda=q^{(N-1)|\lambda|} \det\biggl[\widetilde d[\lambda_{N-i+1}+i,j] \biggr]_{i,j=1,\dots,N}.
$$
\end{proposition}
\begin{proof}
 Let $y_i$ be the roots of $H(t)$, i.e.\
$$
 H(t)=\prod_{i=1}^{m}(1-y_i^{-1}t).
$$
Let, as above, define
$$
 \widehat a_j=-q^{1-j},\quad j\in\mathbb{Z}
$$
and
$$
d[i,j]=e_{m-i+j}(y_1,\dots,y_m\mid \tau^{j-1}\widehat a).
$$
Setting $N=1$ in Proposition \ref{Proposition_dual_Cauchy}, we
conclude that
$$
 d[i,1]= c_{i-1}\prod_{\ell=1}^m y_\ell^{-1}.
$$
 Lemma \ref{lemma_is_a_qToeplitz} yields that $d[i,j]$ is a
$q$--Toeplitz matrix, thus,
$$
 d[i,j]=  \widetilde d[i,j] \prod_{\ell=1}^m y_\ell^{-1}.
$$

Recall that
$$
 s_\lambda^*(x_1,\dots,x_N;q^{-1})=s_\lambda(x_1,\dots,x_N\mid a)
$$
 with
$$
 a_j=-q^{N-j}.
$$
 Definition of factorial Schur functions $s_\lambda(x\mid a)$
implies that
$$
 s_{\lambda}(qy\mid qa)=q^{|\lambda|}s_\lambda(y\mid a).
$$
(Here $qa$ stands for the sequence with $(qa)_j=q \cdot a_j$.)

It follows that
$$
 s_\lambda^*(q^{N-1}x_1,\dots,q^{N-1}x_N;q^{-1})=q^{(N-1)|\lambda|} s_\lambda(x_1,\dots,x_N\mid
\widehat a).
$$
Now let $b_\lambda$ be the coefficient of the expansion
$$
\prod_{i=1}^{N}\prod_{j=1}^m (y_j-x_i)=\sum_{\lambda\in m^N} (-1)^{|\lambda |}
s_\lambda(x_1,\dots,x_N\mid \widehat a) b_\lambda.
$$

Comparing the last formula with the definition \eqref{eq_x5} of
$c_\lambda$ we see that
$$
 c_\lambda=q^{(N-1)|\lambda|} \prod_{\ell=1}^m y_\ell^{-N}
b_\lambda.
$$

On the other hand, Proposition \ref{Proposition_dual_Cauchy} yields
that
$$
b_\lambda=s_{\widehat \lambda'}(y_1,\dots,y_m\mid \widehat a).
$$
Applying Proposition
\ref{Proposition_determinantal_formula_for_factorial_Schur} we
conclude that
$$
 b_\lambda=\det[e_{\widehat \lambda_i-i+j}(y_1,\dots,y_m\mid \tau^{j-1}\widehat
a)]_{i,j=1,\dots,N}.
$$
It is clear that
$$
 \widehat \lambda_i-i=m-(\lambda_{N-i}+i).
$$
Thus,
$$
 e_{\widehat \lambda_i-i+j}(y_1,\dots,y_m\mid \tau^{j-1}\widehat a)=d[\lambda_{N-i}+i,j].
$$
Consequently,
\begin{multline*}
 c_\lambda=q^{(N-1)|\lambda|} \prod_{\ell=1}^m y_\ell^{-N}
\det\biggl[d[\lambda_{N-i}+i,j]\biggr]_{i,j=1,\dots,N}\\=q^{(N-1)|\lambda|} \det\{\widetilde
d[\lambda_{N-i}+i,j]\}_{i,j=1,\dots,N}.
\end{multline*}
\end{proof}

Now we  are ready to prove Theorem \ref{theorem_q_toeplitz}.
\begin{proof}[Proof of Theorem \ref{theorem_q_toeplitz}]
  Recall that $c_l$, $l\ge 0 $ is a sequence of non-negative numbers and  $\sum_l c_l q^{-l(l-1)/2}=1$.

 Let us denote
 $$
  c^{(m)}_l=\begin{cases} \dfrac{c_l}{\sum_{i=0}^{m} c_i q^{-i(i-1)/2}}, \quad
0\le
l\le k,\\
0,\text{ otherwise.}
\end{cases}
 $$
 And let $d^{(m)}[i,j]$, $i>0$, $j>0$ be a unique $q$--Toeplitz
 matrix such that
 $$
  d^{(m)}[i,1]={c^{(m)}_{i-1}}, \quad i=1,2,\dots.
 $$
 Observe that
 $$
  \lim_{m\to\infty} d^{(m)}[i,j]=d[i,j].
 $$
 Consequently, all minors of $d[i,j]$ are limits of the
corresponding minors of $d^{(m)}[i,j]$.

Next, let $H^{(m)}(t)$ be a degree $m$ polynomial such that
$$
 H^{(m)}(t)=\sum_{\ell=0}^{m} c_\ell
\prod_{i=0}^{\ell-1}(q^{-i}-t).
$$

Let $c^{(k)}_\lambda$ be the coefficients of the expansion
$$
 H^{(m)}(x_1)\cdots H^{(m)}(x_N)=\sum_{\lambda\in\mathbb{GT}^+}
 c^{(k)}_\lambda
 (-1)^{|\lambda}{s^*_\lambda(q^{N-1}x_1,\dots,q^{N-1}x_N;q^{-1})}.
$$

Observe that
$$
 H^{(m)}(t)\rightrightarrows \phi(t)
$$
uniformly on compact subsets of $\mathbb C$. Hence,
$$
 H^{(m)}(x_1)\cdots H^{(m)}(x_N)\rightrightarrows \phi(x_1)\cdots
\phi(x_N).
$$
Then Proposition
\ref{proposition_how_to_obtain_interpolation_coefficients} implies
that
$$
 c^{(m)}_\lambda\to c_\lambda.
$$

Applying Proposition \ref{Proposition_q_Toeplitz_for_polynomials} we
conclude that
\begin{multline*}
 c_\lambda=\lim_{m\to\infty} c^{(m)}_\lambda=q^{(N-1)|\lambda|}\lim_{m\to\infty}
\det\biggl[d^{(m)}[\lambda_{N-i}+i,j]\biggr]_{i,j=1,\dots,N}\\=q^{(N-1)|\lambda|}
\det\biggl[d[\lambda_{N-i}+i,j]\biggr]_{i,j=1,\dots,N}.
\end{multline*}

\end{proof}

\end{document}